\documentclass[11 pt,a4paper,leqno,notitlepage]{report}
\usepackage{latexsym,amsmath}
\usepackage[francais]{babel}

\title{\textbf{Syst\`emes aux $q$-diff\'erences singuliers r\'eguliers:
solutions canoniques, classification, matrice de connexion et monodromie}}

\author{\textbf{Jacques SAULOY}\\
{\small Laboratoire Emile Picard, UMR 5580}\\
{\small Universit\'e Paul Sabatier,118,route de Narbonne,31062,Toulouse Cedex,
France.}\\
{\small E-mail : sauloy@picard.ups-tlse.fr}\\
}

\date{22 f\'evrier 1999}

\newcommand{\Div}{\text{div}}

\newcommand{\Diag}{\text{diag}}

\begin{document}

\maketitle

\hrule

{\small
\begin{center} \textbf{R\'esum\'e} \end{center}

G.D. Birkhoff a pos\'e, par analogie avec le cas classique des 
\'equations diff\'erentielles, le probl\`eme de Riemann-Hilbert pour les 
syst\`emes  ``fuchsiens'' aux $q$-diff\'erences lin\'eaires, \`a coefficients 
rationnels. Il  l'a r\'esolu dans le cas g\'en\'erique: l'objet classifiant 
qu'il introduit est
constitu\'e de la matrice de connexion $P$ et des exposants en $0$ et $\infty$.
Nous reprenons sa m\'ethode dans le cas g\'en\'eral, mais en traitant 
sym\'etriquement $0$ et  $\infty$ et sans recours \`a des solutions \`a 
croissance ``sauvage''. Lorsque $q$ tend vers $1$, $P$ tend vers une matrice 
localement constante $\tilde{P}$ telle que les valeurs (en nombre fini) 
$\tilde{P}(a)^{-1}\tilde{P}(b)$ 
sont les  matrices de monodromie du syst\`eme diff\'erentiel limite (suppos\'e 
non r\'esonnant en $0$ et $\infty$) en les singularit\'es de 
$\mathbf{C}^{*}$.
{\footnotesize Ce texte est celui du preprint 148 du Laboratoire
Emile Picard (f\'evrier 1999). Une version un peu plus courte en est parue
aux \emph{Annales de l'Institut Fourier}, \textbf{50}, 4, (2000).} \\

\textbf{Mots cl\'es:} Equations aux $q$-diff\'erences - Matrice de
connexion - Monodromie - Fonctions hyperg\'eom\'etriques basiques.
\textbf{Classification AMS:} 05A30 - 33D - 39A10 - 58 F.}

\bigskip \hrule 

{\small 
\begin{center} \textbf{Abstract} \end{center}

G.D. Birkhoff extended the classical Riemann-Hilbert problem for 
differential equations to the case of ``fuchsian'' 
linear $q$-difference systems with rational coefficients. He solved it in the 
generic case: the classifying object which he introduces is made up of the 
connection matrix $P$, 
together with the exponents at $0$ and $\infty$. We follow his method in the 
general case, but treat symetrically $0$ and $\infty$ and use no ``wildly'' 
growing solutions. When  $q$ tends to $1$, $P$ tends to a locally constant 
matrix $\tilde{P}$ such that the (finitely many) values 
$\tilde{P}(a)^{-1}\tilde{P}(b)$ 
are the monodromy  matrices of the limiting differential system (assumed to be 
non resonant at $0$ and $\infty$) at the singularities on $\mathbf{C}^{*}$.
{\footnotesize This text is that of preprint 148 of the  Laboratoire
Emile Picard (february 1999). A shorter version was published by
the \emph{Annales de l'Institut Fourier}, \textbf{50}, 4, (2000).} \\

\textbf{Keywords:} $q$-difference equations - Connection matrix -
Monodromy - Basic hypergeometric functions.
\textbf{AMS Classification:} 05A30 - 33D - 39A10 - 58 F.}

\bigskip \hrule

\setcounter{chapter}{-1}

\chapter{Introduction}



\begin{quote}

``Je suis convaincu que, tout comme pour les fonctions sp\'eciales solutions
d'\'equations diff\'erentielles, les formules int\'eressantes d\'erivent de
consid\'erations ``g\'eom\'etriques'' simples'' 
(Jean-Pierre Ramis,\cite{RamisJPRTraum}).
\end{quote}

\section{Histoire}

\subsection{Les racines et les branches}

La th\'eorie des \'equations aux $q$-diff\'erences est un sujet romantique,
puisqu'il a connu l'intervention d'Euler, Jacobi, Ramanujan ... et, \`a
un moindre degr\'e, de Gauss et de Cauchy. Les identit\'es d'Euler sur les 
partitions, la th\'eorie des fonctions th\^eta de Jacobi, avec sa magnifique 
formule du triple produit, de nombreuses formules myst\'erieuses de Ramanujan 
y sont li\'ees. On peut trouver une d\'erivation dans cet esprit de nombreuses 
identit\'es d'Euler dans \cite{GR}, et un lien direct entre la formule de 
Rogers-Ramanujan et une \'equation aux $q$-diff\'erences dans \cite{Andrews}. 
Ces deux livres contiennent aussi de nombreuses indications historiques.\\

C'est aussi un sujet en cours de renaissance. Ramis, dans \cite{RamisGrowth}
reprend les probl\`emes de comportement asymptotique des solutions \'etudi\'es
par Poincar\'e, Picard, Valiron et les replace dans le cadre de la th\'eorie
des classes Gevrey; ces recherches ont \'et\'e poursuivies par Zhang dans la
direction d'une th\'eorie de la $q$-sommabilit\'e (\emph{voir} \cite{Zhang}).
Etingof, van der Put, Singer ont \'etendu la th\'eorie de Picard-Vessiot avec 
pour but une th\'eorie de Galois (\emph{voir} \cite{Etingof}, \cite{SVdP}; 
ce dernier ouvrage traite, plus g\'en\'eralement, des \'equations aux 
diff\'erences). En tant que variante multiplicative de la m\'ethode d'Euler 
pour approcher les solutions d'\'equations diff\'erentielles, elles ont 
\'et\'e utilis\'ees dans l'\'etude des int\'egrales it\'er\'ees (\emph{voir} 
\cite{Agreg},\cite{Cartier1}). Il y a une th\'eorie g\'eom\'etrique 
(``holonome'') des identit\'es combinatoires, y compris $q$-combinatoires
et une nouvelle approche de leur calcul formel (Zeilberger, Aomoto, Sabbah,
Chyzak: \emph{voir} \cite{Cartier2}, \cite{Chyzak}). La th\'eorie des groupes 
quantiques commence aussi \`a essaimer dans cette direction (\emph{voir} 
\cite{LesRusses}). Enfin, la tradition qui, depuis Birkhoff et Ore veut 
traiter \'equations diff\'erentielles, aux diff\'erences et aux 
$q$-diff\'erences dans un cadre unifi\'e a \'et\'e r\'ecemment reprise par 
Yves Andr\'e qui a cr\'e\'e pour cela le cadre des connexions non 
commutatives (\emph{voir} \cite{YA}).

\subsection{$q$-analogies, classification et confluence}

La s\'erie hyperg\'eom\'etrique \emph{classique} a \'et\'e g\'en\'eralis\'ee 
par Heine, puis Ramanujan (\emph{voir} \cite{GR}, \cite{Ramanujan}) en la 
s\'erie hyperg\'eom\'etrique \emph{basique}\footnote{On prendra garde que nous faisons, dans tout ce travail, 
l'hypoth\`ese $|q| > 1$. Il est, bien s\^ur, essentiel de supposer 
que $|q| \not= 1$, mais certains des ouvrages cit\'es supposent $|q| < 1$.}
$$\Phi(a,b,c;q,z) = \sum_{n \geq 0}
\frac{(a;p)_{n} (b;p)_{n}}{(c;p)_{n} (p;p)_{n}}z^{n}
\quad |q| > 1 \;,\; p = q^{-1} \;,\; 
(x;p)_{n} = \prod_{0 \leq k \leq n-1}(1 - x p^{k})$$
Si $a,b,c \in \mathbf{C}^{*}$, celle-ci est solution de l'\'equation aux 
$q$-diff\'erences \`a coefficients rationnels:
$$
\Phi(q^{2} z) - 
\frac{(a + b)z - (1+c/q)}{abz - c/q} \Phi(q z) + 
\frac{z-1}{abz-c/q} \Phi(z) = 0
$$
Nous \'ecrirons plut\^ot
$$
\sigma_{q}^{2}\Phi - 
\frac{(a + b)z - (1+c/q)}{abz - c/q} \sigma_{q}\Phi + 
\frac{z-1}{abz-c/q} \Phi = 0
$$
o\`u l'on note $\sigma_{q}f(z) = f(qz)$ pour toute
fonction sur la sph\`ere de Riemann $\mathbf{S} = \mathbf{P}^{1}\mathbf{C}$.
Ceci a motiv\'e la r\'esolution et la classification de telles \'equations 
fonctionnelles, sous le nom d'\'equations aux $q$-diff\'erences. (\emph{voir}
\cite{Adams},\cite{Birkhoff}, \cite{Carmichael}, \cite{GR}, \cite{RamisGrowth},
\cite{RamisJPRTraum}).\\ 

R\'eciproquement, on sait que, lorsque l'on fait tendre $q$ vers $1$, sous des
conditions appropri\'ees portant sur les param\`etres, $\Phi(a,b,c;q,z)$ 
converge vers la s\'erie hyperg\'eom\'etrique: c'est au moins
\'evident terme \`a terme, et cela reste vrai dans un sens moins formel
(voir le traitement complet de cet exemple en 4.4.2). De telles $q$-analogies
sont pr\'esentes \`a un degr\'e plus ou moins important dans la plupart des
exemples que nous avons cit\'es au d\'ebut de cette section. On peut 
l\'egitimement consid\'erer pas mal de ces ph\'enom\`enes comme des cas de
\emph{confluence}.\\ 

Par analogie avec le cas classique des \'equations diff\'erentielles,
G.D. Birkhoff a pos\'e le probl\`eme de Riemann-Hilbert pour les syst\`emes 
lin\'eaires aux $q$-diff\'erences \`a coefficients rationnels 
(\emph{voir} \cite{Birkhoff}):
\begin{equation}
\sigma_{q}X = AX   ,\hspace{3cm}   A \in GL_{n}(\mathbf{C}(z))
\end{equation}

Le groupe de jauge $GL_{n}(\mathbf{C}(z))$ op\`ere \`a gauche sur les solutions
vectorielles ou matricielles de telles \'equations et donc sur les \'equations
elles-m\^emes, ce qui conduit \`a introduire la notion d'\emph{\'equivalence
m\'eromorphe} : $A \sim (\sigma_{q}F)^{-1}AF$ lorsque 
$F \in GL_{n}(\mathbf{C}(z))$. Il s'agit alors de classifier ces syst\`emes 
modulo cette relation d'\'equivalence.\\

Notre but, dans ce travail, est de tenter de clarifier la g\'eom\'etrie
sous-jacente. Dans ce premier rapport, nous \'etudions la 
\emph{r\'esolution}, la \emph{classification} et la \emph{confluence} 
d'une classe d'\'equations aux $q$-diff\'erences, les \'equations 
\emph{fuchsiennes} (ainsi que la classe, qui n'est plus large qu'en apparence, 
des \'equations \emph{singuli\`eres r\'eguli\`eres}). Dans un rapport 
ult\'erieur (\cite{JS2}), nous commencerons l'\'edification d'une 
\emph{th\'eorie de Galois} de ces \'equations sur ces bases, th\'eorie
qui devrait
donc diff\`erer sensiblement de celle,plus alg\'ebrique et plus
g\'en\'erale, de \cite{SVdP} (voir ci-dessous 0.2.3 et 0.5).


\section{R\'esolution}


\subsection{Cas g\'en\'eral}

Le cas o\`u $A(0)$ est singuli\`ere est illustr\'e par l'\'equation
$\sigma_{q}f = zf$ dont une solution typique est la fonction $\Theta$: ce cas
se caract\'erise par l'apparition de solutions \`a croissance sauvage 
(\emph{voir} \cite{RamisGrowth}). A l'oppos\'e, si $A(0) = I_{n}$, le produit
$\underset{i \geq 1}{\prod} A(q^{-i}z)$ d\'efinit une solution holomorphe en 
$0$ et m\'eromorphe sur $\mathbf{C}$. Le premier cas est consid\'er\'e comme
\emph{irr\'egulier}, le second cas comme \emph{r\'egulier}.\\

Le cas interm\'ediaire est le cas singulier r\'egulier, que nous \'etudierons 
plus particuli\`erement: c'est
celui o\`u l'on peut se ramener, via une transformation de jauge rationnelle,
\`a supposer l'\'equation fuchsienne, c'est \`a dire telle que $A(0)$ et
$A(\infty)$ sont inversibles. Les exemples les plus \'el\'ementaires sont: 
$\sigma_{q}f = c f , c \in \mathbf{C}^{*}$ et $\sigma_{q}f = f + 1$. Cette 
derni\`ere se ram\`ene en effet \`a une \'equation lin\'eaire d'ordre $2$,
\`a matrice unipotente. Les solutions de telles \'equations sont les briques 
de base pour la r\'esolution de toutes les \'equations fuchsiennes.


\subsection{Cas fuchsien}

On l'appelle fuchsien parce que la m\'ethode de Frobenius marche et permet 
d'exhiber des solutions de forme tr\`es simple et \`a croissance mod\'er\'ee. 
Comme dans le cas fuchsien classique (des \'equations diff\'erentielles),
nous pourrons m\^eme d\'emontrer la r\'eciproque. On peut cependant se rendre
compte en \'etudiant la confluence qu'il y a d\'ej\`a dans ce cas des aspects
qui le rapprochent du cas irr\'egulier des \'equations diff\'erentielles.\\

Pour commencer par l'\'etude locale, supposons $A(0) \in GL_{n}(\mathbf{C})$ : 
un tel syst\`eme est dit \emph{fuchsien en $0$}. En parfaite analogie
avec la m\'ethode de Frobenius pour les \'equations diff\'erentielles 
lin\'eaires, on construit alors une solution matrice fondamentale
$X \in GL_{n}(\mathbf{K'_{0}})$ \`a coefficients dans une extension
$\mathbf{K'_{0}}$ de $\mathbf{C}(\{z\})$ qui est un sous-corps de
$\mathcal{M}(\mathbf{C}^{*})$ et qui contient :

\begin{itemize}

\item{Pour chaque $c \in Sp(A(0))$, un \'el\'ement $e_{q,c}$ tel que
$\sigma_{q}e_{q,c} = c e_{q,c}$. Il joue un r\^ole analogue \`a celui de
$z^{\gamma}$ et nous le nommerons \emph{caract\`ere d'exposant $c$}.}

\item{Dans le cas o\`u $A(0)$ n'est pas semi-simple, un \'el\'ement $l_{q}$ 
tel que
$\sigma_{q}l_{q} = l_{q} + 1$. Il joue un r\^ole analogue \`a celui du 
logarithme.}

\end{itemize}

L'\'equation fonctionnelle (1), r\'e\'ecrite\footnote{
Le Petit Larousse Illustr\'e donne comme \'egalement correctes les formes
\emph{r\'e\'ecrire} et \emph{r\'ecrire}. On trouvera, dans la Gazette des
Math\'ematiciens de janvier 1999, une argumentation en faveur de ce dernier.
Mais au fond, c'est en fonction du plaisir de l'oeil que l'on tranche ...} 
sous la forme
$X(z) = A(z/q)X(z/q)$, entraine alors que les solutions sont en fait 
\emph{\`a coefficients dans le sous-corps
$\mathbf{K_{0}} = \mathcal{M}(\mathbf{C})(l_{q},(e_{q,c})_{c \in \mathbf{C}^{*}})$
de $\mathbf{K'_{0}}$}. En opposition \`a ce qui se passe avec les \'equations
diff\'erentielles, l'op\'erateur aux $q$-diff\'erences $\sigma_{q}$ propage 
dans $\mathbf{C}$ la propri\'et\'e des solutions d'\^etre m\'eromorphes en 
$0$. Ce n'est
\'evidemment pas le cas des automorphismes infinit\'esimaux de la th\'eorie des
\'equations diff\'erentielles et c'est l\`a une diff\'erence fondamentale.\\

Ceci se voit d'ailleurs d\'ej\`a dans le cas r\'egulier, qui, pour nous, est
caract\'eris\'e par le fait que  la structure de Jordan en $0$ est 
triviale, autrement dit, $A(0) = I_{n}$. L'\'equation peut alors
\^etre r\'esolue sans sortir de $\mathcal{M}(\mathbf{C})$ : 
notant $\Sigma'$ le lieu 
polaire de $A$, il est ais\'e de voir que le produit infini 
$A(z/q) A(z/q^{2}) ...$ converge uniform\'ement sur tout compact de 
$\mathbf{C} - q^{\mathbf{N}^{*}}\Sigma'$ et a des p\^oles sur 
$q^{\mathbf{N}^{*}}\Sigma'$. C'est \'evidemment une solution fondamentale
(d\'eterminant non nul) de l'\'equation si on la restreint \`a 
$\mathbf{C} - q^{\mathbf{N}^{*}}\Sigma$, $\Sigma$ \'etant $\Sigma'$
augment\'e des z\'eros de $\det A$: c'est ce dernier ensemble qu'il faudra
consid\'erer comme le v\'eritable lieu singulier de l'\'equation.


\subsection{Le choix des briques de base}

Classiquement (\emph{voir} \cite{Adams},\cite{Carmichael}), on utilise 
respectivement les fonctions multivalu\'ees $z^{\frac{\log(c)}{\log(q)}}$ et 
$\frac{\log(z)}{\log(q)}$ dans les r\^oles de $e_{q,c}$ et de $l_{q}$. 
Il faut alors
prendre en compte conjointement les effets de l'automorphisme $\sigma_{q}$ et
ceux de la monodromie. Singer et Van der Put proposent dans \cite{SVdP} 
d'instancier les $e_{q,c}$ et $l_{q}$ par des \emph{symboles} alg\'ebriquement
ind\'ependants (sauf en ce qui concerne des relations oblig\'ees). Mais il est
alors difficile d'interpr\'eter les solutions comme des fonctions. Nous 
montrons au chapitre 1 comment r\'esoudre (1) \`a l'aide de fonctions 
\emph{uniformes}.\\

La question peut \^etre formul\'ee ainsi: quelle $q$-analogies suivre ? 
Il y en a tant ! L'analogie suivie ici (motiv\'ee, entre autres, par 
l'exemple des s\'eries hyperg\'eom\'etriques basiques) est 
$(\sigma_{q}-1)/(q-1) \leftrightarrow \delta = z (d/dz)$ 
(\emph{voir} \cite{RamisGrowth}).\\

D\'esignons par $\Theta_{q}$ 
la fonction theta de Jacobi: $\sum_{n \in Z}(-1)^{n}q^{-n(n-1)/2}z^{n}$ 
(\emph{voir} \cite{RamisGrowth}). Elle satisfait la formule du triple produit
de Jacobi:
$$
\Theta_{q}(z) = (p;p)_{\infty} (z;p)_{\infty} (pz^{-1};p)_{\infty}
\quad (z;p)_{\infty} = \prod_{i \geq 0} (1 - p^{i}z)
$$
On en
d\'eduit l'\'equation fonctionnelle $\Theta_{q}(qz) = -qz \Theta_{q}(z)$.
Nous prendrons alors
$$
l_{q}(z) = z\Theta_{q}'(z)/\Theta_{q}(z) \quad \text{    et    } \quad
e_{q,c}(z) = \Theta_{q}(z)/\Theta_{q}(c^{-1}z)
$$ 
qui sont \emph{m\'eromorphes sur $\mathbf{C}^{*}$}. 
Les p\^oles (simples) de $l_{q}$ sont les \'el\'ements de la spirale 
discr\`ete 
$q^{\mathbf{Z}}$; les z\'eros (simples) et p\^oles (simples) de $e_{q,c}$ sont 
les \'el\'ements de $q^{\mathbf{Z}}$ et de $cq^{\mathbf{Z}}$ respectivement. 
Nous avons donc remplac\'e la \emph{ramification} des choix classiques par des 
\emph{spirales discr\`etes de p\^oles}.
Ces p\^oles et leur g\'eom\'etrie vont en fait contr\^oler toute la situation, 
apportant un d\'ebut d'illustration de l'importance de 
``consid\'erations g\'eom\'etriques simples''.\\

Nous obtenons ainsi des solutions fondamentales canoniques 
m\'eromorphes sur $\mathbf{C}^{*}$ et en donnons des formes normales. Dans le 
cas o\`u le syst\`eme est \emph{non r\'esonnant en $0$}, c'est \`a dire o\`u
deux valeurs propres distinctes de $A(0)$ ne sont pas congrues modulo
$q^{\mathbf{Z}}$, il y a m\^eme une forme canonique. Nous \'etendons de plus
la r\'esolution au cas d'un syst\`eme \emph{singulier r\'egulier en $0$}, c'est
\`a dire, par d\'efinition, m\'eromorphiquement \'equivalent \`a un syst\`eme 
fuchsien en $0$. Dans tous les cas, on peut \'ecrire $X^{(0)} = M^{(0)}N^{(0)}$
o\`u:\\ 

(i) $M^{(0)}$ est m\'eromorphe sur $\mathbf{C}$; les p\^oles de $M^{(0)}$
sur $\mathbf{C}^{*}$ forment une union finie de demi-spirales logarithmiques 
discr\`etes: $q^{\mathbf{N}^{*}}\mathcal{S}(A)$, o\`u $\mathcal{S}(A)$ est
l'ensemble form\'e des p\^oles de $A$ et des p\^oles de $A^{-1}$.\\

(ii) $N^{(0)}$ est construite \`a partir des $e_{q,c}$ et de $l_{q}$;
les p\^oles de $N^{(0)}$ sur $\mathbf{C}^{*}$ forment une union finie de
spirales logarithmiques discr\`etes d\'etermin\'ees par la structure de Jordan
de $A(0)$: $q^{\mathbf{Z}}Sp(A(0))$ si $A(0)$ est semi-simple, 
$q^{\mathbf{Z}}Sp(A(0)) \cup q^{\mathbf{Z}}$ dans le cas g\'en\'eral.


\subsection{Le corps des constantes}

Toutes les solutions dans $GL_{n}(\mathbf{K}'_{0})$ s'obtiennent alors sous la
forme $X^{(0)}V^{(0)}$ o\`u $V^{(0)}$ est \`a coefficients dans le ``corps des
constantes'' $(\mathbf{K}'_{0})^{\sigma_{q}}$.\\

On voit facilement que (pour notre choix standard de caract\`eres)
celui ci est \'egal 
\`a $\mathcal{M}(\mathbf{C})^{\sigma_{q}}$, c'est \`a dire au corps des
fonctions m\'eromorphes sur $\mathbf{C}^{*}/q^{\mathbf{Z}}$. Le changement de
variable $z = e^{2 \imath \pi x}$ permet d'identifier le groupe
$\mathbf{C}^{*}/q^{\mathbf{Z}}$ au tore complexe 
$\mathbf{C}/(\mathbf{Z} + \mathbf{Z}\tau)$, avec la notation 
$q = e^{- 2 \imath \pi \tau} \;,\; Im(\tau) > 0$. Le corps
$\mathcal{M}(\mathbf{C})^{\sigma_{q}}$ s'identifie alors au corps des fonctions
elliptiques $\mathcal{M}(E)$.\\

Notons que la r\'ealisation des caract\`eres comme fonctions 
m\'eromorphes sur $\mathbf{C}^{*}$ entraine \emph{n\'ecessairement} la 
pr\'esence d'un ``cocycle'' non trivial de fonctions elliptiques, les 
$\frac{e_{q,cd}}{e_{q,c}e_{q,d}}$: c'est cette circonstance qui fait du 
corps des 
constantes le corps des fonctions elliptiques $\mathcal{M}(\mathbf{E})$ 
tout entier pour notre choix standard de carct\`eres et, dans tous les cas, 
une extension transcendante de $\mathbf{C}$. Ceci est la diff\'erence 
principale avec \cite{SVdP}, o\`u les caract\`eres sont des symboles et peuvent
\^etre astreints \`a v\'erifier : $e_{q,cd} = e_{q,c}e_{q,d}$, ce qui entraine 
que leur corps des constantes est $\mathbf{C}$. Cela a, entre autres, des 
cons\'equences pour la th\'eorie de Galois, puisqu'une approche classique de
celle-ci est la th\'eorie de Picard-Vessiot, qui exige que le corps des
constantes ne bouge pas.


\section{Classification}


\subsection{Relations d'\'equivalence}

Si l'on admet qu'une transformation de jauge rationnelle $X = F Y$,
$F \in GL_{n}(\mathbf{C}(z))$, n'affecte pas
essentiellement la nature analytique des solutions, on est naturellement
conduit \`a consid\'erer l'\'equation de matrice $A$ satisfaite par $X$
comme \'equivalente \`a l'\'equation de matrice $B$ satisfaite par $Y$:
$A \sim B = (\sigma_{q}F)^{-1}AF$. C'est la notion d'\'equivalence
m\'eromorphe sur $\mathbf{S}$, ou d'\'equivalence rationnelle.\\

Si l'on veut l'implication et non l'\'equivalence, apparait la version 
moins sym\'etrique: $(\sigma_{q}F) B = A F \;,\; F \in M_{n}(\mathbf{C}(z))$, 
qui donnera lieu dans 
\cite{JS2} \`a une notion plus g\'en\'erale de morphismes et \`a une 
version plus fonctorielle de la classification. En termes de logique, 
on peut comprendre cette \'egalit\'e comme  une ``weakest precondition'': 
quelle condition sur $F$ garantit que $FY$ est solution de (1) ? On en 
trouve une description plus intrins\`eque dans 
\cite{SVdP} grace \`a la notion de ``difference modules''.\\

Une autre extension de cette relation d'\'equivalence apparaitra lors des
\'etudes locales: si l'on impose seulement \`a $F$ d\^etre m\'eromorphe
sur $\mathbf{C}$ (ou, dualement, sur $\mathbf{S} - \{0\}$), on obtient encore
d'utiles
th\'eor\`emes de r\'eduction qui permettront de normaliser les solutions.


\subsection{Riemann-Hilbert-Birkhoff}

A un syst\`eme fuchsien en $0$ et en $\infty$, et \`a 
des solutions fondamentales $X^{(0)}$ et $X^{(\infty)}$ en $0$ et en $\infty$, 
G.D. Birkhoff associe la matrice de connexion 
$P =  (X^{(\infty)})^{-1}X^{(0)}$, 
dont les coefficients sont $\sigma_{q}$-invariants, multivalu\'es et peuvent 
\^etre 
vus comme des fonctions m\'eromorphes sur $\mathbf{C}$, automorphes pour le 
r\'eseau $(\mathbf{Z} + \mathbf{Z}\tau)$ (o\`u $\mathbf{Z}\tau$ agit 
trivialement et $\mathbf{Z}$ par monodromie). Il code le syst\`eme 
(1) \`a l'aide de $P$ et des exposants en $0$ et en $\infty$ et montre que l'on
peut r\'esoudre le probl\`eme inverse dans le cas o\`u $A(0)$ et $A(\infty)$
sont semi-simples.\\

Malgr\'e les hypoth\`eses fuchsiennes sym\'etriques en $0$ et $\infty$,
Birkhoff commence par dissym\'etriser le probl\`eme, puis fait intervenir des
fonctions \`a croissance sauvage, du type $q^{t} \mapsto q^{t(t-1)/2}$ 
(\emph{voir} \cite{RamisGrowth}). Selon \cite{SVdP}, la construction de la 
matrice de connexion pr\'esente des ambiguit\'es dans certains cas. Il est en
tout cas probable que la matrice de connexion de Birkhoff n'a pas de tr\`es
bonnes propri\'et\'es multiplicatives. Ces difficult\'es sont surmont\'ees par
van der Put et Singer, qui donnent une classification dans \cite{SVdP}. Mais
leurs briques de base sont des symboles formels, ils introduisent les fonctions
a posteriori et l'espace sur lequel celles-ci vivent est assez compliqu\'e.


\subsection{Comment classifier avec des fonctions uniformes}

Nous montrons au chapitre 2 que nos solutions m\'eromorphes conduisent \`a
une matrice de connexion \emph{elliptique}, c'est \`a dire \`a coefficients
dans $\mathcal{M}(\mathbf{E})$. Notre m\'ethode fait jouer un r\^ole
sym\'etrique \`a $0$ et $\infty$ et \'evite le recours \`a des fonctions \`a 
croissance sauvage. 
Nous obtenons alors le th\'eor\`eme de classification des
syst\`emes singuliers r\'eguliers (par la matrice de connexion, assortie des
structures de Jordan en $0$ et $\infty$), ainsi que des pr\'ecisions sur la 
possibilit\'e de prescrire les p\^oles.


\section{Confluence}

Pour donner un sens g\'eom\'etrique \`a ces constructions, 
nous \'etudions la classique 
confluence d'\'equations aux $q$-diff\'erences vers des \'equations 
diff\'erentielles lorsque $q$ tend vers $1$.


\subsection{Confluence des solutions}

Nous nous restreignons pour cela au cas fuchsien non r\'esonnant.
Le syst\`eme admet alors une unique solution 
$X^{(0)} = M^{(0)}e_{A(0)}$,
o\`u $M^{(0)} \in GL_{n}(\mathcal{M}(\mathbf{C}))$, $M^{(0)}(0) = I_{n}$.
Sous les m\^emes hypoth\`eses en $\infty$, on introduit la solution
$X^{(\infty)} = M^{(\infty)}e_{A(\infty)}$ et la matrice de connexion
$P =  (X^{(\infty)})^{-1}X^{(0)}$.\\

Nous supposerons que \emph{$q \to 1$ le long d'une spirale logarithmique} :
$q = q_{0}^{\epsilon}$, c'est \`a dire 
$\tau = \tau_{0} \epsilon$, $q_{0} = e^{-2i\pi\tau_{0}}$, o\`u l'on a fix\'e
$\tau_{0}$ tel que $Im(\tau_{0}) > 0$ et o\`u $\epsilon \rightarrow 0^{+}$.
Les spirales et demi-spirales logarithmiques discr\`etes $q^{\mathbf{Z}}$,
$q^{\mathbf{N}}$ se condensent alors en des coupures spirales 
$q_{0}^{\mathbf{R}}$,
$q_{0}^{\mathbf{R}^{+}}$ de $\mathbf{S}$. 
Nous supposons de plus que $A = A_{\epsilon} = I_{n}+ (q-1)B_{\epsilon}$, o\`u 
$B_{\epsilon} \to \tilde{B}$. 
Sous des hypoth\`eses appropri\'ees, les solutions fondamentales canoniques de 
$(1)$ tendent alors vers celles de
\begin{equation}
\delta \tilde{X} = \tilde{B} \tilde{X}
\end{equation}
obtenues par la m\'ethode de Frobenius (chapitre 3). Bien que nos solutions
canoniques \`a (1) soient sans monodromie, ce n'est \'evidemment pas possible 
en g\'en\'eral pour les solutions de (2). En fait, les spirales et 
demi-spirales discr\`etes de
p\^oles des solutions de (1) se transforment, quand $\epsilon \to 0^{+}$,
en des \emph{coupures spirales} 
qui rendent les solutions de (2) \emph{uniformes}.


\subsection{Confluence de la matrice de connexion et monodromie}

Supposons maintenant les m\^emes conditions satisfaites en $\infty$, et notons
$X_{\epsilon}^{(0)}$, $X_{\epsilon}^{(\infty)}$ et $P_{\epsilon}$ les solutions
canoniques et la matrice de connexion associ\'ees. On montre alors, au 
chapitre 4, que la matrice $P_{\epsilon}$
tend vers une matrice \emph{localement constante} $\tilde{P}$; celle-ci est 
d\'efinie sur $\mathbf{S}$ priv\'e de toutes les coupures des solutions en $0$ 
et $\infty$, c'est \`a dire sur une union finie d'ouverts connexes: $\tilde{P}$
prend donc un nombre fini de valeurs $\tilde{P}_{0}$,...,$\tilde{P}_{r}$. En 
prolongeant les d\'eterminations des solutions en $0$ et en $\infty$ le long de
chemins qui \'evitent les diverses coupures, on conclut alors que les 
$\tilde{P}_{j}^{-1} \tilde{P}_{j-1} (1 \le j \le r)$ 
sont les matrices de monodromie de (2) en les singularit\'es autres que $0$
et $\infty$ (la monodromie en $0$ et $\infty$ est fournie par les structures 
de Jordan de $\tilde{B}(0)$ et $\tilde{B}(\infty)$).\\

Le r\'esultat obtenu \'evoque d'ailleurs de fa\c{c}on frappante la
description du groupe de Galois dans \cite{Etingof} (pour le cas r\'egulier)
et dans \cite{SVdP} (pour le cas singulier r\'egulier g\'en\'eral). 
Ces auteurs d\'ecrivent en effet le groupe de Galois
comme engendr\'e par les $P(a)^{-1}P(b)$, $P$ \'etant la matrice de
connexion. Nous retrouverons dans \cite{JS2} un \'enonc\'e similaire, 
\'equivalent \`a celui de \cite{Etingof} dans le cas r\'egulier mais
malheureusement beaucoup plus compliqu\'e que celui de \cite{SVdP} dans le
cas g\'en\'eral.


\subsection{Exemples}

Le jeu du franchissement des coupures, qui proviennent elles-m\^emes 
des barri\`eres de p\^oles, illustre \`a merveille le r\^ole pr\'epond\'erant
de la g\'eom\'etrie des p\^oles! Ceux-ci se constituent en barri\`eres en
forme de spirales logarithmiques discr\`etes, qui se condensent, lorsque
$q \to 1$, en des coupures qui donnent lieu \`a des automorphismes de
monodromie.\\ 

Nous esp\'erons de m\^eme comprendre les automorphismes galoisiens
\emph{avant confluence} (donc au niveau des \'equations aux $q$-diff\'erences)
comme des automorphismes de monodromie; si nous pouvons observer de telles
transformations, nous supposerons que l'on a tourn\'e autour d'un point, et
donc qu'il y avait un point !\\

Le ph\'enom\`ene d\'ecrit ci-dessus a d'abord \'et\'e observ\'e sur les trois
exemples d\'ecrits en 4.4. Il peut \^etre agr\'eable de les \'etudier en
parall\`ele avec la lecture de tout ce rapport. Ils sont plut\^ot une
illustration de la th\'eorie qu'une application, et sont trait\'es ``\`a la
main''. D'ailleurs, comme pour les \'equations diff\'erentielles, on utilise
dans la pratique des \'equations plut\^ot que des syst\`emes, on en
recherche les exposants \`a l'aide de l'\'equation caract\'eristique, etc...


\section{Perspectives galoisiennes}

L'\'etude pr\'esent\'ee ici est plut\^ot algorithmique et calculatoire; elle
fait pourtant d\'ej\`a apparaitre un important soubassement g\'eom\'etrique.
Les spirales de p\^oles; les diviseurs des caract\`eres et des fonctions 
th\^eta; le cocycle des caract\`eres: tout montre que la courbe elliptique 
$\mathbf{C}^{*}/q^{\mathbf{Z}}$ joue un r\^ole crucial.\\ 

La confluence peut \^etre vue comme une situation de d\'eg\'en\'erescence de
courbes elliptiques; on l'observe d\'ej\`a ici explicitement en ce qui
concerne la d\'eg\'en\'erescence des r\'eseaux dans $\mathbf{C}$.\\

Nous n'avons aucun espoir d'une th\'eorie de Picard-Vessiot: notre corps des 
constantes est bien trop gros. Pourtant, l'apparition de d\'efauts d'unicit\'e
dans la classification qui sont des matrices de $GL_{n}(\mathbf{C})$ nous
laisse esp\'erer que nous pourrons (comme l'ont fait van der Put et Singer
dans \cite{SVdP}) construire un groupe de Galois sur $\mathbf{C}$: nous
verrons dans \cite{JS2} que c'est bien le cas.\\

Enfin, le r\^ole central (pour la classification) de l'uniformisation de
la courbe elliptique $\mathbf{C}^{*}/q^{\mathbf{Z}}$ \`a l'aide de fonctions
th\^eta, indique la possibilit\'e d'\'etendre ces r\'esultats au cas
$p$-adique. Sur une suggestion de Yves Andr\'e, ceci sera fait autant que 
possible dans le cadre plus fonctoriel (en fait, tannakien) de \cite{JS2}.
On verra alors qu'il n'est pas indiff\'erent d'avoir travaill\'e sur le
mod\`ele ``multiplicatif'' $\mathbf{C}^{*}/q^{\mathbf{Z}}$ de la courbe
elliptique plut\^ot que sur son mod\`ele ``additif'' 
$\mathbf{C}/(\mathbf{Z} + \mathbf{Z}\tau)$.

\section*{Remerciements}

Ce travail fait partie d'une th\`ese sous la direction de Jean-Pierre Ramis.
La partie pr\'esent\'ee dans ce rapport a donn\'e lieu \`a la publication de
deux notes aux Comptes Rendus de l'Acad\'emie des Sciences (janvier 1999) et
\`a la soumission d'un article.\\

Je remercie chaleureusement Jean-Pierre pour son soutien et pour ses conseils, 
et surtout pour m'avoir fait partager sa vision g\'eom\'etrique de ce monde 
fascinant, ainsi que d'autres de ses r\^eves. J'esp\`ere avoir traduit un peu
de la beaut\'e de cette vision.\\

C'est \'egalement avec plaisir que je remercie Mich\`ele Loday-Richaud pour
une discussion sur la classification des \'equations diff\'erentielles
singuli\`eres r\'eguli\`eres qui n'a pas peu contribu\'e \`a me d\'esembrumer;
et Jean-Claude Sikorav pour des lumi\`eres sur le th\'eor\`eme de Birkhoff 
et sur d'autres ph\'enom\`enes complexes, et aussi pour \^etre \emph{toujours} 
pr\^et \`a parler de maths, \emph{quelles qu'elles soient}.

\newpage

\section{Notations g\'en\'erales}

\emph{On prendra garde \`a la convention $|q| > 1$ adopt\'ee ici et qui
varie d'une source \`a l'autre.}\\

On fixe\footnote{
Les notations li\'ees \`a la variation de $q$ sont sp\'ecifiques
de la deuxi\`eme partie et seront donn\'ees alors.}
un nombre complexe $\tau$ de partie imaginaire $> 0$, et l'on pose 
$q = e^{-2i\pi\tau}$.\\

\subsection*{Formules de $q$-combinatoire}

On notera $p = 1/q$ et 
$$(x;p)_{n} = \prod_{0 \leq k \leq n-1}(1 - x p^{k})$$
$$(x;p)_{\infty} = \prod_{0 \leq k}(1 - x p^{k})$$

\subsection*{Op\'eration de $q^{\mathbf{Z}}$}

On notera $\mathbf{S}$ la sph\`ere de Riemann, c'est \`a dire la droite 
projective complexe $\mathbf{P}^{1}\mathbf{C}$
et $\mathcal{H}$ le demi-plan de Poincar\'e
$\{z \in \mathbf{C} \;/\; Im(z) > 0\}$.\\

Le groupe cyclique $q^{\mathbf{Z}}$ agit sur $\mathbf{S}$ par homographies 
$z \mapsto q z$, avec
deux points fixes $0$ et $\infty$ (les seuls ``vrais points'' de la 
th\'eorie !); il y a donc des actions sur leurs voisinages ouverts invariants 
$\mathbf{C}$, $\mathbf{S} - \{0\}$ et l'intersection de ceux-ci
$\mathbf{C}^{*}$.\\

Diverses parties (en g\'en\'erales finies) $\Sigma \subset \mathbf{C}^{*}$
engendreront des spirales logarithmiques discr\`etes \'egalement invariantes
$q^{\mathbf{Z}}\Sigma$, voire des demi-spirales positivement ou
n\'egativement invariantes $q^{\mathbf{N}}\Sigma$, $q^{\mathbf{N}^{*}}\Sigma$, 
$q^{\mathbf{-N}}\Sigma$ etc...\\

Si $f$ est une fonction \`a valeurs scalaires, 
vectorielles ou matricielles sur une partie convenable (c'est \`a dire, en
fait, $q$-invariante) de $\mathbf{S}$, on note $\sigma_{q}f(z) = f(qz)$.

\subsection*{Corps des ``constantes''}

Le changement de variables $z = e^{2i\pi x}$
permet d'identifier le groupe $\mathbf{C}^{*}/q^{\mathbf{Z}}$ \`a la courbe 
elliptique $\mathbf{E} = \mathbf{C}/(\mathbf{Z} + \mathbf{Z}\tau)$ et 
$\mathcal{M}(\mathbf{C}^{*})^{\sigma_{q}}$ (corps des fonctions m\'eromorphes 
sur $\mathbf{C}^{*}$ et $\sigma_{q}$-invariantes) \`a 
$\mathcal{M}(\mathbf{E})$ 
(corps des fonctions elliptiques).\\

Soyons plus explicite sur ce point tr\`es important.
L'application $x \mapsto z = e^{2i\pi x}$ identifie 
$\mathcal{M}(\mathbf{C}^{*})$ au sous-corps de $\mathcal{M}(\mathbf{C})$
form\'e des fonctions $1$-p\'eriodiques. L'action de $\sigma_{q}$ sur
$\mathcal{M}(\mathbf{C}^{*})$ s'\'etend \`a $\mathcal{M}(\mathbf{C})$ par
l'action sur la variable : $x \mapsto x - \tau$ (qui est un rel\`evement de
$z \mapsto qz$). L'invariance sous $\sigma_{q}$ d'un \'el\'ement de
$\mathcal{M}(\mathbf{C}^{*})$ se traduit donc, dans $\mathcal{M}(\mathbf{C})$,
par la $\tau$-p\'eriodicit\'e; $\mathcal{M}(\mathbf{C}^{*})^{\sigma_{q}}$ est
ainsi identifi\'e au sous-corps de $\mathcal{M}(\mathbf{C})$ form\'e des
fonctions admettant $\mathbf{Z} + \mathbf{Z}\tau$ pour r\'eseau des 
p\'eriodes,
c'est \`a dire au corps $\mathcal{M}(\mathbf{E})$ des fonctions m\'eromorphes
sur le tore complexe $\mathbf{E} = \mathbf{C}/(\mathbf{Z} + \mathbf{Z}\tau)$
(fonctions elliptiques).

\subsection*{Fonctions de base}

Elles sont issues de la th\'eorie des fonctions th\^eta de jacobi.
La plus importante:
$$
\Theta_{q}(z) = \sum_{n \in Z}(-1)^{n}q^{-\frac{n^{2}-n}{2}}z^{n} =
\Theta_{q}(z) = (p;p)_{\infty} (z;p)_{\infty} (pz^{-1};p)_{\infty}
$$
D'un usage plus \'episodique:
$$
\Theta_{q}^{+}(z) = (z;p)_{\infty} = \prod_{r \geq 0} (1 - q^{-r}z)
$$
La famille des caract\`eres est \emph{en g\'en\'eral} (mais pas toujours !)
d\'efinie par
$$
e_{q,c} = \Theta_{q}/\Theta_{q,c} \text{   avec:   }
\Theta_{q,c}(z) = \Theta_{q}(c^{-1}z)
$$
Enfin, le succ\'edan\'e du logarithme:
$$
l_{q}(z) = z \Theta_{q}'(z)/\Theta_{q}(z)
$$
Pour diverses manipulations alg\'ebriques, on aura recours aux polyn\^omes de
Newton en $l_{q}$, que nous noterons:
$$
l_{q}^{(k)} = \frac{1}{k!} \prod_{0 \leq i \leq k-1}(l_{q} - i) =
\begin{pmatrix} l_{q} \\ k\\ \end{pmatrix} \quad k \in \mathbf{N}
$$

\subsection*{Corps de fonctions}

De fa\c{c}on g\'en\'erale, nous noterons $\mathcal{M}(X)$ le corps des
fonctions m\'eromorphes sur une surface de Riemann $X$.\\

Il y a tout d'abord les trois corps de coefficients:
$$
\mathbf{k} = \mathbf{C}(z) \subset
\mathbf{k}_{0} = \mathcal{M}(\mathbf{C}) \subset
\mathbf{k}'_{0} = \mathbf{C}(\{z\})
$$
ainsi que $\mathbf{k_{\infty}}$ et $\mathbf{k'_{\infty}}$ d\'efinis de 
mani\`ere \'evidente en posant $w = 1/z$. Ces sorps admettent les extensions:
$$
\mathbf{L}_{0} = 
\mathbf{k}_{0} (\Theta'_{q},(\Theta_{q,c})_{c \in \mathbf{C}^{*}}) \subset
\mathbf{L}'_{0} = 
\mathbf{k}'_{0} (\Theta'_{q},(\Theta_{q,c})_{c \in \mathbf{C}^{*}})
$$
et les sous-extensions (plus utiles)
$$
\mathbf{K}_{0} = 
\mathbf{k}_{0} (l_{q},(e_{q,c})_{c \in \mathbf{C}^{*}}) \subset
\mathbf{K}'_{0} = 
\mathbf{k}'_{0} (l_{q},(e_{q,c})_{c \in \mathbf{C}^{*}})
$$
Notons les inclusions 
$\mathbf{k}_{0} \subset \mathbf{K}_{0} \subset \mathbf{L}_{0} \subset
\mathcal{M}(\mathbf{C}^{*})$.

\subsection*{Matrices}

Matrices servant \`a la r\'eduction de Jordan:
$$
\xi_{\lambda,m} =
\begin{pmatrix}
\lambda & 1 & 0 & ... & 0 & 0\\
0 & \lambda & 1 & ... & 0 & 0\\
... & ... & ... & ... & .. & ...\\
... & ... & ... & ... & .. & ...\\
0 & 0 & 0 & ... & \lambda & 1\\
0 & 0 & 0 & ... & 0 & \lambda\\
\end{pmatrix}
\quad (taille \ m)
$$
$$
L_{m} = \xi_{1,m}^{l} =
\begin{pmatrix}
l^{(0)} & l^{(1)} & ... & l^{(m-1)}\\
0 & l^{(0)} & ... & l^{(m-2)}\\
... & ... & ... & ...\\
0 & 0 & ... & l^{(0)}\\
\end{pmatrix}
\quad
l^{(k)} = \frac{1}{k!} \prod_{0 \leq i \leq k-1}(l - i) =
\begin{pmatrix} l \\ k\\ \end{pmatrix} \quad k \in \mathbf{N}
$$
Dans cette derni\`ere, $l$ sera en g\'en\'eral instanci\'e par $l_{q}$.\\ 

Matrices servant \`a des algorithmes de r\'eduction du type pivot de Gauss:
$$
T_{i,\underline{\alpha}} =
\begin{pmatrix}
1 & 0 & ... & 0 & ... & 0 \\
0 & 1 & ... & 0 & ... & 0 \\
... & ... & ... & ... & ... & ... \\
\alpha_{1} & \alpha_{2} & ... & \alpha_{i} & ... & \alpha_{n} \\
... & ... & ... & ... & ... & ... \\
0 & 0 & ... & 0 & ... & 1 \\
\end{pmatrix}
\quad
D_{i,v} =
\begin{pmatrix}
1 & 0 & ... & 0 & ... & 0 \\
0 & 1 & ... & 0 & ... & 0 \\
... & ... & ... & ... & ... & ... \\
0 & 0 & ... & v & ... & 0 \\
... & ... & ... & ... & ... & ... \\
0 & 0 & ... & 0 & ... & 1 \\
\end{pmatrix}
$$

Enfin, matrices servant \`a la r\'eduction de Jordan ``renormalis\'ee'' au
chapitre 4:
$$
L'_{c,m}=
\begin{pmatrix}
l_{q}^{(0)} & (\eta/c)^{1} l_{q}^{(1)} & ... & (\eta/c)^{m-1} l_{q}^{(m-1)}\\
0 & l_{q}^{(0)} & ... & (\eta/c)^{m-2} l_{q}^{(m-2)}\\
... & ... & ... & ...\\
0 & 0 & ... & l_{q}^{(0)}\\
\end{pmatrix}
= (\xi'_{1,\eta/c,m})^{l_{q}} 
$$
$$
\xi'_{\lambda,\mu,m} =
\begin{pmatrix}
\lambda & \mu & 0 & ... & 0 & 0\\
0 & \lambda & \mu & ... & 0 & 0\\
... & ... & ... & ... & .. & ...\\
... & ... & ... & ... & .. & ...\\
0 & 0 & 0 & ... & \lambda & \mu\\
0 & 0 & 0 & ... & 0 & \lambda\\
\end{pmatrix}
= \mu \xi_{\lambda/\mu,m} 
$$
Ici, $\eta = q - 1$.


\addtocontents{toc}{\protect\contentsline{part}{\protect\numberline{}Solutions canoniques et classification}{\thepage}}
\part*{Solutions canoniques et classification}

\chapter{Th\'eorie de Frobenius: construction de solutions locales}

\section{Th\'eorie des fonctions}

On d\'efinit dans ce chapitre un catalogue de fonctions\footnote{Ici, 
fonctions de $z$ seul: le 
comportement en fonction de $q$ sera \'etudi\'e aux chapitres 3 et 4}
de base, dont on \'etablit les 
propri\'et\'es n\'ecessaires (analyticit\'e, p\^oles, \'equations 
fonctionnelles...). Le but est essentiellement d'avoir des
fonctions jouant un r\^ole analogue \`a celui de de $\log(z)$ et de 
$z^{\gamma}$ dans le cas des \'equations
diff\'erentielles lin\'eaires. Mais nos fonctions de base n'ont pas de
ramification : elles sont toutes dans $\mathcal{M}(\mathbf{C}^{*})$.\\

Avec les germes de fonctions holomorphes, elles serviront \`a
construire les solutions au voisinage de $0$ des \'equations aux 
$q$-diff\'erences fuchsiennes, et, plus g\'en\'eralement, singuli\`eres 
r\'eguli\`eres (cf. Introduction). Nous parlons de germes, car c'est la 
forme que nous presctivons a priori pour plagier la m\'ethode de Frobenius;
en r\'ealit\'e, il apparaitra a posteriori que l'on obtient des coefficients
\emph{m\'eromorphes sur $\mathbf{C}^{*}$}.


\subsection{$\Theta_{q}^{+}$ et $\Theta_{q}$}

On introduit, classiquement la fonction th\^eta (\emph{voir} 
\cite{GR}, \cite{RamisGrowth}, \cite{RamisJPRTraum}):
\begin{equation}
\Theta_{q}(z) = \sum_{n \in Z}(-1)^{n}q^{-\frac{n^{2}-n}{2}}z^{n}
\end{equation}
Elle est holomorphe sur $\mathbf{C}^{*}$ et v\'erifie les \'equations 
fonctionnelles :
$$\Theta_{q}(qz) = -qz\Theta_{q}(z)$$
$$\Theta_{q}(1/z) = -(1/z)\Theta_{q}(z)$$
Outre les deux
\'equations ci-dessus, elle poss\`ede de nombreuses propri\'et\'es 
multiplicatives, dont la plus notable est la formule du triple produit de
Jacobi (pour les notations, cf. Introduction):
\begin{equation}
\Theta_{q}(z) = (p;p)_{\infty} (z;p)_{\infty} (pz^{-1};p)_{\infty}
\end{equation}

En vue de l'\'etude des exemples, aux chapitres 3 et 4 (d\'eformation de
$(1 - z/z_{0})^{\alpha}$), on notera alors
$$
\Theta_{q}^{+}(z) = (z;p)_{\infty} = \prod_{r \geq 0} (1 - q^{-r}z)
$$
Celle-ci est caract\'eris\'ee par l'\'equation fonctionnelle :
$$
\Theta_{q}^{+}(qz) = (1 - qz)\Theta_{q}^{+}(z)
$$
Les z\'eros de $\Theta_{q}$ (resp. de $\Theta_{q}^{+}$) sont simples; ce sont
les \'el\'ements de $q^{\mathbf{Z}}$ (resp. de $q^{\mathbf{N}}$).\\

Les changements de variable $q = e^{-2i\pi\tau}$, $z = e^{2i\pi x}$
ram\`enent $\Theta_{q}$ aux fonctions th\^eta classiques de Jacobi. Cela sera
exploit\'e aux chapitres 3 et 4 pour \'etudier la d\'ependance en $q$ des
solutions d'\'equations aux $q$-diff\'erences.


\subsection{Caract\`eres}

Nous nommerons \emph{caract\`eres} les analogues, dans notre th\'eorie, des
fonctions $z^{\gamma}$ de la th\'eorie classique des \'equations 
diff\'erentielles. Un caract\`ere d'exposant $c$ est une solution 
m\'eromorphe dans 
$\mathbf{C}^{*}$ de l'\'equation fonctionnelle:
$$ 
\sigma_{q}f = cf ,\hspace{3cm}   (c \in \mathbf{C}^{*}) 
$$
Notons qu'il suffit d'exiger que $f$ soit m\'eromorphe dans un voisinage
\'epoint\'e de $0$ dans $\mathbf{C}^{*}$ ; la relation $f(z) = cf(z/q)$ permet
alors de prolonger $f$ r\'ecursivement en un unique caract\`ere.\\

Il d\'ecoulera de 1.1.5 que l'on ne peut
choisir les $e_{q,c}$ dans $\mathbf{C}(\{z\})$, mais cela se voit aussi
directement par l'unicit\'e de la s\'erie de Laurent.\\

L'\'equation fonctionnelle de $\Theta_{q}$ fournit une famille de caract\`eres
au comportement ``raisonnable'' en $0$ : on pose 
$e_{q,c} = \Theta_{q}/\Theta_{q,c}$, o\`u l'on a not\'e
$\Theta_{q,c}(z) = \Theta_{q}(c^{-1}z)$. De l'\'egalit\'e 
$\Theta_{q,c}(qz) = -qc^{-1}z\Theta_{q,c}(z)$ s'ensuit alors que $e_{q,c}$ 
est un caract\`ere. Il v\'erifie en outre l'\'equation fonctionnelle :
$$ 
e_{q,c}(1/z) = c/e_{q,c}(c z) 
$$
Son diviseur des z\'eros et des p\^oles (sur $\mathbf{C}^{*}$)\footnote{
Exceptionnellement, nous m\'elangeons la notation additive des diviseurs avec
la notation multiplicative de $\mathbf{C}^{*}$, ce qui peut pr\^eter \`a
confusion. Par la suite, nous ne consid\`ererons que des diviseurs sur
$\mathbf{E}$ et nous nous en tiendrons \`a la notation additive.} 
est
$$ 
\Div_{\mathbf{C}^{*}}(e_{q,c}) = 
\sum_{n \in Z}[q^{n}] - \sum_{n \in Z}[q^{n}c] 
$$
En tant que fonction automorphe pour $q^{\mathbf{Z}}$, $e_{q,c}$ a \'egalement
un diviseur des z\'eros et des p\^oles sur le tore
$\mathbf{E} = \mathbf{C}^{*}/q^{\mathbf{Z}} = 
\mathbf{C}/(\mathbf{Z} + \mathbf{Z}\tau)$:
$$ 
\Div_{\mathbf{E}}(e_{q,c}) = [0] - [\overline{c}]
$$
o\`u $\overline{c}$ d\'esigne la classe de $c$ dans $\mathbf{E}$.\\

Les \emph{relations en famille} entre les $e_{q,c}$ nous seront \'egalement
utiles :
$$ e_{q,1}(z) = z $$
$$ e_{q,qc}(z) = -c^{-1}ze_{q,c}(z) \in \mathbf{C}^{*}ze_{q,c}(z)$$
$$ e_{q,c^{-1}}(z) = 1/e_{q,c}(c z) = c^{-1}e_{q,c}(1/z) $$

Le choix d'une famille de caract\`eres est une \'etape pas du tout innocente
et qui appelle quelques commentaires.

\begin{enumerate}

\item{L'algorithme de r\'esolution locale d'une
\'equation aux $q$-diff\'erences fuchsienne en $0$ (cf. section 1.2) est
ind\'ependant du choix d'une famille particuli\`ere de caract\`eres index\'ee
par $\mathbf{C}^{*}$. Cette libert\'e de choix sera exploit\'ee au chapitre
3 pour controler la polarit\'e des solutions.}

\item{Une famille un peu plus compliqu\'ee aurait donn\'e la relation 
un peu plus agr\'eable: $e_{q,qc}(z) = ze_{q,c}(z)$. C'est le cas, par
exemple, si l'on prend 
$$
e_{q,c}(z) = 
\frac{\Theta_{q}(z)\Theta_{q}(cz^{2})}{\Theta_{q}(cz)\Theta_{q}(z^{2})} =
\frac{e_{q,c^{-1}}(z)}{e_{q,c^{-1}}(z^{2})}
$$
(variante: $\frac{e_{q,c}(z^{2})}{e_{q,c}(z)}$).
On ne peut malheureusement pas faire beaucoup plus simple!\\ 

En revanche, il ne nous est pas possible d'exiger la relation encore plus 
agr\'eable : 
$e_{q,c}e_{q,d} = e_{q,cd}$
pour $c,d \in \mathbf{C}^{*}$; ceci, contrairement aux axiomes impos\'es \`a
la famille des \emph{symboles} $e(c)$ dans \cite{SVdP},p. 150.\\

Soit en effet une famille de caract\`eres $e'_{c}$ satisfaisant cet axiome. 
Alors la fonction $f_{c} = e'_{c}/e_{q,c}$ est invariante sous $\sigma_{q}$ et 
peut donc \^etre consid\'er\'ee comme une fonction elliptique. De plus, 
$f_{c}f_{d}/f_{cd} = e_{q,cd}/e_{q,c}e_{q,d}$ par hypoth\`ese. 
Le diviseur de cette fonction sur $\mathbf{E}$ est
$[\overline{c}] + [\overline{d}] - [0] - [\overline{cd}]$. 
A un \'el\'ement
$c \in \mathbf{C}^{*}$, on associe le diviseur de degr\'e $0$ sur $\mathbf{E}$:
 $[\overline{c}] - [0] - \Div_{\mathbf{E}}(f_{c}) \in Div^{0}(\mathbf{E})$. 
L'application $\mathbf{C}^{*} \rightarrow Div^{0}(\mathbf{E})$ ainsi d\'efinie 
est alors un morphisme de groupes (\`a cause de l'\'egalit\'e ci-dessus). Il
est donc trivial, puisque la source est un groupe divisible et le but un groupe
ab\'elien libre. Mais, en le composant avec l'\'epimorphisme d'\'evaluation 
$Div^{0}(\mathbf{E}) \rightarrow \mathbf{E}$,
$\sum m_{i}[\alpha_{i}] \mapsto \sum m_{i}\alpha_{i}$, (dont le noyau est
form\'e des diviseurs des fonctions elliptiques et contient donc en particulier
$\Div_{\mathbf{E}}(f_{c})$), on obtient la projection canonique de
$\mathbf{C}^{*}$ sur $\mathbf{C}^{*}/q^{\mathbf{Z}}$, contradiction.\\

Un peu plus directement, on peut dire que $c \mapsto \Div_{\mathbf{E}}(e'_{c})$
serait un morphisme de groupes $\mathbf{C}^{*} \rightarrow Div^{0}(\mathbf{E})$
donc trivial; mais en m\^eme temps, chaque $\Div_{\mathbf{E}}(e'_{c})$ doit
\^etre \'equivalent \`a $\Div_{\mathbf{E}}(e_{q,c})$, contradiction.\\

En fait, il est facile de voir que l'application 
$(\overline{c},\overline{d}) \mapsto 
e_{q,c}e_{q,d}/e_{q,cd} \pmod {\mathbf{C}^{*}}$ 
est
bien d\'efinie et que c'est le cocycle associ\'e \`a la suite exacte (non 
scind\'ee !):
$$
\{1\} \rightarrow \mathcal{M}(\mathbf{E})/\mathbf{C}^{*} \rightarrow
Div^{0}(\mathbf{E}) \rightarrow \mathbf{E} \rightarrow \{0\}
$$
(\emph{voir} \cite{CartanEilenberg}).
De plus, toute famille de caract\`eres d\'efinit un cocycle \'equivalent.}

\item{Du point pr\'ec\'edent, il d\'ecoule que l'extension de $\mathbf{C}$
engendr\'ee par les $e_{q,c}e_{q,d}/e_{q,cd}$ est transcendante. Elle n'est pas
n\'ecessairement \'egale \`a $\mathcal{M}(\mathbf{E})$: par exemple en 
prenant les $e'_{c}(z) = e_{q,c'}(z^{2})$ ($c'$ \'etant une racine carr\'ee 
arbitraire de $c$), on obtient la p\'eriode $\tau/2$ et il est facile de
r\'ealiser ainsi toutes sortes d'isog\'enies. Le th\'eor\`eme de 
Riemann-Hurwitz (\emph{voir par exemple} \cite{Chevalley})
nous garantit que le corps obtenu est de genre $0$ ou $1$.
J'ignore s'il est possible d'obtenir le genre $0$.}

\end{enumerate}


\subsection{``Logarithme''}

Pour prendre en compte la partie unipotente des \'equations aux 
$q$-diff\'erences, nous aurons besoin, en parfaite analogie avec le r\^ole du 
logarithme dans la th\'eorie classique des \'equations diff\'erentielles, de 
r\'esoudre l'\'equation (non homog\`ene !):
\begin{equation}
\sigma_{q}f = f + 1
\end{equation}

Comme pour les caract\`eres, on obtient une solution \`a croissance 
``raisonnable'' en $0$ \`a l'aide de la fonction $\Theta_{q}$. 
D\'erivant logarithmiquement l'\'equation fonctionnelle satisfaite par 
$\Theta_{q}$, on trouve
$$ 
q \Theta_{q}'(qz)/\Theta_{q}(qz) = 1/z + \Theta_{q}'(z)/\Theta_{q}(z)
$$
Ceci conduit \`a introduire la solution :
$$ 
l_{q}(z) = z \Theta_{q}'(z)/\Theta_{q}(z) 
$$

De la formule du triple produit, on d\'eduit un d\'eveloppement de $l_{q}$ en
une s\'erie convergeant normalement sur tout compact de
$\mathbf{C}^{*} - q^{\mathbf{Z}}$ :
$$ 
l_{q}(z) = \sum_{r \geq 0} -q^{-r}z/(1 - q^{-r}z) +
\sum_{r \geq 1} q^{-r}z^{-1}/(1 - q^{-r}z^{-1}) 
$$

Outre l'\'equation requise, le ``logarithme'' $l_{q}$  v\'erifie
$$ 
l_{q}(1/z) = 1 - l_{q}(z) 
$$
Il est holomorphe sur $\mathbf{C}^{*} - q^{\mathbf{Z}}$ et admet des p\^oles
simples sur $q^{\mathbf{Z}}$.\\

En vue de l'\'etude des exemples, aux chapitres 3 et 4 (d\'eformation de
$\log(1 - z/z_{0})$), on consid\`ere \'egalement la 
partie holomorphe en $0$ de $l_{q}$ :
$$ 
l_{q}^{+}(z) = z (\Theta_{q}^{+})'(z)/\Theta_{q}^{+}(z) 
$$
Elle est holomorphe sur $\mathbf{C}^{*} - q^{\mathbf{N}}$, et s'y d\'eveloppe 
en une s\'erie normalement convergente sur tout compact :
$$ 
l_{q}(z) = \sum_{r \geq 0} -q^{-r}z/(1 - q^{-r}z) 
$$
Ses p\^oles sont simples, ce sont les \'el\'ements de $q^{\mathbf{N}}$. Elle
satisfait aux relations
$$ 
l_{q}(z) = l_{q}^{+}(z) - l_{q}^{+}(1/z) + 1/(1-z) 
$$
$$ l_{q}^{+}(qz) = -qz/(1-qz) + l_{q}^{+}(z) 
$$
Cette derni\`ere relation, avec l'\'egalit\'e $l_{q}^{+}(0) = 0$ et 
l'holomorphie (ou m\^eme simplement la continuit\'e en $0$ !), caract\'erise 
$l_{q}^{+}$. Donnons enfin son d\'eveloppement en s\'erie enti\`ere au
voisinage de $0$ (de rayon de convergence 1) :
$$ 
l_{q}(z) = \sum_{n \geq 1} -q^{n}z^{n}/(q^{n} - 1) 
$$


\subsection{Les constantes}

Dans cette section et la suivante, on note $\mathbf{k}'_{0}$ le corps 
$\mathbf{C}(\{z\})$ des s\'eries de Laurent en $0$, et $\mathbf{L}'_{0}$ 
l'extension de $\mathbf{k}'_{0}$ 
engendr\'ee par les
$\Theta_{q,c}$ ($c \in \mathbf{C}^{*}$) et $\Theta_{q}'$. Les \'el\'ements de
$\mathbf{L}'_{0}$ sont donc des germes de fonctions m\'eromorphes sur des 
voisinages 
\'epoint\'es de $0$ dans $\mathbf{C}^{*}$. On va \'etudier le sous-corps
$(\mathbf{L}'_{0})^{\sigma_{q}}$ des ``constantes'' de $\mathbf{L}'_{0}$, 
c'est \`a dire des invariants sous
l'action canonique de $\sigma_{q}$. Un germe 
$f \in (\mathbf{L}'_{0})^{\sigma_{q}}$ 
est tel
que $f(z) = f(z/q)$, ce qui permet r\'ecursivement de le prolonger en une
fonction m\'eromorphe sur $\mathbf{C}^{*}$, d'o\`u un morphisme injectif :
$(\mathbf{L}'_{0})^{\sigma_{q}} \hookrightarrow 
\mathcal{M}(\mathbf{C}^{*})^{\sigma_{q}}$:
on va voir que c'est en fait une \emph{\'egalit\'e}.\\

On a vu dans l'introduction comment l'application $x \mapsto z = e^{2i\pi x}$ 
permet d'identifier $\mathcal{M}(\mathbf{C}^{*})^{\sigma_{q}}$ au sous-corps 
de $\mathcal{M}(\mathbf{C})$ form\'e des
fonctions admettant $\mathbf{Z} + \mathbf{Z}\tau$ pour r\'eseau des 
p\'eriodes,
c'est \`a dire au corps $\mathcal{M}(\mathbf{E})$ des fonctions m\'eromorphes
sur le tore complexe $\mathbf{E} = \mathbf{C}/(\mathbf{Z} + \mathbf{Z}\tau)$
(fonctions elliptiques).\\

R\'eciproquement, la th\'eorie des fonctions elliptiques (\cite{Lang},
\cite{Mumford} p. 24, \cite{WW} p. 474) nous dit que tout \'el\'ement $f$ de
$\mathcal{M}(\mathbf{E})$ est, \`a un facteur constant pr\`es, de la forme
$$
\frac{\prod_{1 \leq i \leq r} \theta_{1}(x - \alpha_{i})}
{\prod_{1 \leq i \leq r} \theta_{1}(x - \beta_{i})}
$$

Dans cette \'ecriture,
$\sum_{1 \leq i \leq r} [\alpha_{i}] - \sum_{1 \leq i \leq r} [\beta_{i}]$ est
un rel\`evement dans $\mathbf{C}$ de $\Div_{\mathbf{E}}(f)$ choisi de sorte que
$\sum_{1 \leq i \leq r} \alpha_{i} = \sum_{1 \leq i \leq r} \beta_{i}$. Modulo
l'identification ci-dessus, cela signifie que tout \'el\'ement de
$\mathcal{M}(\mathbf{C}^{*})^{\sigma_{q}}$ est, \`a un facteur constant pr\`es,
de la forme 
$$
\frac{\prod_{1 \leq i \leq r} \Theta_{q,a_{i}}(z)}
{\prod_{1 \leq i \leq r} \Theta_{q,b_{i}}(z)}
$$
o\`u 
$\prod_{1 \leq i \leq r} a_{i} = \prod_{1 \leq i \leq r} b_{i}$. Ceci montre
que 
$\mathcal{M}(\mathbf{C}^{*})^{\sigma_{q}} = (\mathbf{L}'_{0})^{\sigma_{q}}$ 
et que ce corps
des constantes est engendr\'e par les produits de cette forme.\\

On voit m\^eme, par r\'ecurrence sur $r$, que l'on peut se restreindre aux 
$\Theta_{q,a}\Theta_{q,b}/\Theta_{q}\Theta_{q,ab}$. Si l'on a fait les
choix standards $e_{q,c} = \Theta_{q}/\Theta_{q,c}$, on en d\'eduit que les
$e_{q,ab}/e_{q,a}e_{q,b}$ ($a,b \in \mathbf{C}^{*}$) sont une famille de
g\'en\'erateurs.\\

Introduisons \`a pr\'esent

\begin{itemize}

\item{L'extension $\mathbf{K}'_{0}$ de $\mathbf{k}'_{0}$ 
engendr\'ee par les $e_{q,c}$ ($c \in \mathbf{C}^{*}$) et $l_{q}$;}

\item{l'extension $\mathbf{K}_{0}$ de 
$\mathbf{k}_{0} = \mathcal{M}(\mathbf{C})$ engendr\'ee par
ces m\^emes \'el\'ements: $\mathbf{K}_{0}$ est donc un sous-corps de
$\mathbf{K}'_{0}$.}

\end{itemize} 

On peut d'ailleurs se restreindre aux $c \neq 1$ et tels 
que $1 \leq |c| < |q|$, car $e_{q,1}$ et les $e_{q,qc}/e_{q,c}$ sont dans 
$\mathbf{k}_{0}$. Nous donnerons au 1.1.5 une vue d'ensemble des propri\'et\'es
alg\'ebriques de ces deux corps, \`a l'exception de l'\'etude de leurs
automorphismes qui sera abord\'ee dans \cite{JS2}.\\

On montrera au 1.2.4 que les coefficients des solutions locales en $0$ des
\'equations aux $q$-diff\'erences fuchsiennes en $0$ peuvent \^etre 
recherch\'ees dans $\mathbf{K}'_{0}$, et que l'\'equation fonctionnelle assure 
alors a posteriori qu'elles sont en fait \`a coefficients dans 
$\mathbf{k}'_{0}$.\\

Les corps $\mathbf{K}_{0}$ et $\mathbf{K}'_{0}$ sont des sous-corps de 
$\mathbf{L}$. Sous l'hypoth\`ese d'un choix standard des caract\`eres, 
il r\'esulte de ce qui pr\'ec\`ede qu'ils contiennent $\mathcal{M}(\mathbf{E})$
et sont donc des extensions de $\mathbf{L}^{\sigma_{q}}$. On en d\'eduit la\\

\textbf{\underline{Proposition}}\\

\emph{Tous ces corps de constantes sont \'egaux au corps des fonctions
elliptiques:
$$
(\mathbf{K}_{0})^{\sigma_{q}} = (\mathbf{K}'_{0})^{\sigma_{q}} = 
(\mathbf{L}'_{0})^{\sigma_{q}} = \mathcal{M}(\mathbf{E})
$$}
\hfill $\Box$


\subsection{Relations alg\'ebriques}

Les relations
$\Theta_{q,c} = e_{q,c}\Theta_{q}$ et 
$(\Theta_{q})'(z) = l_{q}(z) \Theta_{q}(z)/z$ montrent que 
$\mathbf{L}'_{0} = \mathbf{K}'_{0}(\Theta_{q})$. Cette extension est 
transcendante pure :\\

\textbf{\underline{Proposition 1}}\\

\emph{$\Theta_{q}$ est transcendant sur $\mathbf{K}'_{0}$ (et donc sur
$\mathbf{K}_{0}$).}\\

\textbf{\underline{Preuve}}\\

Elle repose sur les propri\'et\'es de croissance
de $\Theta_{q}$ au voisinage de $0$ (\emph{voir} \cite{RamisGrowth}). Mais ici,
comme dans les d\'emonstrations qui suivront, on a donn\'e aux arguments une
allure plus alg\'ebrique qui devrait en faciliter la r\'eutilisation.\\

En it\'erant l'\'equation fonctionnelle, on voit que 
$\Theta_{q}(q^{-n}z) = (-1/z)^{n}q^{n(n-1)/2}\Theta_{q}(z)$. Par ailleurs, 
pour un \'el\'ement $f \in \mathbf{K}'_{0}$ et un $z_{0} \in \mathbf{C}^{*}$ 
convenable 
fix\'e, la suite des $f(q^{-n}z_{0})$ a une croissance au pire simplement
exponentielle, comme on le voit en consid\'erant s\'epar\'ement les s\'eries
de Laurent des \'el\'ements de $\mathbf{k}'_{0}$, des caract\`eres $e_{q,c}$ 
et de $l_{q}$.\\

Si $\Theta_{q}$ \'etait alg\'ebrique sur $\mathbf{K}'_{0}$, on aurait
$\Theta_{q} \in \mathbf{K}'_{0}[\Theta_{q}^{-1}]$, et, fixant 
$z_{0} \in \mathbf{C}^{*} - q^{\mathbf{Z}}$, la relation \'ecrite plus haut 
entrainerait $q^{n^{2}} = O(q^{Cn})$ pour une certaine constante $C$, ce qui ne
se peut.\hfill $\Box$\\

Nous aurons besoin, pour le cas unipotent, d'introduire la notation
$$
l_{q}^{(k)} = \frac{1}{k!} \prod_{0 \leq i \leq k-1}(l_{q} - i) =
\begin{pmatrix} l_{q}\\ k\\ \end{pmatrix} \quad k \in \mathbf{N}
$$ 
Par convention, $l_{q}^{(k)} = 0$ pour $k < 0$. La formule de Pascal pour 
les coefficients
binomiaux montre que, quelque soit $k \in \mathbf{Z}$, 
$\sigma_{q} l_{q}^{(k)} = l_{q}^{(k)} + l_{q}^{(k-1)}$. Outre son application
ci-dessous, le th\'eor\`eme qui suit est fondamental pour la r\'esolution des
syst\`emes aux $q$-diff\'erences.\\

\textbf{\underline{Th\'eor\`eme}}\\

\emph{Les $l_{q}^{(k)}$ ($k \geq 0$) sont lin\'eairement
ind\'ependants sur $\mathbf{k}'_{0}$; autrement dit, $l_{q}$ est transcendant
sur $\mathbf{k}'_{0}$.}\\

\textbf{\underline{Preuve}}\\

Supposons donn\'ee une relation
$$
l_{q}^{(k+1)} = f_{0} l_{q}^{(0)} + ... + f_{k} l_{q}^{(k)}
$$
o\`u $f_{0}, ..., f_{k} \in $ et o\`u $k$ est minimum. On 
applique $\sigma_{q} - Id$ aux deux membres de cette \'egalit\'e, et l'on 
trouve :
$$
l_{q}^{(k)} = \sum_{0 \leq i \leq k} 
            ((\sigma_{q} f_{i} - f_{i})l_{q}^{(i)} + \sigma_{q} f_{i} l_{q}^{(i-1)}
)$$
Par minimalit\'e de $k$, $\sigma_{q} f_{k} - f_{k} = 1$ (sinon, on en tirerait 
une
\'ecriture de $l_{q}^{(k)}$ comme $\mathbf{k}'_{0}$-combinaison lin\'eaire des
$l_{q}^{(i)}$, $i < k$). Mais cette \'equation est impossible avec
$f_{k} \in \mathbf{k}'_{0}$.\hfill $\Box$\\

\textbf{\underline{Corollaire}}\\
 
\emph{Les seuls caract\`eres de $\mathbf{k}'_{0}(l_{q})$ sont 
les $cz^n$ ($c \in \mathbf{C}^{*}$, $n \in \mathbf{Z}$)}\\

\textbf{\underline{Preuve}}\\

Un \'el\'ement $f \in \mathbf{k}'_{0}(l_{q})^{*}$ s'\'ecrit, de
fa\c{c}on unique, $f = uA(l_{q})/B(l_{q})$, o\`u 
$u \in (\mathbf{k}'_{0})^{*}$, 
et o\`u
$A,B \in \mathbf{k}'_{0}[X]$ sont des polynomes unitaires et premiers 
entre eux.\\

Si 
maintenant $\sigma_{q}f = cf$ ($c \in \mathbf{C}^{*}$), la transcendance de
$l_{q}$ sur $\mathbf{k}'_{0}$ entraine la relation formelle :
$$
(\sigma_{q}u)A^{\sigma_{q}}(X+1)/B^{\sigma_{q}}(X+1) =
cuA^{\sigma_{q}}(X)/B^{\sigma_{q}}(X)
$$
dans laquelle $P^{\sigma_{q}}$ d\'esigne le r\'esultat de l'application de 
$\sigma_{q}$ aux
coefficients de $P \in \mathbf{k}'_{0}[X]$. Vues les hypoth\`eses sur 
$A$ et $B$, 
on a 
alors 
$$
\begin{cases}
\sigma_{q}u = cu \\
A^{\sigma_{q}}(X+1) = A(X) \\
B^{\sigma_{q}}(X+1) = B(X)
\end{cases}
$$

La premi\`ere \'egalit\'e montre que 
$u \in \mathbf{C}^{*}z^{n}$ et $c = q^{n}$ pour un $n \in \mathbf{Z}$ (car $u$
est une s\'erie de Laurent $\sum u_{k} z^{k}$ et il s'agit d'identifier terme 
\`a terme $\sum u_{k} q^{k} z^{k} = \sum c u_{k} z^{k}$). On n'a donc plus
qu'\`a prouver que la deuxi\`eme \'egalit\'e entraine que $A = 1$, le cas de 
$B$
\'etant le m\^eme.\\ 

Mais, si $A = X^{d} + a X^{d-1} + ...$, o\`u $d \geq 1$,
on voit que $\sigma_{q}a = a + d$, ce qui est tout \`a fait impossible pour un
\'el\'ement $a$ de $\mathbf{k}'_{0}$ (l'identification des s\'eries de Laurent 
donnerait, 
pour les termes de degr\'e $0$, $a_{0} = a_{0} + d$).\hfill $\Box$\\

Revenant \`a la structure alg\'ebrique de $\mathbf{K}_{0}$ et 
$\mathbf{K}'_{0}$, on va \`a pr\'esent 
d\'emontrer
que toutes les relations alg\'ebriques entre les $e_{q,c}$ sur 
$\mathbf{k}'_{0}(l_{q})$
d\'ecoulent de leurs relations multiplicatives. C'est une cons\'equence
imm\'ediate du\\

\textbf{\underline{Lemme}}\\

\emph{Soient $f_{1},...,f_{r}$ des fonctions
m\'eromorphes non nulles sur un voisinage \'epoint\'e de $0$ dans
$\mathbf{C}^{*}$. Soient $c_{1},...,c_{r}$ des \'el\'ements de 
$\mathbf{C}^{*}$ deux \`a deux non congrus modulo le sous-groupe 
$q^{\mathbf{Z}}$. On suppose que $\sigma_{q}f_{i} = c_{i}f_{i}$ (et donc
$f_{i} \in \mathcal{M}(\mathbf{C}^{*})$, pour $1 \leq i \leq r$. Alors les
$f_{i}$ sont lin\'eairement ind\'ependants sur $\mathbf{k}'_{0}(l_{q})$}\\

\textbf{\underline{Preuve}}\\

Dans le cas contraire, on aurait une relation
$$
P_{1}(l_{q})f_{1} + ... + P_{r}(l_{r})f_{r} = 0
$$
les $P_{i}$ \'etant des
\'el\'ements non tous nuls de $\mathbf{k}'_{0}[X]$. Si $r = 1$, la conclusion 
est triviale,
si $r = 2$, elle d\'ecoule du lemme pr\'ec\'edent.\\ 

Dans le cas g\'en\'eral, on
peut supposer la relation de taille minimale et les $P_{i}$ tous non nuls.
Appliquant $\sigma_{q}$, on en tire une nouvelle relation
$$
P_{1}(l_{q} + 1)c_{1}f_{1} + ... + P_{r}(l_{r} + 1)c_{r}f_{r} = 0
$$
n\'ecessairement proportionnelle \`a la premi\`ere (sinon on en obtiendrait
une plus courte par \'elimination).\\ 

Mais ceci signifie exactement que, pour
$1 \leq i < j \leq r$, 
$$
((P_{i}/P_{j})(l_{q}))^{\sigma_{q}} = (c_{j}/c_{i})(P_{i}/P_{j})(l_{q})
$$ 
et
entraine, d'apr\`es le lemme 1, que 
$c_{i} \equiv c_{j} \pmod {q^{\mathbf{Z}}}$,
contredisant l'hypoth\`ese.
\hfill $\Box$\\

La preuve de ce ``lemme d'ind\'ependance des caract\`eres'' \`a la Artin
permet en fait de donner des g\'en\'erateurs de toutes les relations
entre les $f_{i}$.\\

\textbf{\underline{Proposition 2}}\\

\emph{Supposons :
$\lambda e_{q,c_{1}}...e_{q,c_{r}} = e_{q,d_{1}}...e_{q,d_{s}}$, o\`u, pour
$1 \leq i \leq r$ et $1 \leq j \leq s$, on a $c_{i} \neq 1$, $d_{j} \neq 1$ et
$c_{i} \neq d_{j}$. Alors, $r = s$, et, \`a r\'eindexation pr\`es,
$d_{i} = q^{m_{i}}c_{i}$ ($1 \leq i \leq r$), avec 
$m_{1} + ... + m_{r} = 0$}\\

\textbf{\underline{Preuve}}\\

La consid\'eration des diviseurs de z\'eros et de
p\^oles montre que 
$$
\sum_{1 \leq i \leq r} ([0] - [\overline{c_{i}}]) =
\sum_{1 \leq j \leq s} ([0] - [\overline{d_{j}}])
$$
d'o\`u $r = s$ et, \`a
r\'eindexation pr\`es, $d_{i} = q^{m_{i}}c_{i}$ ($1 \leq i \leq r$).\\ 

On prouve
alors par r\'ecurrence, \`a partir des relations donn\'ees au 1.1.2, que
$$
e_{q,q^{m}c} = (- c_{i} z)^{m} q^{m(m-1)/2} e_{q,c}
$$
d'o\`u l'on tire :
$$
e_{q,d_{i}}/e_{q,c_{i}} = (- c_{i} z)^{m_{i}} q^{m_{i}(m_{i}-1)/2}
$$ 
L'\'egalit\'e $m_{1} + ... + m_{r} = 0$ vient alors imm\'ediatement.
\hfill $\Box$\\

R\'eciproquement, on voit que ces relations entrainent l'\'egalit\'e de
d\'epart avec
$$
\lambda = q^{\sum_{1 \leq i \leq r}m_{i}^{2}/2}\prod_{1 \leq i \leq r} c_{i}^{-m{i}}
$$



\section{Forme g\'en\'erale des solutions au voisinage de $0$}


\subsection{Pr\'eliminaires}

Afin de mettre en valeur le caract\`ere essentiellement alg\'ebrique et
algorithmique de la r\'esolution locale, nous pr\'esenterons celle-ci de
fa\c{c}on quelque peu axiomatique \emph{dans section 1.2 seulement}.\\

On se donne donc une extension $\mathbf{K}$ de 
$\mathbf{k}'_{0} = \mathbf{C}(\{z\})$ munie d'un automorphisme \'etendant 
$\sigma_{q}$, que l'on notera simplement $\sigma$ ; celui-ci agit donc 
\'egalement sur $\mathbf{K}^{n}$ et sur $GL_{n}(\mathbf{K})$. On exigera que 
$\mathbf{K}$ contienne :

\begin{itemize}

\item{pour chaque $c \in \mathbf{C}^{*}$, un \'el\'ement 
$e_{c} \in \mathbf{K}^{*}$ tel que $\sigma e_{c} = c e_{c}$}

\item{un \'el\'ement $l$ tel que $\sigma l = l + 1$ (et donc 
$l \in \mathbf{K}^{*}$)}

\end{itemize}

Le corps $\mathbf{K}$, l'automorphisme $\sigma$ et les symboles $e_{c}$ et $l$
peuvent \^etre respectivement instanci\'es par $\mathbf{K}'_{0}$, $\sigma_{q}$,
les $e_{q,c}$ et $l_{q}$ \'etudi\'es au 1.1, mais ce n'est pas fondamental, le
contenu de cette section \'etant plut\^ot de nature alg\'ebrique. Cependant,
les raisonnements du 1.1.5 s'appliquent sous les seules conditions ci-dessus.
Par exemple, ni $l$ ni les $e_{c}$ ne sont \'el\'ements de  
$\mathbf{k}'_{0}$ et
$l$ est m\^eme transcendant sur $\mathbf{k}'_{0}$. On introduit encore
$$
l^{(k)} = \frac{1}{k!} \prod_{0 \leq i \leq k-1}(l - i) =
\begin{pmatrix} l \\ k\\ \end{pmatrix} \quad k \in \mathbf{N}
$$

Les $l^{(k)} \;,\; k \in \mathbf{N}$ sont donc lin\'eairement ind\'ependants
sur $\mathbf{k}_{0}$.\\

On \'etudie le syst\`eme aux $q$-diff\'erences lin\'eaire homog\`ene
\begin{equation}
\sigma X = A X
\end{equation}
dans lequel $A \in GL_{n}(\mathbf{C}\{z\})$ (en particulier,
$A(0) \in GL_{n}(\mathbf{C})$).On recherche des solutions vectorielles
$X \in \mathbf{K}^{n}$. S'il y en a suffisamment, on peut former un 
``syst\`eme fondamental de solutions'' et donc aussi une ``solution 
fondamentale'', c'est \`a dire une solution matricielle
$X \in GL_{n}(\mathbf{K})$. Toute autre solution $X'$ de la m\^eme \'equation
s'\'ecrit donc $XV$ avec $V \in GL_{n}(\mathbf{K})$. Mezalor, les \'egalit\'es
$A = (\sigma X) X^{-1} = (\sigma XV) (XV)^{-1}$ montrent que $\sigma V = V$,
autrement dit, $V$ est \`a coefficients dans le corps $\mathbf{K}^{\sigma}$ des
invariants de $\mathbf{K}$ sous l'action de $\sigma$: ce corps tient lieu de
corps des constantes de la th\'eorie.\\

La m\'ethode de r\'esolution s'apparente \`a la m\'ethode de Frobenius pour les
\'equations diff\'erentielles (\cite{Ince}, \cite{Wasow}): notion d'exposant; 
m\'ethode des coefficients ind\'etermin\'es pour la r\'esolution formelle dans 
le cas non-r\'esonnant; m\'ethode des majorants pour \'etablir la convergence;
g\'en\'ericit\'e du cas semi-simple, o\`u $l$ ne joue aucun r\^ole.\\

La r\'esolution est enti\`erement algorithmique et fournit, dans tous les cas,
un syst\`eme fondamental de solutions sous une forme analytique pr\'ecise. On
donnera au 1.3 une version matricielle de la r\'esolution, permettant, dans le
cas non-r\'esonnant, le choix d'une solution fondamentale canonique.


\subsection{R\'esolution dans le cas non-r\'esonnant}

On notera dans cette section $A(z) = A_{0} + A_{1}z + A_{2}z^{2} + ... $.
Ainsi, les $A_{k}$ sont \'el\'ements de $M_{n}(\mathbf{C})$ et
$A_{0} = A(0) \in GL_{n}(\mathbf{C})$.\\

On fixe une norme quelconque sur $\mathbf{C}^{n}$, que l'on note $\|\ \|$, et
l'on note de la m\^eme fa\c{c}on la norme subordonn\'ee sur 
$M_{n}(\mathbf{C})$, d\'efinie par 
$\|B\| = \sup \{ \|BX\| \;/\; X \in \mathbf{C}^{n} \;,\;\|X\| = 1 \}$.\\

On fera, \emph{dans cette section}, l'hypoth\`ese de non-r\'esonnance 
suivante:\\

\emph{Soient $\lambda \not= \mu \in Sp(A_{0})$ deux valeurs propres distinctes de
$A_{0}$; alors $\lambda \not\equiv \mu \pmod{q^{\mathbf{Z}}}$ (rappelons que
$\lambda \;,\; \mu \in \mathbf{C}^{*}$).}\\

Cette hypoth\`ese sera lev\'ee \`a la 
section 1.2.3. Elle est analogue \`a la condition classique de non-r\'esonnance
pour les \'equations diff\'erentielles lin\'eaires (valeurs propres distinctes
ne diff\'erant pas d'un entier). On admet cependant la possibilit\'e de valeurs
propres multiples: outre les caract\`eres, elles feront apparaitre des 
``logarithmes'' $l$ dans les solutions.\\

\textbf{\underline{Th\'eor\`eme}}\\

\emph{Soit $c \in Sp(A_{0})$ une valeur propre de
$A_{0}$ (de sorte que $c \in \mathbf{C}^{*}$). On d\'ecompose le sous-espace 
caract\'eristique associ\'e de $\mathbf{C}^{n}$ en sous-espaces cycliques:
$$
E_{c} = F_{1} \oplus ... \oplus F_{r}
$$
de dimensions $f_{i} = dim_{\mathbf{C}} F_{i}$. A tout choix de vecteurs
cycliques des $F_{i}$ est alors canoniquement associ\'e une famille
$(e_{c}Y_{c,i,j})_{1 \leq i \leq r,0 \leq j < f_{i}}$ de solutions 
$\mathbf{C}(\{z\})$-lin\'eairement ind\'ependantes, o\`u
$$
Y_{c,i,j} \in \sum_{k = 0}^{j} (\mathbf{C}(\{z\})l^{(k)})\mathbf{C}^{n}
$$
}

Pour prouver ce th\'eor\`eme, on op\`ere le changement d'inconnue
$X = e_{c}Y$, qui ram\`ene \`a l'\'equation 
$\sigma Y =c^{-1}AY$. Si l'on prend pour $c$ l'un des exposants en 
$0$, c'est \`a dire un \'el\'ement de $Sp(A_{0})$, cela revient donc \`a
supposer que l'exposant consid\'er\'e est $c = 1$, ce que nous ferons
jusqu'\`a la fin de la preuve. L'hypoth\`ese de non-r\'esonnance nous
dit alors qu'aucun $q^{n} \;(n \in \mathbf{Z} - \{0\})$ n'est valeur
propre de $A_{0}$.\\

\textbf{\underline{Lemme 1}}\\

\emph{On se donne
$$
Y = Y_{0} + z Y_{1} + ... \in (\mathbf{C}\{z\})^{n}
$$
Soit $X_{0} \in \mathbf{C}^{n}$ tel que $X_{0} = A_{0}X_{0} - Y_{0}$.
L'\'equation $\sigma X = A X - Y$ admet alors une unique solution 
formelle
$$
X = X_{0} + z X_{1} + ... \in (\mathbf{C}[[z]])^{n}
$$
De plus, cette solution converge.}\\

\textbf{\underline{Preuve du lemme 1}}\\

Elle comporte les deux \'etapes habituelles !

\begin{itemize}

\item{Premi\`ere \'etape: existence et unicit\'e formelles.\\

L'\'equation se r\'eecrit: 
$$
\forall m \geq 0 \;,\; q^{m}X_{m} = 
\sum_{i+j = m} A_{i} X_{j} - Y_{m}
$$
Pour $m = 0$, c'est la condition pos\'ee en hypoth\`ese. Pour $m \geq 1$,
cela \'equivaut \`a 
$$
(q^{m}Id_{n} - A_{0})X_{m} = \sum_{i = 1}^{m} A_{i}X_{m-i} - Y_{m}
$$
L'hypoth\`ese de non-r\'esonnance nous garantit alors pr\'ecis\'ement
l'existence et l'unicit\'e d'un tel $X_{m} \in \mathbf{C}^{n}$. On 
conclut, par r\'ecurrence, \`a l'existence et \`a l'unicit\'e d'une 
solution formelle
$$
X = X_{0} + z X_{1} + ... \in (\mathbf{C}[[z]])^{n}
$$}

\item{Deuxi\`eme \'etape: convergence de la solution formelle; bien s\^ur,
par la m\'ethode des s\'eries majorantes !\\

On notera ici 
$x_{m} = \|X_{m}\| \;,\; y_{m} = \|Y_{m}\| \;,\, a_{m} = \|A_{m}\|$ et, 
pour $m \geq 1$, $b_{m} = \|(q^{m}Id_{n} - A_{0})^{-1}\|$. Comme
$\underset{m \to \infty}{\lim} b_{m} = 0$, on peut de plus introduire
$b = \underset{m \geq 1}{\sup} b_{m} < \infty$. On a alors, pour $m \geq 1$:
$$
x_{m} \leq b(\sum_{j = 0}^{m-1} a_{m-j} x_{j} + y_{m})
$$
On introduit la suite des majorants $\overline{x}_{m}$ en posant
$\overline{x}_{0} = x_{0}$ et, pour $m \geq 1$, la relation de 
r\'ecurrence:
$$
\overline{x}_{m} = b(\sum_{j = 0}^{m-1} a_{m-j} \overline{x}_{j} + y_{m})
$$
On a ainsi 
$\forall m \in \mathbf{N} \;,\; 0 \leq x_{m} \leq \overline{x}_{m}$. La
s\'erie g\'en\'eratrice 
$\overline{\mathcal{X}}(z) = \sum_{m \geq 0} \overline{x}_{m} z^{m}$ 
majore $\mathcal{X}(z) = \sum_{m \geq 0} x_{m} z^{m}$. Elle est, de plus,
solution de l'\'equation
$$
\overline{\mathcal{X}}(z) = x_{0} +
b \sum_{m \geq 1} y_{m} z^{m} +
(\sum_{i \geq 1} a_{i} z^{i})\overline{\mathcal{X}}(z)
$$
C'est donc le d\'eveloppement en s\'erie enti\`ere \`a l'origine de la 
fonction holomorphe
$$
\overline{\mathcal{X}}(z) = 
\frac{x_{0} + b \sum_{m \geq 1} y_{m} z^{m}}
{1 - z\sum_{i \geq 0} a_{i+1} z^{i}}
$$
Le rayon de convergence de celle-ci est non nul, il en va donc de m\^eme
de la s\'erie major\'ee $\mathcal{X}(z)$; ceci entraine bien la
convergence de $\sum_{m \geq 0} z^{m} X_{m}$.
\hfill $\Box$\\
}

\end{itemize}

Remarquons que le lemme, appliqu\'e au cas $Y = 0$, permet de prolonger
tout vecteur propre $X_{0}$ de $A_{0}$ en une solution convergente
$X = X_{0} + z X_{1} + ...$ de $\sigma X = A X$ : il r\`egle donc
compl\`etement le cas o\`u $A_{0}$ est semi-simple.\\

\textbf{\underline{Lemme 2}}\\

\emph{On se donne
$X_{0}^{(0)} \;,\; ... \;,\; X_{0}^{(j)} \in \mathbf{C}^{n}$ tels que
$$
\begin{cases}
X_{0}^{(0)} = A_{0} X_{0}^{(0)}\\
X_{0}^{(k)} = A_{0} X_{0}^{(k)} - X_{0}^{(k-1)} \quad (1 \leq k \leq j)
\end{cases}
$$
Il existe alors une unique solution 
$X = l^{0}X^{(0)} \;+\;...\;+\; l^{j}X^{(j)}$ de l'\'equation 
$\sigma X = A X$ o\`u chaque 
$X^{(k)} \in (\mathbf{C}\{z\})^{n}$ ($0 \leq k \leq j$) est convergent
et de la forme $X^{(k)} = X_{0}^{(k)} + z X_{1}^{(k)} + ...$.}\\

\textbf{\underline{Preuve du lemme 2}}\\

Rempla\c{c}ant $X$ par 
l'expression
prescrite, et utilisant la $\mathbf{C}(\{z\})$-lin\'eaire ind\'ependance
des $l^{k}$, on est ramen\'e au syst\`eme d'\'equations:
$$
\sigma X{(0)} = A X^{(0)}
$$
$$
\sigma X^{(k)} = A X^{(k)} - \sigma X^{(k-1)} 
\quad (1 \leq k \leq j)
$$
On d\'etermine $X^{(0)}$ en appliquant le lemme 1 avec 
$X_{0} = X_{0}^{(0)}$ et $Y = 0$. On d\'etermine de m\^eme $X^{(k}$
($1 \leq k \leq j$) en prenant $X_{0} = X_{0}^{(k)}$ et 
$Y = - \sigma X^{(k-1)}$.\hfill $\Box$\\

On peut maintenant reprendre la preuve du th\'eor\`eme. On choisit dans
chaque $F_{i}$ ($1 \leq i \leq r$) un vecteur cyclique $X_{0}^{(i)}$.
Autrement dit, $(A_{0} - Id_{n})^{f_{i}}X_{0}^{(i)} = 0$ et les
$X_{0}^{(i,k)} := (A_{0} - Id_{n})^{f_{i}-1-k}X_{0}^{(i)}$
($0 \leq k \leq f_{i} - 1$) forment une base de $F_{i}$. On applique
alors le lemme 2 successivement pour chaque $i \in \{1,...,r\}$ avec 
$j = f_{i} - 1$ et $X_{0}^{(k)} = X_{0}^{(i,k)}$ 
($0 \leq k \leq f_{i} - 1$).\\

L'ind\'ependance lin\'eaire sur $\mathbf{C}(\{z\})$ des solutions ainsi
obtenues d\'ecoule facilement de celle des $l^{k}$ et de la
$\mathbf{C}$-ind\'ependance lin\'eaire des termes constants. Ceci 
ach\`eve la preuve du th\'eor\`eme.\hfill $\Box$\\

\textbf{\underline{Corollaire}}\\

\emph{On obtient ainsi un syst\`eme fondamental de solutions.}\\

En effet, lorsque $c$ varie dans $Sp(A_{0})$, leur nombre total obtenu est
$\sum dim(E_{c}) = n$ et leur ind\'ependance lin\'eaire se d\'eduit du lemme 2
de 1.1.5.\\

\textbf{\underline{Remarque}}\\

On a utilis\'e en fait la propri\'et\'e suivante, un peu moins forte que la
non-r\'esonnance: aucun $q^{n} \;(n \in \mathbf{N} - \{0\})$ n'est valeur 
propre de $A_{0}.$ Il s'ensuit que les arguments de cette section s'appliquent 
en fait \`a tout exposant $c$ tel que $cq \;,\; cq^{2} \;,\;...$ ne sont pas
exposants sans qu'il soit besoin d'appliquer l'algorithme de  pr\'eparation de 
la section suivante.


\subsection{L'algorithme de pr\'eparation}

On l\`eve ici l'hypoth\`ese de non-r\'esonnance : les valeurs propres de
$A_{0}$ sont quelconques (non nulles!). On utilisera des transformations
de jauge $X = M Y$ qui remplacent l'\'equation $\sigma X = A X$ par
l'\'equation $\sigma Y = B Y$ o\`u $B = (\sigma M)^{-1} A M$.\\

\textbf{\underline{Th\'eor\`eme}}\\

\emph{En alternant des transformations
de jauge \`a matrices constantes 
$Q_{1} \;,\;...\;,\; Q_{r} \in GL_{n}(\mathbf{C})$ et \`a matrices de
``shearing'' (voir ci-dessous) $S_{1} \;,\; ... \;,\; S_{r}$, on ram\`ene
l'\'equation en $X$ \`a une \'equation non r\'esonnante \`a l'origine en
$Y$, o\`y $X = Q_{1}S_{1} \;...\; Q_{r}S_{r} Y$.}\\

\textbf{\underline{Preuve}}\\

Supposons que $Sp(A_{0})$ contienne une 
classe non triviale modulo $q^{\mathbf{Z}}$: $\{c \;,\;...\;,\;cq^{m}\}$
($m \in \mathbf{N}^{*}$).\\

Il existe une matrice $Q \in GL_{n}(\mathbf{C})$ telle que l'on peut 
\'ecrire
$$
Q^{-1}A_{0}Q =
\begin{pmatrix}
a_{0} & 0\\
0 & d_{0}\\
\end{pmatrix}
$$
o\`u le bloc $a_{0}$ est triangulaire sup\'erieur, de taille $\mu$ et n'a
que des $cq^{m}$ sur la diagonale et le bloc $d_{0}$ est triangulaire 
sup\'erieur, de taille $\nu$ et n'a aucun $cq^{m}$ sur la diagonale. On
introduit la matrice dite \emph{de shearing} (c'est \`a dire, de cisaillement) 
(\emph{voir} \cite{Wasow}, p.94)
$$
S =
\begin{pmatrix}
zId_{\mu} & 0\\
0 & Id_{\nu}\\
\end{pmatrix}
$$

Ecrivons maintenant, avec des blocs de tailles correspondantes,
$$
Q^{-1}AQ =
\begin{pmatrix}
a & b\\
c & d\\
\end{pmatrix}
$$
de sorte que
$$
B = (\sigma(QS))^{-1} A (QS) =
\begin{pmatrix}
a' & b'\\
c' & d'\\
\end{pmatrix} =
\begin{pmatrix}
q^{-1}a & q^{-1}z^{-1}b\\
zc & d\\
\end{pmatrix}
$$

L'hypoth\`ese faite sur $Q^{-1}A_{0}Q$ montre que les coefficients de 
$b'$ (et, de fa\c{c}on \'evidente, ceux de $a' \;,\; c' \;,\; d'$) 
appartiennent \`a $\mathbf{C}\{z\}$. Le terme constant de $B_{0}$ est donc de 
la forme
$$
\begin{pmatrix}
q^{-1} a_{0} & *\\
0 & d_{0}\\
\end{pmatrix}
$$

Ainsi, le spectre de $B_{0}$ est le m\^eme que celuui de $A_{0}$, sauf en ce
qui concerne les $cq^{m}$ qui ont \'et\'e remplac\'es par $cq^{m-1}$.
\hfill $\Box$\\

On voit m\^eme que le nombre $r$ d'\'etapes $QS$ n\'ecessaires est major\'e par
la somme des $(card(C_{i}) - 1)$ o\`u les $C_{i}$ sont les classes d'exposants
modulo $q^{\mathbf{Z}}$.\\

\textbf{\underline{Remarque 1}}\\

A cause de la n\'ecessit\'e de cette 
pr\'eparation, il semble difficile de pr\'evoir sur la forme d'une matrice
quelconque la pr\'esence de contributions ``logarithmiques'' dans les
solutions.\\

\textbf{\underline{Exemple 1}}

$$
A =
\begin{pmatrix}
q & z\\
0 & 1\\
\end{pmatrix}
\qquad
A_{0} =
\begin{pmatrix}
q & 0\\
0 & 1\\
\end{pmatrix}
$$
On prend $Q = Id_{2}$ et
$$
S =
\begin{pmatrix}
z & 0\\
0 & 1\\
\end{pmatrix}
\Rightarrow
B = (\sigma(QS))^{-1} A (QS) =
\begin{pmatrix}
1 & q^{-1}\\
0 & 1\\
\end{pmatrix}
= B_{0}
$$
La solution d'exposant $1$ comportera donc une partie logarithmique. On obtient
comme matrice fondamentale (cf prochaine section):
$$
\begin{pmatrix}
1 & zl/q\\
0 & 1\\
\end{pmatrix}
$$

\textbf{\underline{Exemple 2}}

$$
A' =
\begin{pmatrix}
q & z^{2}\\
0 & 1\\
\end{pmatrix}
\qquad
A'_{0} =
\begin{pmatrix}
q & 0\\
0 & 1\\
\end{pmatrix}
$$
La m\^eme transformation donne
$$
B' = (\sigma(QS))^{-1} A' (QS) =
\begin{pmatrix}
1 & q^{-1}\\
0 & 1\\
\end{pmatrix}
$$
Donc $B'_{0} = Id_{2}$: la solution d'exposant $1$ ne comporte pas de partie 
logarithmique. On obtient comme matrice fondamentale:
$$
\begin{pmatrix}
1 & z^{2}/(q^{2} -q)\\
0 & 1\\
\end{pmatrix}
$$

\textbf{\underline{Remarque 2}}\\

En revanche, l'algorithme de pr\'eparation ne
modifie pas les p\^oles de $A(z)$ dans $\mathbf{C}^{*}$, dont on verra qu'ils
sont d\'eterminants pour comprendre la signification g\'eom\'etrique des
solutions.\\

De m\^eme, le spectre de $A_{0}$ ne changeant pas modulo $q^{\mathbf{Z}}$, la
partie ``caract\`eres'' des solutions n'est pas essentiellement affect\'ee
(\`a \'equivalence m\'eromorphe pr\`es).


\subsection{Solutions canoniques, matrices fondamentales}

Dans le cas non-r\'esonnant, on a construit une base de $\mathbf{K}^{n}$
form\'ee de solutions de l'\'equation $\sigma X = A X$. Cette base est
uniquement d\'etermin\'ee par la donn\'ee d'une d\'ecomposition de chacun des
sous-espaces caract\'eristiques $E_{c} \subset \mathbf{C}^{n}$ 
($c \in Sp(A_{0})$) en sous-espaces cycliques et du choix d'un vecteur cyclique
pour chacun de ces derniers (ceci, pour l'action de $A_{0}$ sur 
$\mathbf{C}^{n}$).\\

Mais ces donn\'ees \'equivalent en fait \`a la mise de $A_{0}$ sous forme de
Jordan :
$$
Q_{0}^{-1}A_{0}Q_{0} =
\begin{pmatrix}
c_{1} \xi_{1,m_{1}} & 0 & ... & 0\\
0 & c_{2} \xi_{1,m_{2}} & ... & ...\\
... & ... & ... & ...\\
0 & 0 & ... & c_{k} \xi_{1,m_{k}}\\
\end{pmatrix}
$$
$$
\xi_{\lambda,m} =
\begin{pmatrix}
\lambda & 1 & 0 & ... & 0 & 0\\
0 & \lambda & 1 & ... & 0 & 0\\
... & ... & ... & ... &.. & ...\\
... & ... & ... & ... &.. & ...\\
0 & 0 & 0 & ... & \lambda & 1\\
0 & 0 & 0 & ... & 0 & \lambda\\
\end{pmatrix}
\quad (taille \ m)
$$
Le choix de blocs de la forme $c \xi_{1,m}$ au lieu des habituels $\xi_{c,m}$
est rendu possible par cette circonstance que les exposants $c_{i}$ ne sont pas
nuls.\\

Les colonnes de $Q_{0}$ sont les vecteurs $X_{0}^{(i,k)}$ attach\'es aux divers
blocs. Leurs prolongements holomorphes $X^{(i,k)} = X_{0}^{(i,k)} + ...$
construits \`a la section 1.2.2 sont les colonnes d'une matrice
$Q = Q_{0} + z Q_{1} + ... \in GL_{n}(\mathbf{C}\{z\})$. Les r\'esultats de
1.2.2 prennent la forme matricielle suivante. Notant $X$ la matrice solution
fondamentale dont les colonnes sont les solutions $X_{1} \;,;...\;,\; X_{n}$
construites en 1.2.2, on a: 
$$
X = Q
\begin{pmatrix}
L_{m_{1}} & 0 & ... & 0\\
0 & L_{m_{2}} & ... & 0\\
... & ... & ... & ...\\
0 & 0 & ... & L_{m_{k}}\\
\end{pmatrix}
\begin{pmatrix}
e_{c_{1}}I_{m_{1}} & 0 & ... & 0\\
0 & e_{c_{2}}I_{m_{2}} & ... & 0\\
... & ... & ... & ...\\
0 & 0 & ... & e_{c_{k}}I_{m_{k}}\\
\end{pmatrix}
$$
On a ici introduit les ``blocs logarithmiques''
$$
L_{m} = \xi_{1,m}^{l} =
\begin{pmatrix}
l^{(0)} & l^{(1)} & ... & l^{(m-1)}\\
0 & l^{(0)} & ... & l^{(m-2)}\\
... & ... & ... & ...\\
0 & 0 & ... & l^{(0)}\\
\end{pmatrix}
$$

Il sera parfois commode d'\'ecrire $Q = Q_{0} H$, o\`u 
$H \in GL_{n}(\mathbf{C}\{z\})$ et $H(0) = Id_{n}$. On peut alors synth\'etiser
les r\'esultats obtenus en un\\

\textbf{\underline{Th\'eor\`eme}}\\

\emph{A toute d\'ecomposition de Jordan de
$A_{0}$ de matrice de passage $Q_{0}$ est associ\'ee une matrice solution
fondamentale $X = Q_{0} H L C \in GL_{n}(\mathbf{K})$ de l'\'equation
$\sigma X = AX$, o\`u:
\begin{enumerate}
\item{$H \in GL_{n}(\mathbf{k}_{0})$ et $H(0) = Id_{n}$}
\item{$L$ est une diagonale de blocs logarithmiques $L_{m}$ dont les tailles 
sont celles des blocs de Jordan}
\item{$C$ est une diagonale de caract\`eres $e_{c}$, correspondant aux
exposants $c \in Sp(A_{0})$ de l'\'equation \`a l'origine. Les matrices $L$ et 
$C$ commutent.}
\end{enumerate}
L'\'ecriture qui pr\'ec\`ede est valide dans le cas non-r\'esonnant. Dans le
cas g\'en\'eral, cette \'ecriture doit \^etre pr\'ec\'ed\'ee de 
$Q_{1}S_{1}...Q_{r}S_{r}$, o\`u les $Q_{i} \in GL_{n}(\mathbf{C})$ et o\`u les
$S_{i}$ sont des matrices de shearing
$$
S_{i} =
\begin{pmatrix}
z Id_{\mu_{i}} & 0\\
0 & Id_{\nu_{i}}\\
\end{pmatrix}
$$
}\\

Il sera important de diff\'erencier la partie m\'eromorphe en $0$ de la 
solution,
c'est \`a dire $Q_{0} H$, de la partie ``log-car'', dont les coefficients
sont des combinaisons lin\'eaires des caract\`eres $e_{c}$ sur 
$\mathbf{C}[l]$. Outre la commutation de $L$ et $C$, cette partie est 
solution d'une \'equation aux $q-$diff\'erences dont la matrice est \emph{une
matrice de Jordan arbitraire}: en effet, avec les notations ci-dessus,
$(\sigma L)L^{-1}$ est la diagonale de blocs $\xi_{1,m_{1}}$,...,
$\xi_{k,m_{k}}$ et $(\sigma C)C^{-1}$ est la diagonale de blocs (diagonaux !)
$c_{1}I_{m_{1}}$,...,$c_{k}I_{m_{k}}$. R\'eciproquement, il est clair que la
donn\'ee de cette matrice de Jordan sp\'ecifie compl\`etement $LC$.\\

Nous appellerons solutions canoniques les solutions ainsi construites. Elles
sont au moins aussi canoniques que les bases utilis\'ees pour obtenir une
d\'ecomposition de Jordan; dans le cas semi-simple, on pourra les rigidifier
un peu plus en introduisant une relation d'ordre arbitraire entre les valeurs
propres, c'est \`a dire sur $\mathbf{C}^{*}$, voire sur 
$\mathbf{C}^{*}/q^{\mathbf{Z}}$.\\

\textbf{\underline{Remarque}}\\

Le ``Pochhammerien'' $\det(X)$ prend la forme
$\lambda \times z^{\mu_{1} + ... + \mu_{r}} \times 
\prod_{1 \leq i \leq r} e_{c_{i}} \;,\; \lambda \in \mathbf{C}^{*}$. Ceci est
compatible avec l'\'equation $\sigma(\det X) = (\det A) (\det X)$; mais ce
n'est en g\'en\'eral pas une solution canonique ! Par exemple, si $A$ est la
matrice constante $\Diag(c,d)$, on trouve $e_{c}e_{d}$ au lieu de $e_{cd}$.
Cette diff\'erence est essentielle dans le cas o\`u les $e_{c}$ sont 
r\'ealis\'es par de ``vraies fonctions'' (voir les remarques \`a la
fin de 1.1.2).



\section{Solutions canoniques dans le cas non-r\'esonnant}


\subsection*{Modification de la notion de solution canonique}

Nous reprenons ici les notations concr\`etes de notre situation, o\`u les
symboles $l$, $e_{c}$, etc... sont instanci\'es par les fonctions $l_{q}$,
$e_{q,c}$, etc...\\

On a vu ci-dessus que la construction d'une ``solution canonique''
$X = Q_{0} H L C$ m\'eromorphe dans un voisinage \'epoint\'e de $0$ dans
$\mathbf{C}^{*}$ d\'ependait tout de m\^eme du choix de la matrice de passage
$Q_{0}$ \`a la forme de Jordan, $J_{0}$. Nous allons expliciter cette 
d\'ependance.\\

On peut toujours imposer une forme normale de la r\'eduite de Jordan en 
supposant les exposants tels que $1 \leq c \leq |q|$ (en vertu de la relation
$e_{q,qc} \in \mathbf{C}^{*}ze_{q,c}$) et rang\'es selon un ordre total
arbitraire d\'efini sur $\mathbf{C}^{*}/q^{\mathbf{Z}}$; et que les blocs de 
Jordan correspondant \`a un m\^eme exposant sont rang\'es par ordre croissant
de tailles.\\

Si l'on prend une \emph{autre} matrice de passage $Q'_{0}$ \`a la \emph{m\^eme}
forme de Jordan $J_{0}$, les \'egalit\'es 
$(Q_{0})^{-1}A_{0}Q_{0} = (Q'_{0})^{-1}A_{0}Q'_{0} = J_{0}$ montrent que la 
matrice $S = (Q_{0})^{-1}Q'_{0}$ commute avec $J_{0}$, donc avec $L$ et $C$. La
solution canonique construite \`a l'aide de $Q'_{0}$ est $X' = Q'_{0} H' L C$.
Notant $(Q_{0})^{-1}AQ_{0} = J$ et $(Q'_{0})^{-1}AQ'_{0} = J'$, de sorte que
$J(0) = J'(0) = J_{0}$, on voit que $H$ et $H'$ sont respectivement 
caract\'eris\'es par les \'equations avec conditions initiales:
$$
\left \lbrace
\begin{array}{l}
\sigma_{q}H = J H (J_{0})^{-1}\\
H(0) = Id_{n}\\
\end{array}
\right.
\qquad
\left \lbrace
\begin{array}{l}
\sigma_{q}H' = J' H' (J'_{0})^{-1}\\
H'(0) = Id_{n}\\
\end{array}
\right.
$$

Le fait que ces syst\`emes avec conditions initiales \emph{caract\'erisent} 
leurs solutions respectives, c'est \`a dire que chacun admet une unique
solution formelle et que celle-ci converge, est cons\'equence de 1.2.2: en lieu
et place de l'endomorphisme $X \mapsto AX$ de $\mathbf{C}^{n}$, on l'applique
\`a l'endomorphisme $H \mapsto J H (J_{0})^{-1}$ (par exemple) de l'espace des
matrices; celui-ci admet bien $1$ comme exposant en $0$, et $Id_{n}$ comme
vecteur propre associ\'e \`a cet exposant. On peut \'egalement le d\'eduire de 
ce qui suit.\\ 

On voit alors que $H' = S^{-1} H S$ car cette derni\`ere matrice est
solution du syst\`eme qui caract\'erise $H'$. On trouve donc que $X' = X S$.\\

Posons alors $\overline{X} = X (Q_{0})^{-1} = Q_{0} H L C (Q_{0})^{-1}$ et
$\overline{X}' = X' (Q'_{0})^{-1} = Q'_{0} H' L C (Q'_{0})^{-1}$: ce sont des
solutions de l'\'equation $\sigma_{q}X = AX$, construites par un proc\'ed\'e
canonique \`a partir de deux jordanisations de $A_{0}$; et l'\'equation
$X' = X S$ dit pr\'ecis\'ement que $\overline{X} = \overline{X}'$. On a donc
ici une solution r\'eellement canonique, qui ne d\'epend que du choix d'une
forme de Jordan; et on a vu que celui-ci admettait une forme normale.\\

Le but de ce paragraphe et du suivant est de d\'efinir proprement et 
d'\'etudier cette solution canonique:

\begin{itemize}

\item{D\'efinir proprement: on va voir que l'on peut en effet reporter le choix
de la jordanisation (et la preuve que la solution construite n'en d\'epend pas)
au cas des solutions d'\'equations $\sigma_{q}X = AX$ \`a matrice constante
$A \in GL_{n}(\mathbf{C})$. L'argument ci-dessus, qui suppose que l'on a
``rigidifi\'e'' la r\'eduite de Jordan est commode, mais pas tr\`es \'el\'egant
et l'on pourra s'en passer. 
}

\item{Etudier: en vue de la seconde partie (o\`u l'on fait tendre $q$ vers $1$)
il faudra \'etablir quelques propri\'et\'es d'alg\`ebre lin\'eaire, avec ou 
sans param\`etre, pour lesquelles, pour bien connues qu'elles soient, 
je n'ai pas de r\'ef\'erence.L'\'etude des propri\'et\'es tensorielles est
report\'ee \`a la suite de ce travail (sur la th\'eorie de Galois).}

\end{itemize}



\subsection{Comment se ramener aux \'equations \`a coefficients constants}

La m\'ethode de r\'esolution choisie ici s'inspire de \cite{SVdP}, p. 154, 
o\`u le r\'esultat est toutefois formul\'e en termes plus intrins\`eques.\\

Supposons (ce sera fait \`a partir du 1.3.2) que l'on sache r\'esoudre les
\'equations $\sigma_{q}X = AX$ \`a matrice constante 
$A \in GL_{n}(\mathbf{C})$. Alors le cas g\'en\'eral se traite comme suit: on
part de l'\'equation $\sigma_{q}X = AX$, o\`u l'on suppose seulement que
$A_{0} = A(0) \in GL_{n}(\mathbf{C})$. On impose $X = F Y$ et 
$\sigma_{q}Y = A_{0}Y$. L'\'equation d\'eduite $\sigma_{q}(FY) = A(FY)$
\'equivaut \`a $\sigma_{q}F A_{0} = A F$. On montre plus loin que cette
derni\`ere \'equation admet une unique solution convergente astreinte \`a la
condition initiale $F(0) = Id_{n}$, \`a la seule condition que $A$ soit
non-r\'esonnante. On en d\'eduit la solution canonique d\'esir\'ee.\\

Cela revient donc \`a prouver l'existence d'une transformation de jauge
$X = FY$ qui \'etablisse l'\'equivalence m\'eromorphe $A \sim A_{0}$; on peut
de plus exiger que cette transformation soit tangente \`a l'identit\'e. Il faut
prendre garde qu'il ne s'agit pas l\`a d'un cas particulier de la relation
d'\'equivalence m\'eromorphe entre \'equations aux $q$-diff\'erences, car $F$
n'est ici pas m\'eromorphe sur la sph\`ere de Riemann.


\subsubsection{Transformation de jauge qui fait passer de $A$ \`a $A_{0}$}

On montre ici que l'\'equation $\sigma_{q} X = A X$
est m\'eromorphiquement (et m\^eme holomorphiquement) \'equivalente
\`a $\sigma_{q}Y = A_{0}Y$ dans le cas non-r\'esonnant. On reprend les 
notations et hypoth\`eses de 1.2.2. Nous aurons besoin du\\

\textbf{\underline{Lemme}}\\

\emph{Soit $M_{0}$ un \'el\'ement de 
$M_{n}(\mathbf{C})$. L'endomorphisme
$\Phi_{M_{0},\lambda} : M \mapsto \lambda M M_{0} - M_{0} M$ de
$M_{n}(\mathbf{C})$ a pour spectre le multi-ensemble
$\lambda Sp(M_{0}) - Sp(M_{0})$.}\\

Autrement dit, si l'on note 
$Sp(M_{0}) = \{\lambda_{1} \;,\;...\;,\; \lambda_{n} \}$, 
le spectre de $\Phi_{M_{0},\lambda}$ est 
$\{\lambda \lambda_{i} - \lambda_{j} \;/\; 1 \leq i \;,\; j \leq n\}$ en 
comptant les multiplicit\'es.\\

\textbf{\underline{Preuve}}\\

L'assertion du lemme revient \`a dire que le polynome
caract\'eristique de $\Phi_{M_{0},\lambda}$ est
$\underset{1 \leq i,j \leq n}{\prod} (T - \lambda \lambda_{i} + \lambda_{j})$. Cela 
\'equivaut \`a un syst\`eme d'identit\'es polynomiales en les coefficients de 
$M_{0}$ et il suffit donc de le v\'erifier pour les \'el\'ements d'une partie 
Zariski-dense de $M_{n}(\mathbf{C})$, par exemple pour les matrices 
semi-simples.\\

Par ailleurs, si l'on note, pour $N_{0} \in GL_{n}(\mathbf{C})$, $i_{N_{0}}$
l'automorphisme int\'erieur $M \mapsto N_{0} M (N_{0})^{-1}$ de
$M_{n}(\mathbf{C})$, on a l'\'egalit\'e :
$\Phi_{N_{0} M_{0} (N_{0})^{-1},\lambda} = 
i_{N_{0}} \circ \Phi_{M_{0},\lambda} \circ (i_{N_{0}})^{-1}$ : si $M_{0}$ est 
semblable \`a $M_{1}$, alors $\Phi_{M_{0},\lambda}$ est semblable \`a
$\Phi_{M_{1},\lambda}$. On peut donc supposer $M_{0}$ diagonale :
$M_{0} = \Diag(\lambda_{1} \;,\; ... \;,\; \lambda_{n})$. Mais, dans ce cas,
l'effet de $\Phi_{M_{0},\lambda}$ sur les matrices \'el\'ementaires est
$E_{i,j} \mapsto (\lambda \lambda_{j} - \lambda_{i})E_{i,j}$ et la conclusion
est alors imm\'ediate.\hfill $\Box$\\

Ce lemme est d'ailleurs valable sur tout corps commutatif.\\

\textbf{\underline{Th\'eor\`eme}}\\

\emph{Il existe une unique transformation
de jauge formelle , tangente \`a l'identit\'e 
$F = Id_{n} + z F_{1} + ... \in GL_{n}(\mathbf{C}\{z\})$ telle que
$(\sigma_{q}F)^{-1}AF = A_{0}$. Cette transformation est convergente.}\\

\textbf{\underline{Preuve}}\\

Encore une fois, il y a une \'etape formelle puis
le recours aux s\'eries majorantes.

\begin{itemize}

\item{Etape 1: existence et unicit\'e de la solution formelle.
L'\'equation prescrite \'equivaut au syst\`eme avec conditions initiales:
$$
\left \lbrace
\begin{array}{l}
A F = (\sigma_{q}F) A_{0}\\
F_{0} = Id_{n}\\
\end{array}
\right.
\quad \iff \quad
\left \lbrace
\begin{array}{l}
F_{0} = Id_{n}\\
\forall m \geq 1 \;,\; q^{m}F_{m}A_{0} = A_{0}F_{m} +...+ A_{m}F_{0}\\
\end{array}
\right.
$$

Le lemme ci-dessus, joint \`a l'hypoth\`ese de non-r\'esonnance, garantit que,
pour $m \geq 1$, l'endomorphisme $\Phi_{A_{0},q^{m}}$ est inversible. Les
coefficients $F_{m}$ sont donc uniquement d\'efinis par les relations
$$
\left \lbrace
\begin{array}{l}
F_{0} = Id_{n}\\
\forall m \geq 1 \;,\; 
F_{m} = (\Phi_{A_{0},q^{m}})^{-1}(A_{1}F_{m-1} +...+ A_{m}F_{0})\\
\end{array}
\right.
$$

Cette partie de la preuve garde un sens et reste valable sur un corps
commutatif quelconque.
}

\item{Etape 2: convergence des solutions formelles.
Puisque $|q| > 1$, l'endomorphisme $\Phi_{A_{0},q^{m}}$ est \'equivalent, 
lorsque $m \to \infty$, \`a la multiplication \`a droite qui \`a
$M$ associe $M \times (q^{m}A_{0})$. Choisissant, dans 
$End_{\mathbf{C}}(M_{n}(\mathbf{C}))$ la norme associ\'ee \`a une norme 
d'alg\`ebre quelconque sur $M_{n}(\mathbf{C})$, on conclut que
$\underset{m \to \infty}{\lim} \|(\Phi_{A_{0},q^{m}})^{-1}\| = 0$.. 

Il s'ensuit en particulier que
$b = \underset{m \geq 1}{\sup} \|(\Phi_{A_{0},q^{m}})^{-1}\| < \infty$. Notant
$a_{m} = \|A_{m}\|$ et $f_{m} = \|F_{m}|$, on obtient les in\'egalit\'es:
$\forall m \geq 1 \;,\; 0 \leq f_{m} \leq b(a_{1}f_{m-1} + ... +a_{m}f_{0})$.
Comme la s\'erie $\sum a_{m} z^{m}$ converge, la preuve se termine comme dans 
le lemme 1 de 1.2.2.\hfill $\Box$
}

\end{itemize}

Il est \`a noter que, si $A$ est m\'eromorphe sur $\mathbf{C}$ (dans la suite,
elle le sera m\^eme sur $\mathbf{S}$), $F$ l'est aussi: son \'equation
fonctionnelle se r\'e\'ecrit en effet $\sigma_{q}F = A F A_{0}^{-1}$ d'o\`u
la ``propagation des p\^oles'' d\'ej\`a not\'ee.



\subsection{Solution canonique d'une \'equation \`a coefficients constants}

Pour la construire, nous nous appuierons sur la d\'ecomposition de
Dunford multiplicative d'une matrice inversible.

\begin{enumerate}

\item{
\textbf{\underline{Cas d'une matrice unipotente}}\\

Soit $U = Id_{n} + N$ o\`u $N \in M_{n}(\mathbf{C})$ est nilpotente.\\

\textbf{\underline{Proposition}}\\

\emph{ L'\'equation $\sigma_{q} M = U M$ admet
une unique solution \`a coefficients dans l'anneau
$\mathbf{C}[l_{q}] = \sum_{k \geq 0} \mathbf{C} l_{q}^{(k)}$, soit 
$M = l_{q}^{(0)}M_{0} + l_{q}^{(1)}M_{1} + ...$ telle que $M_{0} = Id_{n}$. Cette 
solution est d\'efinie par $M_{k} = N^{k} \;,\; \forall k \geq 0$. Elle 
commute avec $N$ et $U$.}\\

\textbf{\underline{Preuve}}\\

Compte tenu de ce que 
$\sigma_{q} l_{q}^{(k)} = l_{q}^{(k)} + l_{q}^{(k-1)}$, on r\'esoud 
formellement:
$$
\sum_{k=0}^{d} (l_{q}^{(k)} + l_{q}^{(k-1)})M_{k} = \sigma_{q} M = 
U M = \sum_{k=0}^{d} l_{q}^{(k)}(UM_{k})
$$
Cela donne, vue la $\mathbf{C}$-ind\'ependance lin\'eaire des $l_{q}^{(k)}$,
$M_{k} + M_{k+1} = UM_{k}$, soit $M_{k+1} = N M_{k}$. Le reste s'ensuit.
\hfill $\Box$\\

L'exponentielle d'une matrice nilpotente et le logarithme d'une matrice 
unipotente \'etant d\'efinies sur un corps commutatif de caract\'eristique $0$
quelconque, cette solution vaut donc 
$$
U^{l_{q}} = \exp(l_{q} \log(Id_{n} + N)) = 
\sum_{k \geq 0} \begin{pmatrix} l_{q}\\ k\\ \end{pmatrix} N^{k}
$$
On la notera $e_{q,U}$. Par construction, elle est unipotente. A titre 
d'exemple, si l'on prend $U = \xi_{1,m}$, on retrouve la matrice $L_{m}$ 
de 1.2.4.\\

\textbf{\underline{Corollaire 1}}\\

\emph{Soit $R \in GL_{n}(\mathbf{C})$ et soit
$V = R U R^{-1}$. Alors $e_{q,V} = R e_{q,U} R^{-1}$.}\\

\textbf{\underline{Preuve}}\\

$R e_{q,U} R^{-1}$ satisfait les conditions qui
sp\'ecifient $e_{q,V}$ ! 
\hfill $\Box$\\

\textbf{\underline{Corollaire 2}}\\

\emph{Une matrice commute avec $e_{q,U}$ si et
seulement si elle commute avec $U$.}\\

\textbf{\underline{Preuve}}\\

Les \'egalit\'es $e_{q,U} = \exp(l_{q} \log(U))$ et
$U = \exp(\log(e_{q,U})/l_{q})$ montrent que les matrices $U$ et $e_{q,U}$ sont
chacune un polynome en l'autre, \`a coefficients dans $\mathbf{C}(l_{q})$.
\hfill $\Box$\\
}

\item{
\textbf{\underline{Cas d'une matrice semi-simple}}\\

Soit $D = \Diag(c_{1},...,c_{n}) \in GL_{n}(\mathbf{C})$ une matrice diagonale
inversible;
on note alors $e_{q,D} = \Diag(e_{q,c_{1}},...,e_{q,c_{n}})$. On a donc
$\sigma_{q} e_{q,D} = D e_{q,D} = e_{q,D} D$. De plus, le fait que la correspondance
$c \leftrightarrow e_{q,c}$ est biunivoque montre que $D$ est un polynome en 
$e_{q,D}$ et r\'eciproquement (par exemple, polynomes \`a coefficients dans
$\mathcal{M}(\mathbf{C}^{*})$) et donc qu'ils ont m\^eme commutant.\\

{\small
\textbf{\underline{Disgression sur les conjugaisons entre matrices diagonales}}\\

Si $\Delta = \Diag(c_{1},...,c_{n}) \in M_{n}(\mathbf{C})$ et
$\Delta' = \Diag(c'_{1},...,c'_{n'}) \in M_{n'}(\mathbf{C})$, la relation
matricielle $R\Delta = \Delta'R$ ($R \in M_{n',n}(\mathbf{C})$) est ``de nature
combinatoire''. En effet, elle s'\'ecrit:
$$
\forall i \in \{1,...,n'\} \;,\; \forall j \in \{1,...,n\} \;,\; 
r_{i,j}c_{j} = c'_{i}r_{i,j}
$$
autrement dit:
$$
c_{j} \not= c'_{i} \Rightarrow r_{i,j} = 0
$$
Seules les positions de $R$
``connectant des positions de valeurs \'egales'' de $\Delta$ et $\Delta'$ ont 
le droit de
porter une valeur non nulle, et c'est leur seule contrainte (la valeur non 
nulle port\'ee n'importe pas).\\

Soit $f : \mathbf{C} \rightarrow \mathbf{C}$ une application quelconque. Notons
(abusivement) $f$ son extension de aux matrices diagonales (application de $f$
aux coefficients). Alors, on a la relation \'evidente mais tr\`es utile:
$$
R\Delta = \Delta'R \Rightarrow R f(\Delta) = f(\Delta') R
$$
Cela se prouve par les implications suivantes:
$r_{i,j} \not= 0 \Rightarrow c_{j} = c'_{i} \Rightarrow f(c_{j}) = f(c'_{i})$
Si $f$ est injective, on a m\^eme une \'equivalence logique.\\

Ainsi, si une matrice $D$ admet les deux diagonalisations
$D = R \Delta R^{-1} = R' \Delta' (R')^{-1}$, on a l'\'egalit\'e
$R f(\Delta) R^{-1} = R' f(\Delta') (R')^{-1}$ de sorte que l'on peut poser
$f(D) = R f(\Delta) R^{-1}$ et que cela ne d\'epend pas de la diagonalisation
particuli\`ere choisie. Cela permet d'\'etendre canoniquement aux matrices
semi-simples toute application $f : \mathbf{C} \rightarrow \mathbf{C}$. On peut
d'ailleurs d\'efinir de la m\^eme fa\c{c}on $g(D_{1},D_{2})$ pour
$g : \mathbf{C} \times \mathbf{C} \rightarrow \mathbf{C}$, etc...\\
}

\textbf{\underline{Proposition et d\'efinition}}\\

\emph{Soit $D$ une matrice
diagonalisable. Ecrivons 
$$
D = R \Delta R^{-1} = R' \Delta' (R')^{-1}
$$ 
o\`u $R \;,\; R' \in GL_{n}(\mathbf{C})$ et o\`u $\Delta$ et $\Delta'$ sont
diagonales. Alors 
$$
R e_{q,\Delta} R^{-1} = R' e_{q,\Delta} (R')^{-1}
$$
On note cette matrice $e_{q,D}$.}\\

\textbf{\underline{Preuve}}\\

Cela d\'ecoule imm\'ediatement de ce qui 
pr\'ec\`ede.
\hfill $\Box$\\

On tire facilement des m\^emes remarques les\\

\textbf{\underline{Propri\'et\'es}}\\

\emph{Soit $D \in GL_{n}(\mathbf{C})$ 
diagonalisable:\\
(i) $\sigma_{q}e_{q,D} = D e_{q,D} = e_{q,D} D$ \\
(ii) Si $D' = R D R^{-1}$, $e_{q,D'} = R e_{q,D} R^{-1}$\\
(iii) Chacune des matrices $D$ et $e_{q,D}$ est un polynome en l'autre \`a
coefficients dans $\mathcal{M}(\mathbf{C}^{*})$ et elles ont donc m\^eme
commutant.}\\
}

\item{
\textbf{\underline{Cas g\'en\'eral}}\\

Soit $A = D + N$ la d\'ecomposition de Dunford d'une matrice 
$A \in M_{n}(\mathbf{C})$: $D$ est semi-simple, $N$ est nilpotente et 
$[D,N] = 0$ (les matrices $D$ et $N$ commutent). Ces conditions sp\'ecifient 
uniquement $D$ et $N$. Si $A$ est inversible, on en d\'eduit la
\emph{d\'ecomposition de Dunford multiplicative} $A = DU$, o\`u $D$ est 
semi-simple, $U$ est unipotente et $[D,U] = 0$ (les matrices $D$ et $U$ 
commutent): il suffit de prendre $U = Id_{n} + D^{-1}N$. Il y a encore 
unicit\'e sous ces conditions. Si $R \in GL_{n}(\mathbf{C})$, la 
d\'ecomposition de Dunford multiplicative de $RAR^{-1}$ est alors
$(RDR^{-1})(RUR^{-1})$.\\

\textbf{\underline{D\'efinition}}\\

On appellera \emph{solution canonique dans le cas non-r\'esonnant} la matrice
$e_{q,A} = e_{q,D} e_{q,U}$\\

Ce qui suit est alors imm\'ediat:\\

\textbf{\underline{Propri\'et\'es}}\\

\emph{
(i) $\sigma_{q} e_{q,A} = A e_{q,A} = e_{q,A} A$\\
(ii) La d\'ecomposition de Dunford multiplicative de $e_{q,A}$ est 
$e_{q,D} e_{q,U}$\\
(iii) $[M,A] = 0 \Leftrightarrow [M,D] = [M,U] = 0 \Leftrightarrow
[M,e_{q,D}] = [M,e_{q,U}] = 0 \Leftrightarrow [M,e_{q,A}] = 0$\\
(iv) $e_{q,A}$ est m\'eromorphe sur $\mathbf{C}^{*}$}\\

Pour la premi\`ere et la troisi\`eme \'equivalence de (iii), voir Bourbaki, 
Alg\`ebre Lin\'eaire ou 1.4.
\hfill $\Box$
}

\end{enumerate}


\chapter{Th\'eorie de Fuchs-Riemann-Birkhoff:\\connexion des solutions locales
         sur $\mathbf{P}^{1}\mathbf{C}$}



\section{Equations fuchsiennes sur la sph\`ere de Riemann}

On se place d\'esormais sur la sph\`ere de Riemann 
$\mathbf{S} = \mathbf{P}^{1}\mathbf{C}$. Le corps de base est donc
$\mathcal{M}(\mathbf{S}) = \mathbf{C}(z)$: on le notera $\mathbf{k}$. On se
restreint de plus \`a l'\'etude des \'equations aux $q$-diff\'erences 
lin\'eaires \`a coefficients m\'eromorphes sur $\mathbf{S}$. Ce sont donc les 
\'equations $\sigma_{q}X = AX \;,\; A \in GL_{n}(\mathbf{k})$.\\

Nous consid\'ererons plus pr\'ecis\'ement les \'equations \emph{singuli\`eres 
r\'eguli\`eres} en $0$ et en $\infty$. Cela signifie que $A$ est \'equivalente
\`a une \'equation fuchsienne en $0$ via une transformation de jauge
m\'eromorphe sur $\mathbf{C}$, et, simultan\'ement, \`a une \'equation 
fuchsienne en $\infty$ via une transformation de jauge m\'eromorphe sur 
$\mathbf{S} - \{\infty\}$. Nous prouverons plus loin que $A$ est alors en fait
\'equivalente \`a une \'equation fuchsienne en $0$ et en $\infty$ via une 
transformation de jauge m\'eromorphe sur $\mathbf{S}$, c'est \`a dire
rationnelle.\\



\subsection{Equations singuli\`eres r\'eguli\`eres en $0$: formes normales des 
solutions}

Le d\'eveloppement des fractions rationnelles en s\'eries de Laurent d\'efinit
le plongement 
$\mathbf{k} := \mathbf{C}(z) \hookrightarrow 
\mathbf{k}'_{0} := \mathbf{C}(\{z\}) = \mathbf{C}\{z\}[z^{-1}]$
et permet de consid\'erer toute \'equation aux $q$-diff\'erences lin\'eaires
\`a coefficients rationnels 
$\sigma_{q} X = A X$ ($A \in GL_{n}(\mathbf{C}(z))$) comme une \'equation au
voisinage de $0$, c'est \`a dire \`a coefficients m\'eromorphes au voisinage
de $0$. Dans le cas fuchsien, on peut donc appliquer les r\'esultats du 
chapitre 1. On r\'ealisera les solutions comme des fonctions en interpr\'etant
comme suit les symboles $e_{c}$ et $l$:

\begin{itemize}

\item{$l = l_{q}$ est la fonction $z \Theta_{q}'/\Theta_{q}$ introduite au
1.1.3; c'est donc un \'el\'ement du sous-corps
$\mathbf{C}(z , \Theta_{q} , \Theta_{q}')$ de $\mathcal{M}(\mathbf{C}^{*})$.}

\item{Les $e_{c}$ sont une famille de caract\`eres comme on en a construit au
1.1.2; on aura notamment besoin des propri\'et\'es suivantes:\\

(i) Tous les ``caract\`eres'' $e_{c}$ sont \'el\'ements de
$\mathbf{C}(z)((\Theta_{q,a})_{a \in \mathbf{C}^{*}}) \subset
\mathcal{M}(\mathbf{C}^{*})$.\\
(ii) $e_{1}$ et tous les $e_{qc}/ze_{c}$ (qui sont \emph{a priori} elliptiques)
sont des vraies constantes, c'est \`a dire \'el\'ements de $\mathbf{C}^{*}$.\\
(iii) Pour $c \not \in q^{\mathbf{Z}}$, le diviseur des z\'eros et des p\^oles
de $e_{c}$ dans $\mathbf{C}^{*}$ est une spirale logarithmique discr\`ete de
la forme $q^{\mathbf{Z}}\alpha$.\\

La propri\'et\'e (ii) va servir \`a normaliser les solutions d'\'equations aux
$q$-diff\'erences sous une forme suffisamment rigide. La propri\'et\'e (iii),
que l'on pourrait d'ailleurs \'elargir, permettra aux chapitres 3 et 4 de
donner une allure agr\'eable \`a l'\'etude des p\^oles des solutions.}

\end{itemize}

Toutes ces conditions sont en particulier r\'eunies dans le cas de la famille 
``standard'' des $e_{q,c} = \Theta_{q}/\Theta_{q,c}$, que l'on utilisera donc, 
sauf 
mention expresse du contraire.\\

Nous allons d'abord r\'ecapituler les r\'esultats de 1.3.2. sous les 
hypoth\`eses de ce chapitre.\\

\begin{itemize}

\item{\underline{Cas fuchsien non-r\'esonnant:}\\

On a montr\'e l'existence d'une solution fondamentale de la
forme $X^{(0)} = Q^{(0)}H^{(0)}L^{(0)}C^{(0)}$ (voir 1.2.2), o\`u 
$Q^{(0)} \in GL_{n}(\mathbf{C})$ est une matrice de passage de $A(0)$ \`a sa
r\'eduite de Jordan (sous la forme que nous avons adopt\'ee, form\'ee de blocs
$c\xi_{1,m}$) et o\`u le choix de $Q^{(0)}$ d\'etermine totalement la forme
obtenue.\\

D'autre part, $A$ \'etant maintenant rationnelle, donc m\'eromorphe sur 
$\mathbf{C}$, l'\'equation fonctionnelle qui caract\'erise $H^{(0)}$ montre
que celle-ci est \'egalement m\'eromorphe sur $\mathbf{C}$: ce point sera
pr\'ecis\'e au 2.1.2.
}\\

\item{\underline{Cas fuchsien g\'en\'eral:}\\

L'algorithme de pr\'eparation de 1.2.3 montre qu'il y a alors \'equivalence 
rationnelle avec le cas non-r\'esonnant. On obtient donc dans ce cas une 
solution de la forme $X^{(0)} = U^{(0)}Q^{(0)}H^{(0)}L^{(0)}C^{(0)}$ avec
$U^{(0)} \in GL_{n}(\mathbf{k})$, les autres facteurs \'etant les m\^emes que
ci-dessus.
}\\

\item{\underline{Cas singulier r\'egulier:}\\

Par d\'efinition, il y a \`a nouveau \'equivalence (m\'eromorphe sur 
$\mathbf{C}$) avec le cas pr\'ec\'edent. On obtient donc une solution de la 
forme $X^{(0)} = V^{(0)}U^{(0)}Q^{(0)}H^{(0)}L^{(0)}C^{(0)}$ avec
$V^{(0)} \in GL_{n}(\mathbf{\mathcal{M}(\mathbf{C})})$, les autres facteurs 
\'etant les m\^emes que
ci-dessus.
}\\

\end{itemize}

Dans tous les cas, on peut donc \'ecrire $X^{(0)} = M^{(0)} N^{(0)}$, o\`u
$M^{(0)} \in GL_{n}(\mathcal{M}(\mathbf{C}))$ (voir 2.1.2) et o\`u $N^{(0)}$ 
est une matrice log-car.\\

Rappelons \'egalement que, dans le cas fuchsien non-r\'esonnant, nous avons
obtenu une autre construction de solution, compl\`etement canonique, sous la
forme $M^{(0)} e_{q,A^{(0)}}$; la partie $M^{(0)}$ est encore m\'eromorphe sur
$\mathbf{C}$ (cf 2.1.2).\\

\textbf{\underline{Normalisation}}\\

Du fait que $e_{qc} \equiv z e_{c} \pmod{\mathbf{C}^{*}}$, on peut imposer \`a
la matrice $C^{(0)}$ (diagonale de caract\`eres) de ne contenir que des
caract\`eres d'exposants $c$ tels que $1 \leq |c| < |q|$ et de commuter
tout de m\^eme avec $L$ (voir les remarques de 1.3.2, point 2). De m\^eme, on
peut introduire un ordre arbitraire sur $\mathbf{C}^{*}/q^{\mathbf{Z}}$ et
ranger les caract\`eres en ordre croissant d'exposants; et, pour chaque
caract\`ere, les blocs unipotents associ\'es par ordre croissant de tailles:
ces deux conditions s'obtiennent en effet \`a l'aide de matrices de
permutations. On dira, si toutes ces conditions sont r\'eunies, que la solution
est sous forme normalis\'ee. Il n'y a pas unicit\'e de solutions de cette 
forme: on verra au 2.1.3 ce qu'il en est exactement. En revanche, \`a solution 
donn\'ee, on prouvera que l'\'ecriture normalis\'ee est unique.\\



\subsection{Polarit\'e des solutions fondamentales en $0$}

A strictement parler, l'\'etude locale du chapitre 1 d\'efinit $M^{(0)}$ comme
un germe de fonction m\'eromorphe sur un voisinage de $0$ dans $\mathbf{C}$,
\`a valeurs dans $GL_{n}(\mathbf{C})$ (alors que $N^{(0)}$ est m\'eromorphe
sur $\mathbf{C}^{*}$).\\

On \'etudiera la polarit\'e des solutions en comptant comme singularit\'e d'une
fonction $F$ \`a valeurs dans $GL_{n}(\mathbf{C})$ tout point o\`u elle n'est 
pas d\'efinie (donc, pour une fonction $F$ m\'eromorphe, tout p\^ole de $F$); 
et aussi tout point o\`u elle prend pour valeur une matrice non inversible, 
i.e. de d\'eterminant nul (donc tout z\'ero de $\det F$). Cela revient en 
r\'ealit\'e \`a consid\'erer les p\^oles de la fonction 
$z \mapsto (F(z),(F(z))^{-1})$ \`a
valeurs dans l'espace affine $M_{n}(\mathbf{C}) \times M_{n}(\mathbf{C})$. Ce
point de vue est justifi\'e par notre but, qui est d'\'etudier le lien entre
solutions en $0$ et solutions en $\infty$, ce qui fait naturellement apparaitre
les inverses des fonctions matricielles apparaissant dans les \'equations et
dans les solutions (voir 2.2). Nous introduisons donc la\\

\textbf{\underline{Notation}}\\

\emph{On d\'efinit le lieu singulier de la fonction matricielle $F$ comme:}
$$\mathcal{S}(F) := 
\{ \text{p\^oles de } F \} \cup \{ \text{z\'eros de } \det F \} =
\{ \text{p\^oles de } F \} \cup \{ \text{p\^oles de } F^{-1} \}
$$

\textbf{\underline{Th\'eor\`eme}}\\

\emph{
(i) $M^{(0)}$ admet un unique prolongement m\'eromorphe \`a $\mathbf{C}$ tel
que le prolongement correspondant de $X^{(0)} = M^{(0)}N^{(0)}$ reste 
solution fondamentale.\\
(ii) $X^{(0)}$ est m\'eromorphe sur $\mathbf{C}^{*}$.\\
(ii) Les singularit\'es de $M^{(0)}$ sur $\mathbf{C}^{*}$ forment une union 
finie de demi-spirales logarithmiques discr\`etes: 
$q^{\mathbf{N}^{*}}\mathcal{S}(A)$.
}\\

\textbf{\underline{Preuve}}\\

Soit $K^{(0)} = (\sigma_{q} N^{(0)}) (N^{(0)})^{-1} \in GL_{n}(\mathbf{C})$.
Alors, $M^{(0)}$ satisfait l'\'equation fonctionnelle
$\sigma_{q}M^{(0)} = A M^{(0)} (K^{(0)})^{-1}$. On r\'e\'ecrit celle-ce sous
la forme \'equivalente
$$
M^{(0)} = (\tau_{q}A) (\tau_{q}M^{(0)}) (K^{(0)})^{-1}
$$
que l'on va plut\^ot
consid\'erer comme une d\'efinition r\'ecursive; le cas de terminaison \'etant
celui ou la variable est dans le voisinage de $0$ dans $\mathbf{C}$ o\`u le
\emph{germe} $M^{(0)}$ est \emph{a priori} d\'efini. Notons donc $V$ un disque
ouvert de centre $0$ sur lequel le germe $M^{(0)}$ est holomorphe, sauf 
peut-\^etre en $0$.\\

En it\'erant $r$ fois l'\'equation ci-dessus, on trouve
$$
M^{(0)} = 
(\tau_{q}A) ... (\tau_{q}^{r}A) (\tau_{q}^{r}M^{(0)}) (K^{(0)})^{-r}
$$
ce qui conduit \`a \emph{d\'efinir} la fonction
$(\tau_{q}A) ... (\tau_{q}^{r}A) (\tau_{q}^{r}M^{(0)}) (K^{(0)})^{-r}$
sur $\sigma^{r}V$: elle y est m\'eromorphe avec pour seules singularit\'es 
possibles dans $V^{*}$ celles de $(\tau_{q}A) , ... , (\tau_{q}^{r}A)$, c'est
\`a dire les \'el\'ements de 
$\underset{1 \leq i \leq r}{\bigcup} \mathcal{S}(A)$.\\

Les d\'efinitions de ces fonctions sur la suite croissante des $\sigma^{r}V$
sont compatibles en vertu de l'\'equation fonctionnelle et la r\'eunion de ces
ouverts est $\mathbf{C}$ puisque $|q| > 1$; ceci ach\`eve la d\'emonstration.
\hfill $\Box$\\

Dans le cas fuchsien non r\'esonnant, on sait de plus que $M^{(0)}$ est
holomorphe en $0$.\\

\textbf{\underline{Corollaire}}\\

\emph{Outre les demi-spirales logarithmiques discr\`etes engendr\'ees par les
singularit\'es de $A$, les singularit\'es de la solution fondamentale
$X^{(0)} = M^{(0)} N^{(0)}$ sur $\mathbf{C}^{*}$ forment des spirales
logarithmiques discr\`etes d\'etermin\'ees par la structure de Jordan:
\begin{itemize}
\item{Les spirales $q^{\mathbf{Z}}\alpha$ de z\'eros et de p\^oles des
caract\`eres non triviaux, c'est \`a dire d'exposants
$c \in Sp(A(0)) - q^{\mathbf{Z}}$.}
\item{La spirale $q^{\mathbf{Z}}$ des p\^oles de $l_{q}$ si, apr\`es 
pr\'eparation, $A(0)$ a au moins un bloc unipotent de taille $\geq 2$.}
\end{itemize}
}

Notons que les blocs log ont tous un d\'eterminant \'egal \`a $1$ et ne peuvent
donc faire sortir les valeurs $X^{(0)}$ de $GL_{n}(\mathbf{C})$.



\subsection{Non-unicit\'e des solutions canoniques}

On se donne encore une \'equation $\sigma_{q}X = AX$ \`a coefficients
rationnels, singuli\`ere r\'eguli\`ere en $0$.\\

\textbf{\underline{Th\'eor\`eme}}\\

\emph{Soit $X = MN$ une solution canonique fondamentale en $0$ en forme 
normale. Alors $X' = M'N'$ est une solution canonique fondamentale en $0$ en 
forme normale si et seulement si $N' = N$ et $M' = MR$ o\`u 
$R \in GL_{n}(\mathbf{C})$ et $R$ commute avec $N$.}\\

\newpage

\textbf{\underline{Preuve}}\\

Seule la n\'ecessit\'e de la condition n'est pas \'evidente!!\\

Les matrices $K = (\sigma_{q}N) N^{-1}$ et $K' = (\sigma_{q}N') (N')^{-1}$ sont 
des \'el\'ements de $GL_{n}(\mathbf{C})$ ayant respectivement m\^eme structure
de Jordan (tailles de blocs et positions de valeurs propres \'egales ou 
distinctes: voir la disgression sur les matrices diagonales de 1.3) et donc
aussi m\^eme commutant que $N$ et $N'$. Par ailleurs, 
l'\'equation fonctionnelle entraine les \'egalit\'es 
$\sigma_{q}M = A M K^{-1}$ et $\sigma_{q}M' = A M' (K')^{-1}$.\\

Notons, $R = M^{-1}M' \in GL_{n}(\mathcal{M}(\mathbf{C}))$. On a alors
$N^{-1}RN' = X^{-1}X'$. Cette derni\`ere est 
elliptique puisque $X$ et $X'$ sont solutions d'une m\^eme \'equation 
lin\'eaire aux $q$-diff\'erences.\\ 

Soient $e_{c_{1}}L_{1},...,e_{c_{r}}L_{r}$ 
(resp. $e_{c'_{1}}L'_{1},...,e_{c'_{s}}L'_{s}$) les blocs diagonaux 
(carr\'es) de $N$ (resp. de $N'$), et soient 
$R_{i,j}\;,\; 1 \leq i \leq r \;,\; 1 \leq j \leq s$ les blocs
correspondants (rectangulaires) de $R$. Les blocs de $N^{-1}RN'$ sont les
$(e_{c'_{j}}/e_{c_{i}})L_{i}^{-1}R_{i,j}L'_{j}$, et ils sont elliptiques:
les coefficients de $L_{i}^{-1}R_{i,j}L'_{j}$ sont donc des caract\`eres
d'exposant $c_{i}/c'_{j}$. Mais ils sont \'el\'ements de 
$\mathbf{k}_{0}(l_{q})$, donc de la forme $a z^{p} \;,\; p \in \mathbf{Z}$
(voir 1.1.5). Les conditions de normalisation entrainent alors soit que 
$R_{i,j} = 0$,
soit que $c'_{j} = c_{i}$ et $L_{i}^{-1}R_{i,j}L'_{j}$ est \`a coefficients
dans $\mathbf{C})$.\\

Dans ce dernier cas, notant 
$S_{i,j} = L_{i}^{-1}R_{i,j}L'_{j}$, on d\'eveloppe l'\'egalit\'e 
$L_{i} S_{i,j} = R_{i,j} L'_{j}$ dans la base des $l_{q}^{(k)}$. Des
\'ecritures 
$L_{i} = \sum_{k \geq 0} l_{q}^{(k)} \xi_{0,m_{i}}^{k}$ et
$L'_{i} = \sum_{k' \geq 0} l_{q}^{(k')} \xi_{0,m'_{i}}^{k'}$, on d\'eduit que
$R_{i,j} = S_{i,j}$ puis que 
$\xi_{0,m_{i}} R_{i,j} = R_{i,j}\xi_{0,m'_{i}}^{k'}$. Mais ceci signifie que
$R \in GL_{n}(\mathbf{C})$ et que $R^{-1} K R = K'$. Les conditions impos\'ees
par la normalisation entrainent alors que $K = K'$.
\hfill $\Box$\\

\textbf{\underline{Corollaire}}\\

\emph{L'\'ecriture sous forme normale d'une solution canonique donn\'ee est
unique.}\\

\textbf{\underline{Preuve}}\\

C'est en effet le cas $MN = M'N'$, c'est \`a dire $R = I_{n}$, du th\'eor\`eme.
\hfill $\Box$



\subsection{Solutions \`a l'infini et \'equations singuli\`eres 
            r\'eguli\`eres sur $\mathbf{S}$}

Posons $w = \frac{1}{z}$. Si la fonction (vectorielle ou matricielle) $X$
v\'erifie l'\'equation fonctionnelle $X(qz) = A(z) X(z)$, alors la fonction
$\overline{X}(w) := X(z)$ v\'erifie l'\'equation aux $q$-diff\'erences
$\overline{X}(qw) = (A(1/qw))^{-1}\overline{X}(w)$. La l\'eg\`ere dissym\'etrie
vient du fait que, sous l'effet de $\sigma_{q}$, $z \leftarrow qz$ et donc
$w \leftarrow q^{-1}w$.\\

Notons donc $\overline{A}(w) = (A(1/qw))^{-1}$. Ainsi,
$\overline{A}$ est fuchsienne en $w = 0$ si et seulement si $A$ l'est en 
$\infty$. De m\^eme, $\overline{A}$ est singuli\`ere r\'eguli\`ere en $w = 0$ 
si et seulement si $A$ l'est en $\infty$: car l'\'egalit\'e 
$(\sigma_{q}U)^{-1}AU = B$ \'equivaut \`a l'\'egalit\'e
$(\sigma_{q}\overline{U})^{-1}\overline{A}\overline{U} = \overline{B}$, o\`u
$\overline{U}(w) = U(1/w)$ et $\overline{B}(w) = (B(1/qw))^{-1}$.\\

On peut donc reprendre les r\'esultats pr\'ec\'edents:\\

\textbf{\underline{Th\'eor\`eme}}\\

\emph{Si $A$ est singuli\`ere r\'eguli\`ere sur 
$\mathbf{S}$,
elle admet des solutions canoniques $X^{(0)} = M^{(0)} N^{(0)}$ et
$X^{(\infty)} = M^{(\infty)} N^{(\infty)}$ telles que:\\
(i) $M^{(0)}$ est m\'eromorphe sur $\mathbf{C}$ et $M^{(\infty)}$ est
m\'eromorphe sur $\mathbf{S} - \{0\}$. Le lieu singulier de $M^{(0)}$ est
$q^{\mathbf{N}^{*}}\mathcal{S}(A)$ et celui de $M^{(\infty)}$, en tant que 
fonction de $z$, est $q^{\mathbf{-N}}\mathcal{S}(A)$.\\
(ii) $N^{(0)}$ et $N^{(\infty)}$ sont des matrices log-car, respectivement en
$z$ et en $1/z$. Leurs lieux singuliers sont les spirales logarithmiques
discr\`etes de forme $q^{\mathbf{Z}}\alpha$ engendr\'ees par les valeurs 
propres de matrices fuchsiennes non-r\'esonnantes m\'eromorphiquement
\'equivalentes \`a $A$ (respectivement sur $\mathbf{C}$ et sur 
$\mathbf{S} - \{0\}$).}\\

Cet \'enonc\'e peut \^etre pr\'ecis\'e de mani\`ere \'evidente pour une
\'equation fuchsienne non-r\'esonnante en $0$ et en $\infty$.



\subsection{Double caract\'erisation des \'equations singuli\`eres 
r\'eguli\`eres}

On va montrer que la forme des solutions caract\'erise les \'equations
singuli\`eres r\'eguli\`eres en $0$ (resp. en $\infty$), puis donner une 
propri\'et\'e de r\'eduction simultan\'ee en $0$ et en $\infty$ des \'equations
singuli\`eres r\'eguli\`eres sur $\mathbf{S}$.\\

Appelons \emph{uniformisante} un \'el\'ement $u \in \mathcal{M}(\mathbf{C})$ 
dont $0$ est z\'ero simple, autrement dit dont la valuation $v_{0}(u)$ en $0$
vaut $1$. Nous \'etablirons d'abord un lemme tr\`es utile. Il s'agit en fait 
d'un algorithme de r\'eduction analogue \`a l'algorithme du pivot de Gauss,
valable dans tout anneau de valuation discr\`ete.\\

\textbf{\underline{Lemme}}\\

\emph{Toute matrice $M \in GL_{n}(\mathcal{M}(\mathbf{C}))$ peut s'\'ecrire
sous la forme $M = CR$, avec des matrices 
$C , R \in GL_{n}(\mathcal{M}(\mathbf{C}))$ telles que:\\
(i) $R$ est r\'eguli\`ere en $0$ (i.e. $R(0) \in GL_{n}(\mathbf{C})$).\\
(ii) $C$ est produit: 
\begin{enumerate}
\item{d'une puissance $u^{k} \;,\; -k \in \mathbf{N}$ d'une uniformisante 
$u$ arbitraire}
\item{de matrices $T_{i,\underline{\alpha}}$ d'op\'erations \'el\'ementaires 
sur les lignes}
\item{de matrices de dilatation $D_{i,v}$, o\`u les 
uniformisantes $v$ sont arbitraires (et choisies ind\'ependamment les unes des 
autres).}
\end{enumerate}
}

Pour \^etre pr\'ecis, les matrices mentionn\'ees sont d\'efinies comme suit. La
matrice $T_{i,\underline{\alpha}}$ ($\underline{\alpha} \in \mathbf{C}^{n}$,
$\alpha_{i} \not= 0$) est \'egale \`a l'identit\'e, sauf en ce qui concerne la 
ligne de rang $i$:
$$
T_{i,\underline{\alpha}} =
\begin{pmatrix}
1 & 0 & ... & 0 & ... & 0 \\
0 & 1 & ... & 0 & ... & 0 \\
... & ... & ... & ... & ... & ... \\
\alpha_{1} & \alpha_{2} & ... & \alpha_{i} & ... & \alpha_{n} \\
... & ... & ... & ... & ... & ... \\
0 & 0 & ... & 0 & ... & 1 \\
\end{pmatrix}
$$
L'inverse d'une telle matrice est une matrice de m\^eme forme.
La matrice de dilatation $D_{i,v}$ est diagonale et est \'egale \`a 
l'identit\'e, sauf en ce qui concerne le coefficient de rang $i$:
$$
D_{i,v} =
\begin{pmatrix}
1 & 0 & ... & 0 & ... & 0 \\
0 & 1 & ... & 0 & ... & 0 \\
... & ... & ... & ... & ... & ... \\
0 & 0 & ... & v & ... & 0 \\
... & ... & ... & ... & ... & ... \\
0 & 0 & ... & 0 & ... & 1 \\
\end{pmatrix}
$$

\textbf{\underline{Preuve}}\\

On commence par multiplier $M$ par une puissance $u^{k}$ d'uniformisante, de
mani\`ere \`a la rendre holomorphe en $0$: cela revient \`a se ramener \`a ce
cas, et l'on doit alors prouver le lemme avec $k = 0$ dans la conclusion. Puis
l'on fait une r\'ecurrence sur $v_{0}(\det M)$. Si cette valuation vaut $0$, $M$
est r\'eguli\`ere en $0$ et l'on prend $R = M$ et $C = I_{n}$ (produit vide).\\

Si $v_{0}(\det M) > 0$, la matrice $M(0)$ est singuli\`ere et il existe donc
$\underline{\alpha} \in \mathbf{C}^{n} - \{0\}$ tel que 
$\underline{\alpha} M(0) = 0 \in \mathbf{C}^{n}$. Soit $i$ un indice tel que
$\alpha_{i} \not= 0$ et soit $v$ une uniformisante quelconque. Alors
$M_{1} = D_{i,v}^{-1} T_{i,\underline{\alpha}} M$ est holomorphe en $0$ et
$\det M_{1} = \frac{\alpha_{i}}{v} \det M$; $M_{1}$ est meilleure que $M$:
$v_{0}(\det M_{1}) = v_{0}(\det M) - 1$, ce qui ach\`eve la r\'ecurrence.
\hfill $\Box$\\

Si par exemple, dans ce lemme, on prend syst\'ematiquement l'uniformisante $z$,
on obtient $C$ \`a coefficients dans $\mathbf{C}[z,z^{-1}]$; si l'on prend 
plut\^ot $\frac{z}{1+z}$, on obtient $C \in GL_{n}(\mathbf{C}(z))$ et de plus
r\'eguli\`ere en $\infty$.\\

\textbf{\underline{Th\'eor\`eme}}\\

\emph{Si l'\'equation aux $q$-diff\'erences \`a coefficients rationnels
$\sigma_{q}X = A X$ admet une solution fondamentale de la forme $MN$, o\`u
$M$ est m\'eromorphe sur $\mathbf{C}$ et o\`u $N$ est une matrice log-car,
alors c'est une \'equation singuli\`ere r\'eguli\`ere en $0$.}\\

\newpage

\textbf{\underline{Preuve}}\\

Appliquons ce qui pr\'ec\`ede \`a $M$, que l'on \'ecrit donc $M = CR$. Alors 
$Y = RN$ est solution de l'\'equation
$\sigma_{q}Y = BY$, avec $B = (\sigma_{q}C)^{-1}AC$. La matrice $B$ est 
rationnellement \'equivalente \`a $A$, puisque $C \in GL_{n}(\mathbf{C}(z))$.
Mais $B = (\sigma_{q}(RN))(RN)^{-1} = \sigma_{q}R (\sigma_{q}N N^{-1}) R^{-1}$
est le produit de trois matrices r\'eguli\`eres en $0$, donc l'est elle-m\^eme;
l'\'equation correspondante est donc fuchsienne en $0$.
\hfill $\Box$\\

En appliquant la m\^eme d\'ecomposition \`a des solutions canoniques en forme
normale en $0$ et en $\infty$ \`a la fois, on prouve alors le\\

\textbf{\underline{Th\'eor\`eme}}\\

\emph{Une \'equation singuli\`ere r\'eguli\`ere est \'equivalente, via une
transformation de jauge rationnelle, \`a une
\'equation fuchsienne en $0$ et $\infty$ \`a la fois.}\\

\textbf{\underline{Preuve}}\\

Soient $X^{(0)} = M^{(0)} N^{(0)}$ et
$X^{(\infty)} = M^{(\infty)} N^{(\infty)}$ des solutions canoniques en formes
normales en $0$ et en $\infty$ respectivement. Appliquant le lemme, on \'ecrit
$M^{(0)} = C^{(0)}R^{(0)}$ et $M^{(\infty)} = C^{(\infty)}R^{(\infty)}$. Il
suffit maintenant de trouver $U \in GL_{n}(\mathbf{C}(z))$ telle que $UC^{(0)}$
et $UC^{(\infty)}$ soient r\'eguli\`eres en $0$ et en $\infty$ respectivement:
on pourra alors conclure comme dans le th\'eor\`eme pr\'ec\'edent.\\

On applique pour cela le lemme \`a 
$M = (C^{(\infty)})^{-1} C^{(0)} \in GL_{n}(\mathbf{C}(z))$, que l'on \'ecrit 
$CR$, avec pour choix syst\'ematique des uniformisantes $\frac{z}{1+z}$, de
sorte que $C$ est r\'eguli\`ere en $\infty$. On peut alors prendre 
$U = C^{-1} (C^{(\infty)})^{-1}$: $U C^{(0)} = C^{-1} M = R$  est r\'eguli\`ere
en $0$, $U C^{(\infty)} = C^{-1}$ est r\'eguli\`ere en $\infty$ et, bien sur,
$U \in GL_{n}(\mathbf{C}(z))$.
\hfill $\Box$



\section{Matrice de connexion et classification}



\subsection{La matrice de connexion}

On a associ\'e canoniquement \`a l'\'equation $\sigma_{q}X = AX$ deux matrices
fondamentales de solutions m\'eromorphes sur $\mathbf{C}^{*}$, $X^{(0)}$ et
$X^{(\infty)}$. Suivant Birkhoff (\emph{voir} \cite{Birkhoff}), on notera
$P$ la \emph{matrice de connexion} $(X^{(\infty)})^{-1} X^{(0)}$. C'est un
\'el\'ement de $GL_{n}(\mathcal{M}(\mathbf{C}^{*}))$. De plus:
$$
\sigma_{q}P = (\sigma_{q}X^{(\infty)})^{-1} \sigma_{q}X^{(0)} = 
(AX^{(\infty)})^{-1} AX^{(0)} = (X^{(\infty)})^{-1} X^{(0)} = P
$$
Autrement dit, les coefficients de $P$ sont elliptiques selon les
identifications du chapitre 1, et $P \in GL_{n}(\mathcal{M}(\mathbf{E}))$. Nous
dirons, pour simplifier, que $P$ est elliptique. En particulier, le lieu
singulier de $P$ est $\sigma_{q}$-invariant, donc une r\'eunion de spirales 
logarithmiques discr\`etes de la forme $q^{\mathbf{Z}}\alpha$. Nous verrons
qu'elles sont en nombre fini et les d\'ecrirons plus pr\'ecis\'ement.\\

Naturellement, le choix de $P$ n'est pas univoque puisque celui des solutions
canoniques locales ne l'est pas. Nous allons \'egalement pr\'eciser cela.

\subsubsection{Cas non-r\'esonnant}

Dans le cas non-r\'esonnant, $A(0) , A(\infty) \in GL_{n}(\mathbf{C})$ et les
valeurs propres distinctes de chacune de ces deux matrices sont deux \`a deux
non congrues modulo $q^{\mathbf{Z}}$. Alors, les matrices fondamentales ont 
\'et\'e obtenues sous les formes:
$$
\begin{cases}
X^{(0)} = Q^{(0)}H^{(0)}L^{(0)}C^{(0)} \\
X^{(\infty)} = Q^{(\infty)}H^{(\infty)}L^{(\infty)}C^{(\infty)} \\
\end{cases}
$$

On peut normaliser en imposant la forme plus rigide aux r\'eduites de Jordan. 
Cependant, $Q^{(0)}$ et $Q^{(\infty)}$ peuvent respectivement \^etre 
remplac\'ees par $Q^{(0)}R^{(0)}$ et $Q^{(\infty)}R^{(\infty)}$, o\`u
$R^{(0)}$ (resp. $R^{(\infty)}$) appartient \`a $GL_{n}(\mathbf{C})$ et
commute avec la r\'eduite de Jordan de $A(0)$ (resp. de $A(\infty)$). Cela
revient \`a choisir, au lieu de $X^{(0)}$ et $X^{(\infty)}$ les solutions
$X^{(0)}R^{(0)}$ et $X^{(\infty)}R^{(\infty)}$ et donc \`a remplacer $P$ par
$(R^{(\infty)})^{-1} P R^{(0)}$.\\

Dans tous les cas, on pourra \'ecrire:
$$
P = ((L^{(\infty)}C^{(\infty)})^{-1}) \times
((H^{(\infty)})^{-1}(Q^{(\infty)})^{-1}Q^{(0)}H^{(0)}) \times
(L^{(0)}C^{(0)})
$$

Les singularit\'es des facteurs exterieurs sont exactement celles des facteurs
log-car des solutions, c'est \`a dire les spirales logarithmiques discr\`etes
li\'ees aux structures de Jordan de $A(0)$ et de $A(\infty)$. Les 
singularit\'es du facteur central sont contenues dans les spirales discr\`etes 
engendr\'ees par les singularit\'es de $A$.\\

La forme canonique ci-dessus sera pr\'ef\'erentiellement utilis\'ee pour faire
de l'analyse, par exemple pour \'etudier le comportement de $P$ lorsque
$q \to 1$ aux chapitres 3 et 4.

\subsubsection{Cas g\'en\'eral}
 
Dans le cas g\'en\'eral, on a les formes canoniques (moins pr\'ecises en ce qui
concerne la partie m\'eromorphe):
$$
\begin{cases}
X^{(0)} = M^{(0)} N^{(0)} \\
X^{(\infty)} = M^{(\infty)} N^{(\infty)} \\
\end{cases}
$$
d'o\`u l'on tire 
$$
P = (N^{(\infty)})^{-1} \times ((M^{(\infty)})^{-1} M^{(0)}) \times N^{(0)}
$$

Si l'on normalise en imposant la forme plus rigide aux parties log-car, les
solutions canoniques $X^{(0)}$ et $X^{(\infty)}$ sont soumises exactement \`a
la m\^eme ind\'etermination que ci-dessus, et il en est donc de m\^eme de $P$.
C'est sous cette forme que nous \'etudierons les propri\'et\'es classifiantes
de la matrice de connexion.\\

Dans tous les cas, on codera la matrice de connexion sans perdre trace des
parties log-car (en vue du th\'eor\`eme de classification au 2.2), sous l'une 
des deux formes \'equivalentes suivantes:
$$
(N^{(\infty)},P,N^{(0)}) \text{  ou bien  } 
(N^{(\infty)},M,N^{(0)}) \text{ o\`u } M = (M^{(\infty)})^{-1} M^{(0)}
$$

Les composantes de ces triplets codants ne sont d'ailleurs pas libres. Outre la
structure log-car impos\'ee \`a $N^{(0)}$ et \`a $N^{(\infty)}$, il y a
l'\'equation $\sigma_{q} P = P$ et l'\'equation qui s'en d\'eduit pour $M$:
$\sigma_{q}M = K^{(\infty)} M (K^{(0)})^{-1}$. Ici, l'on note
$K^{(0)} = (\sigma_{q}N^{(0)}) (N^{(0)})^{-1}$ et
$K^{(\infty)} = (\sigma_{q}N^{(\infty)}) (N^{(\infty)})^{-1}$; ce sont des 
\'el\'ements de $GL_{n}(\mathbf{C})$.\\

L'avantage principal du second codage est de faire clairement ressortir les
diff\'erents groupes de singularit\'es, qui auront tant d'importance aux
chapitres 3 et 4.

\subsubsection{Codage \`a la Birkhoff}

On peut rapprocher les codages ci-dessus de l'esprit de Birkhoff en leur 
donnant une allure un peu plus alg\'ebrico-combinatoire. Birkhoff pose en effet
le probl\`eme ``inverse'' \`a partir de la donn\'ee de $P$ et de celle des
exposants en $0$ et en $\infty$: n'\'etudiant que le cas g\'en\'erique 
semi-simple, cela lui tient en effet lieu de structure de Jordan.\\

La donn\'ee de $N^{(0)}$ (resp. de $N^{(\infty)}$) \'equivaut \`a celle de
$K^{(0)}$ (resp. de $K^{(\infty)}$), qui est \`a son tour d\'etermin\'ee par
la donn\'ee des exposants et des tailles des blocs correspondants. On peut de
plus consid\'erer les exposants modulo $q^{\mathbf{Z}}$, c'est \`a dire, comme
des \'el\'ements de $\mathbf{E}$.\\

La donn\'ee d'un coefficient de $P$, qui est une fonction elliptique,
\'equivaut \`a celle de son diviseur des z\'eros et des p\^oles, que l'on peut
consid\'erer comme un diviseur sur $E$ de degr\'e $0 \in \mathbf{Z}$ et 
d'\'evaluation $0 \in \mathbf{E}$; plus celle d'un facteur constant non nul.
Ceci, en vertu de la suite exacte:
$$
\{1\} \rightarrow \mathbf{C}^{*} \rightarrow \mathcal{\mathbf{E}}^{*}
\rightarrow Div^{(0)}(\mathbf{E}) \rightarrow E \rightarrow \{0\}
$$

Ces codages se comportent particuli\`erement bien du point de vue du produit
tensoriel, mais les autres op\'erations s'y traduisent moins facilement: par
exemple, il n'est pas simple d'y lire le lieu des z\'eros de $P$; ni, 
d'ailleurs, d'expliciter un choix syst\'ematique des facteurs constants 
ci-dessus mentionn\'es.



\subsection{Effet de l'\'equivalence rationnelle}

On a associ\'e \`a une \'equation singuli\`ere r\'eguli\`ere 
$\sigma_{q} X = AX$, le triplet $(N^{(\infty)},P,N^{(0)})$. Celui-ci n'est
cependant pas d\'efini de mani\`ere univoque: $P$ peut \^etre remplac\'e par
$(R^{(\infty)})^{-1} P R^{(0)}$ si 
$R^{(0)} \;,\; R^{(\infty)} \in GL_{n}(\mathbf{C})$ commutent respectivement
avec $N^{(\infty)}$ et $N^{(0)}$. Ceci justifie l'introduction des notations
suivantes.\\

Sur l'ensemble des triplets $(N_{1},P,N_{2})$, o\`u 
$P \in GL_{n}(\mathcal{M}(\mathbf{E}))$ et o\`u $N_{1}, N_{2}$ sont des 
matrices log-car en forme normale, on d\'efinit une relation: 
$$
(N_{1},P,N_{2}) \sim (N'_{1},P',N'_{2}) 
\text { si et seulement si }
\begin{cases}
N_{1} = N'_{1} \;,\; N_{2} = N'_{2} \\
R_{1} P' = P R_{2} 
\end{cases} 
$$
o\`u $R_{1},R_{2} \in GL_{n}(\mathbf{C})$ et
$[R_{1},N_{1}] = [R_{2},N_{2}] = 0$.  C'est alors une relation d'\'equivalence, et on notera $\mathcal{F}_{n}$ le 
quotient.\\

L'image dans $\mathcal{F}_{n}$ du triplet $(N^{(\infty)},P,N^{(0)})$ est alors
bien d\'efinie d'apr\`es le 2.1. On va maintenant prouver le\\

\textbf{\underline{Th\'eor\`eme}}\\

\emph{On obtient ainsi une
bijection de l'ensemble $\mathcal{E}_{n}$ des classes d'\'equivalence
rationnelle de syst\`emes singuliers r\'eguliers de rang $n$ dans l'ensemble
$\mathcal{F}_{n}$.}\\

\subsubsection{L'application est bien d\'efinie}

Si $B = (\sigma_{q}U)^{-1}AU$ et si les solutions fondamentales choisies pour 
$A$ sont $M^{(0)}N^{(0)}$ et $M^{(\infty)}N^{(\infty)}$, on peut choisir pour
$B$ les solutions $(UM^{(0)})N^{(0)}$ et $(UM^{(\infty)})N^{(\infty)}$, ce qui
associe \`a $B$ les m\^emes parties log-car et la m\^eme matrice de connexion,
donc le m\^eme triplet codant qu'\`a $A$. Ainsi, on a une application bien
d\'efinie de $\mathcal{E}_{n}$ dans $\mathcal{F}_{n}$.

\subsubsection{Injectivit\'e}

R\'eciproquement, supposons que $A$ et $B$ ont fourni des triplets 
\'equivalents. Quitte \`a modifier le choix des solutions pour l'une de ces
deux \'equations, on peut alors supposer que l'on a
associ\'e \`a $A$ et $B$ le \emph{m\^eme} triplet
$(N^{(\infty)},P,N^{(0)})$. Ainsi, on a des solutions fondamentales
$M_{1}^{(0)}N^{(0)}$ et $M_{1}^{(\infty)}N^{(\infty)}$ pour $A$, et
$M_{2}^{(0)}N^{(0)}$ et $M_{2}^{(\infty)}N^{(\infty)}$ pour $B$. De plus,
$(M_{1}^{(\infty)})^{-1}M_{1}^{(0)} = (M_{2}^{(\infty)})^{-2}M_{1}^{(0)}$ car
ces deux matrices sont \'egales \`a $N^{(\infty)}P(N^{(0)})^{-1}$. On a donc
$M_{1}^{(0)}(M^{(0)})^{-1} = M_{1}^{(\infty)}(M^{(\infty)})^{-1}$. Notant $U$
cette matrice, on voit qu'elle est inversible et m\'eromorphe \`a la fois sur
$\mathbf{C}$ et sur $\mathbf{S} - \{0\}$, donc rationnelle: et il est clair
que $B = (\sigma_{q}U)^{-1}AU$, autrement dit, on a \'equivalence rationnelle.
Ceci prouve l'injectivit\'e.\\

\textbf{\underline{Remarque}}\\

En omettant les parties log-car dans le codage, on n'aurait plus 
l'injectivit\'e. Par exemple, en dimension $1$, $A = 1$ et 
$B = c \in \mathbf{C}^{*} - q^{\mathbf{Z}}$ ne sont pas \'equivalentes, mais
on peut \`a toutes deux associer la matrice de connexion $P = 1$.

\subsubsection{Surjectivit\'e}

Reste \`a \'etablir le seul point non trivial, la surjectivit\'e. On part d'un
triplet $(N^{(\infty)},P,N^{(0)})$ ($P$ elliptique, $N^{(\infty)}$ et $N^{(0)}$
des matrices log-car) et l'on \'ecrit $P = (N^{(\infty)})^{-1}MN^{(0)}$ o\`u 
$M \in GL_{n}(\mathcal{M}(\mathbf{C}^{*}))$. On invoque alors la 
``trivialit\'e des fibr\'es m\'eromorphes sur $\mathbf{S}$'' (\emph{voir} \cite
{SVdP} p.158): cela entraine que $M$ s'\'ecrit $(M^{(\infty)})^{-1}M^{(0)}$,
o\`u $M^{(0)} \in GL_{n}(\mathcal{M}(\mathbf{C}))$ et
$M^{(\infty)} \in GL_{n}(\mathcal{M}(\mathbf{S} - \{0\}))$.\\

Les matrices
$X^{(0)} = M^{(0)}N^{(0)}$ et $X^{(\infty)} = M^{(\infty)}N^{(\infty)}$ sont 
connect\'ees par la matrice elliptique $P$: $X^{(\infty)} = P X^{(0)}$, et 
sont donc solutions de la m\^eme \'equation $\sigma_{q}X = AX$, en prenant 
pour $A$ la matrice
$\sigma_{q}X^{(0)}(X^{(0)})^{-1} = \sigma_{q}X^{(\infty)}(X^{(\infty)})^{-1}$.\\

La premi\`ere \'ecriture implique que 
$A = \sigma_{q}M^{(0)} K^{(0)} (M^{(0)})^{-1}$, o\`u l'on a encore pos\'e
$K^{(0)} = (\sigma_{q}N^{(0)}) (N^{(0)})^{-1} \in GL_{n}(\mathbf{C})$, ce qui 
montre que $A$ est m\'eromorphe sur $\mathbf{C}$. Le raisonnement analogue \`a
l'infini montre que $A$ est m\'eromorphe sur $\mathbf{S} - \{0\}$ et l'on 
conclut que $A \in GL_{n}(\mathbf{C}(z))$.\\

L'\'equation aux $q$-diff\'erences lin\'eaire \`a coefficients rationnels
$\sigma_{q}X = AX$ admettant les deux solutions en formes normales 
$X^{(0)} = M^{(0)}N^{(0)}$ et $X^{(\infty)} = M^{(\infty)}N^{(\infty)}$,
il r\'esulte de 2.1 qu'elle est singuli\`ere r\'eguli\`ere et ceci ach\`eve
la preuve du th\'eor\`eme.
\hfill $\Box$

\subsection{Pr\'ecisions et am\'eliorations diverses}

\begin{enumerate}

\item{Le deuxi\`eme th\'eor\`eme de 2.1.5 affirme que notre \'equation est
\'equivalente \`a une \'equation fuchsienne en $0$ et $\infty$; cela revient
\`a dire que l'on peut obtenir $M^{(0)}$ et $M^{(\infty)}$ respectivement
r\'eguli\`eres en $0$ et en $\infty$.}

\item{En vue des applications \`a la th\'eorie de Galois (\emph{voir} \cite{JS2}),
il est utile de pouvoir pr\'eciser la position des singularit\'es. Celles de 
$N^{(0)}$
et de $N^{(\infty)}$ sont totalement d\'etermin\'ees par la structure de Jordan
de $A(0)$ et de $A(\infty)$.\\

En ce qui concerne $M$, on d\'eduit facilement ce
qui suit de la preuve du th\'eor\`eme 3. Soit une partie finie
$\Sigma \subset \mathbf{C}^{*}$ ; pour simplifier la formulation, nous
supposerons que deux \'el\'ements distincts de $\Sigma$ ne sont pas congrus
modulo $q^{\mathbf{Z}}$.

\begin{itemize}

\item{Sens direct : partant d'un syst\`eme (1) tel que 
$\mathcal{S}(A) \subset \Sigma$, on voit que les singularit\'es de $M$ 
appartiennent \`a $q^{\mathbf{Z}}\Sigma$.}

\item{Sens inverse : supposons que les singularit\'es de $M$ appartiennent \`a 
$q^{\mathbf{Z}}\Sigma$. L'invocation du ``Preliminary Theorem'' de 
\cite{Birkhoff} (lemme de Birkhoff) \`a
la place de l'argument de trivialit\'e des fibr\'es m\'eromorphes permet de
choisir les singularit\'es de $M^{(0)}$ et $M^{(\infty)}$ respectivement dans
$q^{\mathbf{N}^{*}}\Sigma$ et dans $q^{-\mathbf{N}}\Sigma$; le syst\`eme (1)
correspondant est alors tel que $\mathcal{S}(A) \subset \Sigma$.}

\end{itemize}
}

\item{La normalisation un peu artificielle des de la forme des matrices log-car a
pour but d'obtenir un \'enonc\'e pas trop compliqu\'e. Si on relaxe cette
condition, la relation d'\'equivalence qui sert \`a d\'efinir $\mathcal{F}_{n}$
se traduit, pour les parties log-car, par les conditions:
$R_{1}N_{1} = N'_{1}R_{1}$ et $R_{2}N_{2} = N'_{2}R_{2}$.\\

La forme normale des matrices log-car n'\'etant pas pr\'eserv\'ee par le 
produit tensoriel, nous serons n\'ecessairement amen\'es \`a cette d\'efinition
\'elargie dans la partie de ce travail concernant la th\'eorie de Galois
(\emph{voir} \cite{JS2}).}

\end{enumerate}


\addtocontents{toc}{\protect\contentsline{part}{\protect\numberline{}Confluence vers un syst\`eme diff\'erentiel,\\matrice de connexion et
monodromie}{\thepage}}
\part*{Confluence vers un syst\`eme diff\'erentiel, matrice de connexion et
monodromie}

\chapter{Comportement des solutions canoniques lorsque $q \rightarrow 1$}



\section{Comportements asymptotiques lorsque $q \rightarrow 1$}

\subsection*{Conventions et notations g\'en\'erales}

Dor\'enavant $A(z)$ sera aussi consid\'er\'ee comme fonction de $q$, et l'on
cherchera \`a comprendre le comportement des solutions canoniques lorsque
$q \to 1$. Dans tout ce qui suit, on supposera que $q \to 1$ le long d'une 
spirale logarithmique fix\'ee $q_{0}^{\mathbf{R}_{+}}$. Pr\'ecis\'ement:
$$
q_{0} = e^{ - 2 \imath \pi \tau_{0}} \qquad 
q = e^{ - 2 \imath \pi \tau} \qquad
\tau = \tau_{0} \epsilon
$$
o\`u $\tau_{0}$ est fix\'e, \`a partie imaginaire strictement positive et o\`u
$\epsilon$ est un r\'eel strictement positif
que l'on fera tendre vers $0$. Nous emploierons la notation $q \to 1$ comme
abr\'eviation (abusive) pour $\epsilon \to 0^{+}$. Pour $t$ r\'eel, les 
notations $q_{0}^{t}$, $q^{t}$ d\'esigneront alors
respectivement $e^{- 2 \imath \pi \tau_{0} t}$ et $e^{- 2 \imath \pi \tau t}$.
Ainsi $q = q_{0}^{\epsilon}$ et $q^{t} = q_{0}^{\epsilon t}$. Ces notations seront
\'etendues \`a la droite num\'erique achev\'ee $\overline{\mathbf{R}}$ en 
posant $q^{+ \infty} = q_{0}^{+ \infty} = \infty \in \mathbf{S}$ et 
$q^{- \infty} = q_{0}^{- \infty} = 0$.\\

On introduit les demi-spirales logarithmiques compactes: 
$q_{0}^{\overline{\mathbf{R}_{+}}}$,
qui relie $1$ \`a $\infty$; $q_{0}^{\overline{\mathbf{R}_{-}}}$, qui relie $1$ 
\`a $0$: et leur r\'eunion, la spirale logarithmique compacte 
$q_{0}^{\overline{\mathbf{R}}}$, qui relie $0$ \`a $\infty$. Celles-ci nous 
serviront de coupures, autrement dit, on rendra uniformes certaines fonctions
analytiques sur leurs compl\'ementaires, les ouverts simplement connexes
$\Omega_{0} = \mathbf{S} - q_{0}^{\overline{\mathbf{R}_{+}}}$,
$\Omega_{\infty} = \mathbf{S} - q_{0}^{\overline{\mathbf{R}_{-}}}$,
$\Omega = \Omega_{0} \cap \Omega_{\infty} = 
\mathbf{S} - q_{0}^{\overline{\mathbf{R}}}$.\\


\bigskip
\hrule
\bigskip

\unitlength = 1cm

\underline{Figure 1}\\

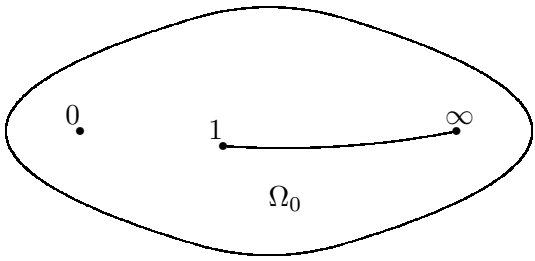
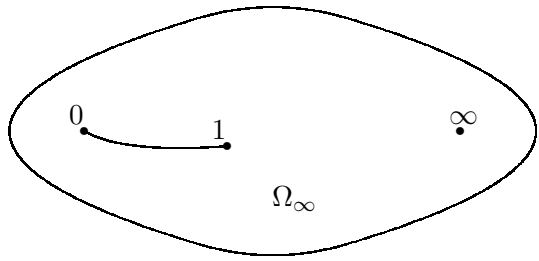
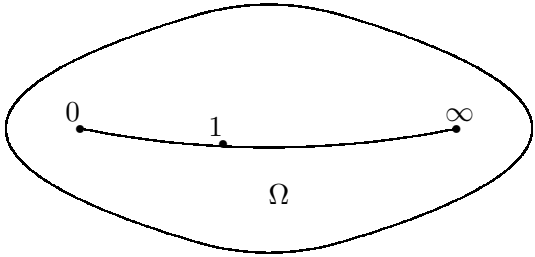
\bigskip

\begin{picture} (15,3)


\qbezier (4,0)(9,1.5)(4,3)     

\qbezier (2,0)(-3,1.5)(2,3)      

\qbezier (2,3)(3,3.3)(4,3)   

\qbezier (2,0)(3,-0.3)(4,0)  


\put(0.5,1.5){\circle*{0.1}}
\put(0.3,1.6){$0$}

\put(2.4,1.3){\circle*{0.1}}
\put(2.2,1.4){$1$}

\put(5.5,1.5){\circle*{0.1}}
\put(5.35,1.6){$\infty$}


\qbezier (2.4,1.3) (4,1.2) (5.5,1.5)


\put(3,0.5){$\Omega_{0}$}


\qbezier (12,0)(17,1.5)(12,3)     

\qbezier (10,0)(5,1.5)(10,3)      

\qbezier (10,3)(11,3.3)(12,3)   

\qbezier (10,0)(11,-0.3)(12,0)  


\put(8.5,1.5){\circle*{0.1}}
\put(8.3,1.6){$0$}

\put(10.4,1.3){\circle*{0.1}}
\put(10.2,1.4){$1$}

\put(13.5,1.5){\circle*{0.1}}
\put(13.35,1.6){$\infty$}


\qbezier (8.5,1.5) (9,1.2) (10.4,1.3)


\put(11,0.5){$\Omega_{\infty}$}

\end{picture}

\bigskip
\bigskip

\begin{picture} (15,3)


\qbezier (4,0)(9,1.5)(4,3)     

\qbezier (2,0)(-3,1.5)(2,3)      

\qbezier (2,3)(3,3.3)(4,3)   

\qbezier (2,0)(3,-0.3)(4,0)  


\put(0.5,1.5){\circle*{0.1}}
\put(0.3,1.6){$0$}

\put(2.4,1.3){\circle*{0.1}}
\put(2.2,1.4){$1$}

\put(5.5,1.5){\circle*{0.1}}
\put(5.35,1.6){$\infty$}


\qbezier (0.5,1.5) (3,1) (5.5,1.5)


\put(3,0.5){$\Omega$}

\end{picture}

\bigskip
\hrule
\bigskip

\newpage

\textbf{\underline{Remarque}}\\

L'hypoth\`ese sur la fa\c{c}on dont $q \to 1$ a pour but d'imposer que les
demi-spirales discr\`etes de p\^oles de la forme $z_{0} q^{\mathbf{N}}$, qui
sont apparues naturellement aux chapitres 1 et 2, se condensent en des coupures
de la forme $z_{0} q_{0}^{\mathbf{R}_{+}}$. Par le changement de variables
$z = e^{2 \imath \pi x}$, cela revient \`a dire que le r\'eseau 
$\mathbf{Z} + \mathbf{Z} \tau = \mathbf{Z} + \mathbf{Z} \tau_{0} \epsilon$ se 
condense en le sous-groupe ferm\'e $\mathbf{Z} + \mathbf{R} \tau_{0}$ de
$\mathbf{R}$. Cette g\'eom\'etrie est encore respect\'ee si $q \to 1$ 
\emph{tangentiellement \`a une spirale logarithmique}, autrement dit si 
$\tau = \tau_{0} \epsilon + o(\epsilon)$: il est donc probable que tous les
r\'esultats qui suivent restent valables sous cette hypoth\`ese plus faible.

\subsection*{Contenu de ce chapitre}

Dans la section 3.1, nous reprenons la ``th\'eorie des fonctions'' de 1.1, mais
du point de vue de l'\'etude asymptotique lorsque $q \to 1$. Dans la section 3.2, nous
\'etudions le comportement asymptotique de solutions d'une famille 
d'\'equations aux $q$-diff\'erences \emph{\`a conditions initiales (valeur en
$0$) fix\'ees}. Cette \'etude est essentiellement ind\'ependante de celle de
3.1. Dans la section 3.3, nous d\'eduisons de 3.1 le comportement de la partie
log-car (sous la forme $L C$ aussi bien que sous la forme $e_{q,A}$). La 
synth\`ese a lieu en 3.4, o\`u nous explicitons le comportement des
solutions canoniques.\\

Une partie des r\'esultats de 3.1 ne sert pas \`a la th\'eorie g\'en\'erale
mais seulement \`a l'\'etude d'exemples \`a la fin du chapitre 4.\\

D\'ecrivons maintenant la strat\'egie mise en oeuvre pour traiter le 
probl\`eme. Nous partons d'une \'equation diff\'erentielle lin\'eaire 
fuchsienne donn\'ee sous la forme d'un syst\`eme d'ordre $1$:
$$
\delta \tilde{X} = \tilde{B} \tilde{X} \qquad (\delta := \frac{d}{dz})
$$
on va montrer comment la matrice de connexion d'une famille (param\'etr\'ee
par $q$) d'\'equations aux $q$-diff\'erences fuchsiennes donne, lorsque $q$
tend vers $1$ le long d'une spirale, la description compl\`ete de la
monodromie dans une base explicite.\\

Pour que les solutions $X$ d'une \'equation 
$$
\sigma_{q} X = A X
$$
(comme au 3.1, les d\'ependances en $q$ ne sont pas n\'ecessairement explicites
dans les notations) tendent vers $\tilde{X}$, on r\'e\'ecrit ce syst\`eme sous
la forme
$$
\delta_{q} X = B X \qquad (\delta_{q} := \frac{\sigma_{q} - Id}{q - 1})
$$
o\`u $B = \frac{A - I_{n}}{q - 1}$ et o\`u l'on suppose que $B \to \tilde{B}$
en un sens convenable. Cette heuristique est justifi\'ee par les exemples
\'etudi\'es par voie directe (cf 4.4), en particulier le cas de la s\'erie
hyperg\'eom\'etrique basique. Nous privil\'egions donc, parmi les 
$q$-analogies classiques, la mise en relation
$\delta \leftrightarrow \frac{\sigma_{q} - Id}{q - 1}$. En fait, on a un peu 
l'id\'ee que $\delta_{q} \to \delta$.\\

La m\'ethode consiste alors \`a \'ecrire un syst\`eme fondamental sous la forme
obtenue au 1.2.4:
$$
X = Q_{0} H L C
$$
puis \`a controler s\'epar\'ement les diff\'erentes composantes:

\begin{enumerate}

\item{$Q_{0}$ est la matrice de passage \`a la r\'eduite de Jordan pour $A(0)$,
donc \'egalement pour $B(0)$; elle devra tendre vers son analogue 
$\tilde{Q}_{0}$ pour $\tilde{B}(0)$. On voudra donc une \emph{d\'eformation
\`a structure de Jordan fix\'ee}. Ceci sera pr\'ecis\'e plus loin.}

\item{$H$ a \'et\'e caract\'eris\'ee comme solution d'une \'equation aux
$q$-diff\'erences matricielle, avec condition initiale $H(0) = I_{n}$. Elle
sera consid\'er\'ee comme solution d'une \'equation aux $q$-diff\'erences 
vectorielle d'ordre $n^{2}$ avec condition initiale fixe (ind\'ependante de
$q$). La conclusion dans ce cas, qui est \`a la fois un cas particulier et une
\'etape vers le cas g\'en\'eral, est l'objet de la section 3.2}

\item{La convergence de la partie log-car $LC$ sera trait\'ee au 3.3 \`a
l'aide des r\'esultats de 3.1.3 et 3.1.4. Cependant, ces derniers, ainsi que
l'exemple, trait\'e au chapitre 4, de la d\'eformation de $\log(1 - z)$, ont
montr\'e la n\'ecessit\'e de renormaliser $l_{q}$ et nous devrons modifier
en cons\'equence la d\'efinition de nos solutions canoniques.}

\end{enumerate}

On pourra alors montrer, au 3.4, que le d\'evissage de la r\'esolution de
l'\'equation aux $q$-diff\'erences et celui de l'\'equation diff\'erentielle
se correspondent, de sorte que notre solution canonique tend vers la solution
canonique de l'\'equation diff\'erentielle obtenue par la m\'ethode de
Frobenius\footnote{A ceci pr\`es qu'on y remplacera les $\log z$ et $z^{\mu}$ usuels
par $\log(-z)$ et $(-z)^{\mu}$}.\\

On traitera \'egalement le cas des solutions b\^aties \`a l'aide des matrices
$e_{q,A}$, qui s'av\`erera un peu plus simple. Par exemple, il ne sera pas
n\'ecessaire de renormaliser ces solutions.\\

Tout ceci concerne les solutions locales en $0$. L'analogue \`a l'infini et
l'\'etude globale, c'est \`a dire celle de la matrice de connexion, feront
l'objet du chapitre 4.


\subsection{Lemmes pr\'eliminaires utiles}

\textbf{\underline{Lemme 1}}\\

\emph{Soit $z_{0} \in \Omega_{0}$. Il existe alors $\rho,\delta > 0$ tels que:
$$\forall z \in \overset{\circ}{D}(z_{0},\rho) \subset \Omega_{0} \;:\;
\inf_{\alpha > 0} |1 - q_{0}^{-\alpha} z| \geq \delta$$}\\

\textbf{\underline{Preuve}}\\

Nions la conclusion du lemme: pour tout $n \in \mathbf{N}^{*}$, posant
$\rho_{n} = \delta_{n} = \frac{1}{n}$, il existe 
$z_{n} \in \overset{\circ}{D}(z_{0},\rho_{n})$ et $\alpha_{n} > 0$ tels que
$|1 - q_{0}^{-\alpha_{n}} z_{n}| < \delta_{n}$. Alors
$z_{n} \underset{n \to \infty}{\rightarrow} z_{0}$ et
$z_{n}q_{0}^{-\alpha_{n}} \underset{n \to \infty}{\rightarrow} 1$, de sorte que
$q_{0}^{\alpha_{n}} \underset{n \to \infty}{\rightarrow} z_{0}$. Celui-ci
appartient donc \`a la spirale $q_{0}^{\overline{\mathbf{R}_{+}}}$, puisque
celle-ci est ferm\'ee: contradiction. \hfill $\Box$\\

\textbf{\underline{Lemme 2}}\\

\emph{
(i) Soit $a \in \mathbf{C}$ tel que $Re(a) >0$. Alors
$$1 \leq \frac{|e^{a} - 1|} {|e^{a}| - 1} \leq \frac{|a|} {Re(a)}$$
(ii) Si $\lambda > 0$:
$$1 \leq \frac{|q_{0}^{\lambda} - 1|} {|q_{0}^{\lambda}| - 1} \leq
\frac{|\tau_{0}|} {Im(\tau_{0})}$$}\\

\textbf{\underline{Preuve}}\\

La minoration dans (i) d\'ecoule de l'in\'egalit\'e du triangle et du fait que
$|e^{a}| = e^{Re(a)} > 1$. La majoration vient du calcul suivant:

$$|e^{a} - 1| \leq \int_{0}^{1} |\frac{d}{ds}(e^{sa})| ds =
\int_{0}^{1} |a| e^{s Re(a)} ds$$

L'appliquant maintenant avec $a = - 2 \imath \pi \tau_{0} \lambda$,
on tire (ii). \hfill $\Box$\\

Notons que (ii) s'appliquera \`a $q = q_{0}^{\epsilon}$ et aux 
$q^{m} \;,\; m \in \mathbf{N}$:\\

\textbf{\underline{Corollaire}}\\

\emph{
Si $m \in \mathbf{N}$:
$$ 
\frac{Im(\tau_{0})} {|\tau_{0}|} \leq 
\left|\frac{q^{m} - 1} {q - 1}\right| \div \frac{|q|^{m} - 1} {|q| - 1} \leq 1
$$
$$ 
\left|\frac{q^{m} - 1} {q - 1}\right| \geq m \frac{Im(\tau_{0})} {|\tau_{0}|}
$$}\\

\textbf{\underline{Preuve}}\\

La minoration dans (i) d\'ecoule du lemme 1 appliqu\'e \`a $q$ et \`a $q^{m}$:

$$
\left.
\begin{array}{c}
1 \leq \frac{|q - 1|} {|q| - 1} \leq \frac{|\tau_{0}|} {Im(\tau_{0})} \\
1 \leq \frac{|q^{m} - 1|} {|q^{m}| - 1} \leq \frac{|\tau_{0}|} {Im(\tau_{0})}
\end{array}
\right\}
\Rightarrow
\frac{Im(\tau_{0})} {|\tau_{0}|} \leq 
\left|\frac{q^{m} - 1} {q - 1} \right| \div \frac{|q|^{m} - 1} {|q| - 1}
$$

La majoration vient de l'in\'egalit\'e du triangle appliqu\'ee \`a
$\sum_{j = 0}^{m-1} q^{j}$. L'in\'egalit\'e (ii) d\'ecoule de la premi\`ere
minoration et du fait que $\sum_{j = 0}^{m-1} |q^{j}| \geq m$.
\hfill $\Box$

\subsection{Comportement de $\Theta_{q}(q^{\alpha}z)$}

Posons $z = e^{2 \imath \pi x}$ et, comme d'habitude, 
$q = e^{- 2 \imath \pi \tau}$ ($Im(\tau) > 0$). Alors

$$\Theta_{q}(z) = 
\sum_{m \in \mathbf{Z}} (-1)^{m} q^{- \frac{m^{2} - m}{2}} z^{m} =
\sum_{m \in \mathbf{Z}} 
e^{\imath \pi \tau m^{2} + 2 \imath \pi m (x + \frac{1}{2} - \frac{\tau}{2})}$$

C'est donc $\theta(x + \frac{1 - \tau}{2},\tau)$, o\`u $\theta$ est la fonction
introduite par Mumford dans \cite{Mumford}. Nous noterons ici plut\^ot
$\tilde{\Theta}(\tau,x) = \Theta_{q}(e^{2i\pi x})$.\\

\textbf{\underline{Proposition}}\\

\emph{$\tilde{\Theta}$ v\'erifie l'\'equation fonctionnelle:
$$
\tilde{\Theta}(\tau,x) =  
\frac{\sqrt{i/\tau}}{e^{(i\pi/\tau)(x-\frac{\tau+1}{2})^{2}}}
\tilde{\Theta}(-1/\tau,-x/\tau)
$$
}

Comme $Re(\imath/\tau) > 0$, on a ici affaire \`a la d\'etermination principale
de la racine carr\'ee sur $\mathbf{C} - \mathbf{R}_{-}$ (telle que $1 \mapsto 1$).\\

\textbf{\underline{Preuve}}\\

C'est un cas tr\`es particulier de l'\'equation fonctionnelle d'automorphie
donn\'ee dans \cite{Mumford} (p. 32), mais nous en donnerons une preuve plus 
directe. Nous nous ram\`enerons au cas r\'eel, qui se traite avec la formule
sommatoire de Poisson (\emph{voir} \cite{Lang}, p. 269 et \cite{WW}, p. 475).\\

La convergence normale sur tout compact montre que $\tilde{\Theta}$ est
analytique en chacune des variables $x \in \mathbf{C}$ et 
$\tau \in \mathcal{H}$ (demi-plan de Poincar\'e, d\'efini par $Im(\tau) > 0$). 
Il
suffit donc d'\'etablir la formule ci-dessus pour $\tau = \imath/\lambda$ avec
$\lambda > 0$ et elle s'\'etendra alors par le principe du prolongement
analytique. En y rempla\c{c}ant $x$ par $x + \frac{\imath/\lambda + 1}{2}$, on
est ramen\'e \`a prouver:

$$\forall \lambda > 0 \;,\;
\tilde{\Theta}(\frac{\imath}{\lambda},x + \frac{\imath/\lambda + 1}{2}) =
\frac{\sqrt{\lambda}}{e^{\pi \lambda x^{2}}}
\tilde{\Theta}(\imath \lambda, \imath \lambda x + \frac{\imath \lambda - 1}{2})$$

soit encore, apr\`es substitution dans les s\'eries concern\'ees:

$$
\sum_{m \in \mathbf{Z}} e^{- \frac{\pi m^{2}}{\lambda} + 2 \imath \pi m (x+1)}
=
\sqrt{\lambda} \sum_{m \in \mathbf{Z}} e^{- \pi \lambda (x+m)^{2}}
$$

Pour cela, on emploiera la formule sommatoire de Poisson. On introduit la
fonction

$$\phi_{\lambda}(x) = e^{- \pi \lambda x^{2}}$$

dont la transform\'ee de Fourier est

$$\hat{\phi}_{\lambda}(x) = \frac{1}{\sqrt{\lambda}} e^{- \pi x^{2} / \lambda}$$

Toutes deux sont dans la classe de Schwarz, c'est \`a dire des fonctions
\`a d\'ecroissance 
rapide sur $\mathbf{R}$ ainsi que toutes leurs d\'eriv\'ees, et, quelque soit
$x$:

$$\sum_{m \in \mathbf{Z}} \hat{\phi}_{\lambda}(m) e^{2 \imath \pi m x} =
\sum_{m \in \mathbf{Z}} \phi_{\lambda}(x+m)$$

soit encore:

$$\sum_{m \in \mathbf{Z}} 
\frac{1}{\sqrt{\lambda}} e^{-\pi \frac{m^{2}}{\lambda} + 2 \imath \pi m x} =
\sum_{m \in \mathbf{Z}} e^{- \pi \lambda (x + m)^{2}}$$

Le membre de droite de cette \'egalit\'e est $1$-p\'eriodique. C'est, aux 
notations  pr\`es, la formule souhait\'ee: celle-ci est donc
\'etablie pour $x$ r\'eel, donc pour $x$ quelconque par une nouvelle
application du principe du prolongement analytique.\hfill $\Box$\\

Il y a \'egalement des preuves ``\'el\'ementaires'' de cette \'equation
fonctionnelle, du type comparaison s\'erie-int\'egrale (\emph{voir par
exemple} \cite{Zuily-Queffelec}, \cite{Godement}).\\

On cherchera, avec cette \'equation fonctionnelle, \`a \'etudier le 
comportement de $\tilde{\Theta}$ lorsque $\tau \to 0$ \`a travers son
comportement lorsque $\tau \to \infty$ (pr\'ecis\'ement: 
$\infty \times \frac{-1}{\tau_{0}}$). Nous voulons un \'equivalent de 
$\Theta_{q}(q^{\alpha}z) = \tilde{\Theta}(\tau,x - \tau \alpha)$ lorsque
$\epsilon = \tau/\tau_{0} \to 0^{+}$, $x$ et $\alpha$ \'etant fix\'es dans
$\mathbf{C}$. L'\'equation fonctionnelle se r\'e\'ecrit

$$\tilde{\Theta}(\tau,x - \tau \alpha) =
\frac{\sqrt{\imath/\tau}}
     {e^{(\imath \pi / \tau)(x - \tau \alpha -\frac{\tau + 1}{2})^{2}}}
\tilde{\Theta}(-1/\tau,-x/\tau + \alpha)$$

Le ``facteur d'automorphie'' vaut

$$\sqrt{\imath/\tau} \; e^{-(\imath \pi / \tau)(x - 1/2)^{2}}
e^{2 \imath \pi (\alpha + 1/2)(x - 1/2)}(1 + O(\tau))$$

avec un $O(\tau)$ uniforme en $x \in \mathbf{C}$. Le facteur th\^eta du 
membre de droite est une somme index\'ee par $m \in \mathbf{Z}$
de termes de la forme

$$e^{2 \imath \pi m (\alpha + \frac{1}{2})} \times
e^{-\frac{\imath \pi}{\tau} (m^{2} + 2 m(x - \frac{1}{2}))}$$

Ecrivons $\tau_{0} = \tau_{1} + \imath \tau_{2}$ (donc $\tau_{2} > 0$), et 
$x = u + \imath v$: $x$ et $u$ ne sont d\'etermin\'es qu'\`a un entier pr\`es, 
puisque c'est la condition $e^{2 \imath \pi x} = z$ qui les s\'pecifie. Il est 
donc loisible de supposer que $u - v \frac{\tau_{1}}{\tau_{2}}$ appartient \`a $[0;1[$,
ce que nous ferons. On a:

$$Re(-\frac{\imath \pi}{\tau} (m^{2} + 2 m(x - \frac{1}{2}))) =
- \frac{\pi Im(\tau_{0})}{\epsilon |\tau_{0}|^{2}}
(m^{2} + 2m(u - \frac{1}{2} - v \frac{\tau_{1}}{\tau_{2}}))$$

\begin{enumerate}

\item{Cas g\'en\'eral: $u - v \frac{\tau_{1}}{\tau_{2}} \in ]0;1[$.\\
Alors, quelque soit $m \in \mathbf{Z} - \{0\}$,
$(m^{2} + 2m(u - \frac{1}{2} - v \frac{\tau_{1}}{\tau_{2}})) > 0$, et tous les
termes correspondants de la s\'erie th\^eta tendent vers $0$, le terme restant
($m = 0$) \'etant constant \'egal \`a $1$. Par vertu de convergence normale,
le facteur th\^eta lui-m\^eme tend vers $1$ dans ce cas:
$\underset{\epsilon \to 0^{+}}{\lim} \tilde{\Theta}(-1/\tau,-x/\tau + \alpha)
= 1$.
}

\item{Cas d\'eg\'en\'er\'e: $u - v \frac{\tau_{1}}{\tau_{2}} = 0$.\\
Dans ce cas, outre le terme constant $1$ correspondant \`a $m = 0$, le terme
d'indice $m = 1$: 
$e^{2 \imath \pi (\alpha + \frac{1}{2})} \times
e^{-\frac{2 \imath \pi v}{\epsilon Im(\tau_{0})}}$ oscille (sauf dans le cas
particulier o\`u $x = 0 \;,\; z = 1$: sa limite est alors 
$1 - e^{2 \imath \pi \alpha}$).
}

\end{enumerate}

On voit que la convergence de $\tilde{\Theta}(-1/\tau,-x/\tau + \alpha)$ vers
$1$ est exponentiellement rapide sur tout compact de $\Omega$, autrement dit,
le terme d'erreur est domin\'e par $e^{-C/\epsilon}$ pour une certaine
constante positive $C$.\\

Pour comprendre la formule, et les deux cas consid\'er\'es, on introduit la 
d\'etermination principale de $\log(-z)$ sur l'ouvert simplement connexe
$\Omega = \mathbf{S} - q_{0}^{\mathbf{\overline{R}}} = 
\mathbf{C}^{*} - q_{0}^{\mathbf{R}}$
qui satisfait $-1 \mapsto 0$. Le ``cas g\'en\'eral'' ci-dessus est celui o\`u 
$x \not \in \tau_{0} \mathbf{R} \pmod{\mathbf{Z}}$, c'est \`a dire
$z \in \Omega$ et le choix de $x$ est alors tel que 
$2 \imath \pi(x - 1/2) = \log(-z)$.
Du point de vue num\'erique, $\Theta_{q}$ a un comportement ``chaotique'' 
le long de la barri\`ere des p\^oles (dans le plan des $x$: le r\'eseau
$\mathbf{Z} + \mathbf{Z} \tau$), lorsque celle-ci se condense en une coupure
spirale (dans le plan des $x$: le sous-groupe ferm\'e non discret
$\mathbf{Z} + \mathbf{R} \tau_{0}$).\\

De la discussion ci-dessus, on d\'eduit le\\

\textbf{\underline{Th\'eor\`eme}}\\

\emph{
(i) Pour $z \in \Omega$ et $\alpha \in \mathbf{C}$, on a, lorsque
$\epsilon \to 0^{+}$:
$$\Theta_{q}(q^{\alpha}z) \sim 
\sqrt{\frac{\imath}{\tau}} \;
e^{-(\imath \pi/\tau)(x-1/2)^{2} + 2\imath \pi(\alpha+1/2)(x-1/2)}$$
(ii) Si de plus $\beta \in \mathbf{C}$:
$$\underset{\epsilon \to 0^{+}}{\lim} 
\frac{\Theta_{q}(q^{\alpha}z)}{\Theta_{q}(q^{\beta}z)} = 
(-z)^{\alpha - \beta}$$}
\hfill $\Box$\\

Cette derni\`ere notation d\'esigne bien sur $e^{(\alpha - \beta)\log(-z)}$.


\subsection{Application aux caract\`eres}

On consid\`erera seulement la famille ``standard'' des
$e_{q,c} = \frac{\Theta_{q}}{\Theta_{q,c}}$ (cf. 1.2), mais il y a des
formules analogues pour des caract\`eres de la forme 
$\prod \Theta_{q,c_{i}}^{r_{i}}$ (qui en sont bien, du moment que
$\sum r_{i} = 0$).\\

\textbf{\underline{Lemme}}\\

\emph{ 
Soit $\gamma \in \mathbf{C}$. Si l'on note $c(q) = q^{\gamma}$, on a, sur 
$\Omega$: $\underset{\epsilon \to 0^{+}}{\lim} e_{q,c(q)}(z) = (-z)^{\gamma}$
(d\'etermination principale, telle que $-1 \mapsto 1$).
}\\

Cela d\'ecoule imm\'ediatement du th\'eor\`eme. On peut toutefois y ajouter 
une interpr\'etation g\'eom\'etrique.\\ 

Le lieu polaire de
$e_{q,c(q)}$ est la spirale discr\`ete $c(q) q^{\mathbf{Z}}$, c'est \`a dire
l'image par $x \mapsto e^{2 \imath \pi x}$ du r\'eseau translat\'e
$(\mathbf{Z} + \mathbf{Z} \tau) + \gamma \epsilon$. Lorsque
$\epsilon \to 0^{+}$, celui-ci se condense en le sous-groupe ferm\'e
$\mathbf{Z} + \mathbf{R} \tau_{0}$ de $\mathbf{C}$, dont l'image dans
$\mathbf{C}^{*}$ est la spirale logarithmique continue $q_{0}^{\mathbf{R}}$:
\emph{la barri\`ere de p\^oles (discr\`ete) se condense en une coupure.}\\

En ajoutant \`a ce lemme la discussion du
3.1.7 sur le comportement au premier ordre, il vient alors:\\

\textbf{\underline{Proposition}}\\

\emph{ 
On suppose donn\'es des exposants $c(q)$ d\'ependant de $q$.
Si $\underset{\epsilon \to 0^{+}}{\lim}(c(q) - 1)/(q - 1) = \gamma$, 
on a, sur $\Omega$, l'\'egalit\'e: 
        $\underset{\epsilon \to 0^{+}}{\lim} e_{q,c(q)}(z) = (-z)^{\gamma}$
}
\hfill $\Box$\\

L'interpr\'etation g\'eom\'etrique est ici la m\^eme, le r\'eseau translat\'e
$(\mathbf{Z} + \mathbf{Z} \tau) + \gamma \epsilon$ \'etant simplement 
remplac\'e par $(\mathbf{Z} + \mathbf{Z} \tau) + \gamma(q) \epsilon$ qui se
condense \'egalement en le sous-groupe ferm\'e
$\mathbf{Z} + \mathbf{R} \tau_{0}$
(ici, $c(q) = q^{\gamma(q)}$ et $\gamma(q) \to \gamma$).


\subsection{Logarithmes}

On a introduit au chapitre 1 la fonction 
$l_{q}(z) = z \frac{\Theta_{q}'(z)}{\Theta_{q}(z)}$, qui nous tient lieu
d'analogue du logarithme. On v\'erifie ici qu'elle en est en effet (\`a une
constante multiplicative pr\`es) une d\'eformation.\\

\textbf{\underline{Proposition}}\\

\emph{On a, sur $\Omega$, l'\'egalit\'e: 
$\underset{\epsilon \to 0^{+}}{\lim} (q-1)l_{q}(z) = \log(-z)$
(d\'etermination principale, telle que $-1 \mapsto 0$).
}\\

\textbf{\underline{Preuve}}\\

Avec toujours les m\^emes notations concernant $q$, $\tau$, $z$ et $x$, on 
introduit $\tilde{L}_{\tau}(x) = l_{q}(z)$. La d\'erivation logarithmique de
$\tilde{\Theta}(\tau,x)$ par rapport \`a $x$ donne:
$$
\tilde{L}_{\tau}(x) = 
\frac{1}{2 \imath \pi} 
\frac {\tilde{\Theta}'_{x} (\tau,x)}{\tilde{\Theta} (\tau,x)}
$$

On va encore avoir recours \`a une \'equation fonctionnelle. On d\'erive
logarithmiquement celle obtenue au 3.1.2. On obtient alors:
$$
-\tau \tilde{L}_{\tau}(x) = 
x - \frac{\tau + 1}{2} + \tilde{L}_{-1/\tau}(-x/\tau)
$$

Lorsque $\epsilon \rightarrow 0^{+}$, $q - 1 \sim - 2 \imath \pi \tau$ et
$(q-1)l_{q}(z)$ est \'equivalent \`a
$$
- 2 \imath \pi \tau \tilde{L}_{\tau}(x) =
2 \imath \pi (x - \frac{\tau + 1}{2}) + 
2 \imath \pi \tilde{L}_{-1/\tau}(-x/\tau)
$$

On choisit 
$x \in ]0;1[ + \mathbf{R} \tau_{0}$, ce qui est loisible puisque $z \in \Omega$
et que $x$ n'est d\'etermin\'e qu'\`a un entier pr\`es. Le terme
$2 \imath \pi (x - \frac{\tau + 1}{2})$ tend vers 
$2 \imath \pi (x - \frac{1}{2})$, qui vaut $\log(-z)$ dans la d\'etermination
choisie. Reste \`a voir que, pour ce choix de $x$,
$\tilde{L}_{-1/\tau}(-x/\tau) \underset{\epsilon \to 0^{+}}{\to} 0$.\\

L'expression de $l_{q}$ comme somme infinie au 1.1.3 (d\'eduite en d\'erivant
logarithmiquement la formule du triple produit de Jacobi) montre que

$$
\tilde{L}_{-1/\tau}(-x/\tau) =
\sum_{r \geq 1} 
   \frac{e^{2 \imath \pi \frac{x - r}{\tau}}}
        {1 - e^{2 \imath \pi \frac{x - r}{\tau}}} -
\sum_{r \geq 0} 
   \frac{e^{2 \imath \pi \frac{- x - r}{\tau}}}
        {1 - e^{2 \imath \pi \frac{- x - r}{\tau}}}
$$

Les termes $\frac{x - r}{\tau}$ ($r \geq 1$) et 
$\frac{- x - r}{\tau}$ ($r \geq 0$) sont tous de la forme 
$\frac{a}{\tau} + $ un r\'eel, o\`u $a < 0$; les exponentielles 
correspondantes ont donc toutes pour limite $0$. La convergence \'etant normale
sur tout compact de $\Omega$, on obtient bien la limite voulue.
\hfill $\Box$\\

Cette proposition (ainsi que les corollaires 4 et 5 des deux lemmes de 3.1.7)
entraine en particulier qu'il faudra ``renormaliser''
les parties log-car des solutions canoniques, en les exprimant en fonction de
$(q-1)l_{q}$ plut\^ot que de $l_{q}$ lui-m\^eme, si l'on veut un comportement
r\'egulier lorsque $q \to 1$ (ceci ne sera d'ailleurs pas n\'ecessaire avec
les solutions du cas non r\'esonnant b\^aties \`a l'aide des matrices
$e_{q,A}$).


\subsection{Comportement de $\Theta_{q}^{+}(q^{\alpha}z)/\Theta_{q}^{+}(z)$}

On en aura besoin en 4.4 pour d\'eformer $(1 - z/z_{0})^{\alpha}$. Pour 
l'obtenir, on proc\`edera comme suit:

\begin{itemize}

\item{Cas o\`u $|z| < 1$: estimation individuelle des coefficients, puis
invocation de la convergence normale.}

\item{Cas o\`u $|z| > 1$: exploitation du lien entre 
$\Theta_{q}^{+}(z) \Theta_{q}^{+}(1/z)$ et $\Theta_{q}(z)$.}

\item{Cas o\`u $|z| = 1$: majoration directe de l'\'ecart \`a un $w$ tel que
$|w| \not= 1$.}

\end{itemize}

On trouve dans \cite{GR} (page 9), sans preuve et sans r\'ef\'erence, la 
r\'eponse lorsque $q$ est \`a valeurs r\'eelles, mais cela ne nous suffit pas.\\

On notera (dans cette section seulement): 
$f_{q,\alpha}(z) = \frac{\Theta_{q}^{+}(q^{\alpha}z)}{\Theta_{q}^{+}(z)}$.
C'est une fonction holomorphe sur $\mathbf{C} - q^{\mathbf{N}}$, admettant en
g\'en\'eral des p\^oles simples en les $q^{n} \;,\; n \in \mathbf{N}$. Le cas
sp\'ecial est celui o\`u $\alpha \in \mathbf{Z}$: les \'eliminations font alors
de cette fonction une fraction rationnelle de degr\'e fixe (ind\'ependant de
$q$); les r\'esultats qui suivent sont triviaux dans ce cas.\\

\textbf{\underline{Th\'eor\`eme}}\\

\emph{$\underset{\epsilon \to 0^{+}}{\lim} f_{q,\alpha}(z) = (1 - z)^{\alpha}$,
cette derni\`ere fonction \'etant la d\'etermination principale, telle que
$0 \mapsto 1$, sur l'ouvert simplement connexe $\Omega_{0}$.}\\

\textbf{\underline{Preuve}}\\

La d\'efinition de $\Theta_{q}^{+}(z)$ donne le d\'eveloppement en produit
infini:
$$
f_{q,\alpha}(z) = 
\prod_{r \geq 0} \frac{1 - q^{\alpha - r} z}{1 - q^{- r} z}
$$

On voit qu'il admet un d\'eveloppement en s\'erie enti\`ere de rayon de 
convergence $1$ (module de la plus proche singularit\'e), disons
$\sum_{m \geq 0} a_{m} z^{m}$.\\ 

De l'\'equation fonctionnelle
$$
\sigma_{q} f_{q,\alpha} = \frac{1 - q^{\alpha + 1} z}{1 - q z} f_{q,\alpha}
$$
on d\'eduit les relations de r\'ecurrence:
$$
a_{0} = 1 \text{  et  } 
a_{m} = a_{m-1} \frac{q^{m} - q^{\alpha + 1}}{q^{m} - 1} \; (m \geq 1)
$$
On en tire donc la valeur exacte de $a_{m}$ et sa limite lorsque $q \to 1$:
$$a_{m} = \prod_{j = 1}^{m} \frac{q^{j} - q^{\alpha + 1}}{q^{j} - 1}
$$
$$
\underset{\epsilon \to 0^{+}}{\lim} a_{m} = 
\prod_{j = 1}^{m} \frac{j - (\alpha + 1)}{j} = 
(-1)^{m} \left( \begin{array}{c} \alpha \\ m \end{array} \right)
$$

Ceci montre d\'ej\`a que $f_{q,\alpha} \to (1 - z)^{\alpha}$ coefficient par
coefficient, donc terme \`a terme (pour un $z$ donn\'e). On va montrer que 
cette convergence est domin\'ee sur tout compact de 
$\overset{\circ}{D}(0,1)$, mettons sur $\overline{D}(0,\rho)$ o\`u
$\rho < 1$. Pour cela, on majore les coefficients:
$$
\left|\frac{a_{m}}{a_{m-1}}\right| = 
\left|1 - \frac{q^{\alpha + 1} - 1}{q^{m} - 1}\right| \leq
1 + \left|\frac{q^{\alpha + 1} - 1}{q^{m} - 1}\right| =
1 + \left|\frac{q^{\alpha + 1} - 1}{q - 1}\right| \div 
\left|\frac{q^{m} - 1}{q - 1}\right|
$$

Pour $q$ voisin de $1$, $|\frac{q^{\alpha + 1} - 1}{q - 1}|$ est uniform\'ement
major\'e; d'autre part, le corollaire du lemme 2 de 3.1.1 dit que
$|\frac{q^{m} - 1}{q - 1}| \geq m \frac{Im(\tau_{0})}{|\tau_{0}|}$ et donc que
$|\frac{a_{m}}{a_{m-1}}| \leq 1 + \frac{C}{m}$, o\`u $C$ ne d\'epend ni de $q$
ni de $m$. Alors, 
$|a_{m}| \leq b_{m} := \prod_{1 \leq j \leq m}(1 + \frac{C}{j})$ et 
$f_{q,\alpha}$ est domin\'ee par $\sum b_{m} z^{m}$, dont le rayon de
convergence est $1$. Ceci ach\`eve le cas o\`u $|z| < 1$.\\

Pour passer au cas o\`u $|z| > 1$, on introduit
$$
f_{q,\alpha}(z)f_{q,-\alpha}(z^{-1}) = 
\frac{\Theta_{q}^{+}(q^{\alpha}z)}{\Theta_{q}^{+}(z)}
\frac{\Theta_{q}^{+}(q^{-\alpha}z^{-1})}{\Theta_{q}^{+}(z^{-1})} =
\frac{1 - q^{\alpha}z^{-1}}{1 - z^{-1}} \frac{\Theta_{q}(q^{\alpha}z)}{\Theta_{q}(z)
}$$
Rappelons en effet que, d'apr\`es la formule du triple produit de Jacobi
(donn\'ee au 1.1.1),
$$
\Theta_{q}^{+}(z) \Theta_{q}^{+}(z^{-1}) = 
\frac{1 - z^{-1}}{(p;p)_{\infty}} \Theta_{q}(z)
$$

Utilisant 4.1.2 et le cas $|z| < 1$, on trouve que, pour $z \in \Omega_{0}$ et
$|z| > 1$:
$$
\underset{\epsilon \to 0^{+}}{\lim} f_{q,\alpha}(z) =
\frac{(-z)^{\alpha}}{(1 - z^{-1})^{- \alpha}}
$$
et il n'est pas tr\`es difficile de v\'erifier que ceci vaut bien 
$(1 - z)^{\alpha}$, achevant la preuve dans le second cas.\\

Soit maintenant $z$ de module $1$; la condition $z \in \Omega$ signifie que
$z \not= 1$. On fixe, selon le lemme 1 de 4.1.1, des r\'eels
$\rho , \delta > 0$ tels que
$$
\forall w \in \overset{\circ}{D}(z,\rho) \;,\; \forall \lambda \geq 0 \;:\;
\left|1 - \frac{w}{q_{0}^{\lambda}}\right| \geq \delta  \text{  et  }
\left|1 - \frac{q^{\alpha}w}{q_{0}^{\lambda}}\right| \geq \delta
$$
Il suffit pour cela de restreindre $q$ \`a un voisinage de $1$ d\'efini par
$\epsilon \leq \epsilon_{0}$, puis de remplacer le $\rho$ donn\'e par le lemme
par $\rho - \underset{\epsilon \leq \epsilon_{0}}{\sup}|q^{\alpha} - 1|$
($\epsilon_{0} > 0$ \'etant fix\'e tel que l'on obtienne $\rho > 0$).\\

Notons temporairement, pour $r \in \mathbf{N}$,
$$
P_{r}(w) = \frac{1 - q^{\alpha - r}w}{1 - q^{-r}w} \text{  de sorte que  }
P_{r}'(w) = \frac{(1 - q^{\alpha})q^{-r}}{(1 - q^{-r}w)^{2}}
$$ 

Tous les $P_{r}$ et $P_{r}'$ sont donc d\'efinis et non nuls sur 
$\overset{\circ}{D}(z,\rho)$. De plus, on peut majorer la d\'eriv\'ee:
$$
|P_{r}'(w)| \leq |1 - q^{\alpha}| |q|^{-r} \frac{1}{\delta^{2}}
$$
et appliquer l'in\'egalit\'e des accroissements finis:
$$
|P_{r}(w) - P_{r}(z)| \leq 
|1 - q^{\alpha}| |q|^{-r} \frac{\rho}{\delta^{2}}
$$

Par ailleurs,
$$
|P_{r}(z)| = 
\left|1 + \frac{q^{-r} z (1 - q^{\alpha})}{1 - q^{-r}z}\right| \geq
1 - \frac{|q|^{-r} |z| |1 - q^{\alpha}|}{\delta}
$$

On suppose maintenant $q$ assez voisin de $1$ pour que
$\frac{|z||1 - q^{\alpha}|}{\delta^{2}} \leq \frac{1}{2}$, d'o\`u les
majorations:
$$
\left|\log \frac{f_{q,\alpha}(w)}{f_{q,\alpha}(z)}\right| \leq
\sum_{r = 0}^{\infty} \left|\log \frac{P_{r}(w)}{P_{r}(z)}\right| \leq
\sum_{r = 0}^{\infty} \left|\frac{P_{r}(w)}{P_{r}(z)} - 1\right| \leq
\sum_{r = 0}^{\infty} \frac{2 \rho}{\delta^{2}} |1 - q^{\alpha}||q|^{-r
}
$$

Cette derni\`ere somme est un $O(\rho)$ uniforme vis \`a vis de $q$ puisque le
facteur vaut 
$\frac{2}{\delta^{2}} \frac{|q| |1 - q^{\alpha}|}{|q| - 1}$ qui est born\'e.
On peut alors prendre $\rho$ petit et $w \in \overset{\circ}{D}(z,\rho)$
de module $\not= 1$, donc tel que $f_{q,\alpha}(w) \to (1 - w)^{\alpha}$: ceci
permet aussi de conclure pour $z$ et ach\`eve la preuve du th\'eor\`eme 
dans le dernier cas.
\hfill $\Box$


\subsection{La fonction $l_{q}^{+}(z) = 
                         z(\Theta_{q}^{+})'(z)/\Theta_{q}^{+}(z)$}

Elle a \'et\'e introduite au 1.1.3. On s'en servira au 4.4 pour d\'eformer 
$\log(1 - z)$.\\

\textbf{\underline{Th\'eor\`eme}}\\

\emph{Pour $z \in \Omega_{0}$,
$\underset{\epsilon \to 0^{+}}{\lim} (q - 1) l_{q}^{+}(z) = \log(1 - z)$
(d\'etermination principale sur $\Omega_{0}$, telle que $0 \mapsto 0$).}\\

\textbf{\underline{Preuve}}\\

La strat\'egie est la m\^eme qu'\`a la section pr\'ec\'edente. Supposons donc
pour commencer que $|z| < 1$, plus pr\'ecis\'ement que 
$z \in \overline{D}(0,\rho)$, $\rho < 1$ fix\'e. Alors
$$
l_{q}^{+}(z) = \sum_{r = 0}^{\infty} \frac{- q^{-r} z}{1 - q^{-r} z} =
- \sum_{r \geq 0} \sum_{m \geq 1} (q^{-r} z)^{m} =
- \sum_{m \geq 1} (\sum_{r \geq 0} q^{-m r}) z^{m} =
- \sum_{m \geq 1} \frac{q^{m}}{q^{m} - 1} z^{m}
$$
En effet, la famille double des $q^{- m r} z^{m}$ est domin\'ee par celle des
$|q|^{- m r} \rho^{m}$, laquelle est sommable, ce qui justifie les 
commutations.\\

Comme $\frac{q^{m}}{q^{m} - 1} z^{m} = z^{m} + \frac{1}{q^{m} - 1} z^{m}$, on
a l'expression plus agr\'eable:
$$
l_{q}^{+}(z) = - \frac{z}{1 - z} - \sum_{m \geq 1} \frac{1}{q^{m} - 1} z^{m}
$$
Ainsi, nous sommes conduits \`a \'etudier la s\'erie enti\`ere:
$$
\sum_{m \geq 1} \frac{q - 1}{q^{m} - 1} z^{m} =
- (q - 1) l_{q}^{+}(z) - (q - 1) \frac{z}{1 - z}
$$
Ses coefficients tendent bien vers les coefficients $\frac{1}{m}$ de 
$- \log(1 - z)$. Reste \`a voir que la convergence est domin\'ee sur
$\overline{D}(0,\rho)$. Mais:
$$
\left|\frac{q - 1}{q^{m} - 1}\right| \leq 
\frac{|\tau_{0}|}{Im(\tau_{0})} \frac{1}{m}
$$
d'apr\`es le corollaire du lemme 2 de 4.1.1. Comme la s\'erie
$\sum_{m \geq 1} \frac{\rho^{m}}{m}$ converge, ce premier cas est r\'egl\'e.\\

Le cas $|z| > 1$ vient alors facilement de la formule:
$$
l_{q}^{+}(z) - l_{q}^{+}(z^{-1}) = l_{q}(z) - \frac{1}{1 - z}
$$
et du comportement de $(q - 1)l_{q}$, \'etudi\'e en 4.1.4: il n'est en effet
pas tr\`es difficile de v\'erifier que
$$
\log(1 - z) - \log(1 - z^{-1}) = \log(-z)
$$

Reste le cas de $z$ de module $1$ (mais diff\'erent de $1$). Soient
$\rho \;,\; \delta > 0$ tels que
$$
\forall w \in \overset{\circ}{D}(z,\rho) \;,\; \forall \lambda \geq 0 \;:\;
\left|1 - \frac{w}{q_{0}^{\lambda}}\right| \geq \delta
$$
Leur existence est assur\'ee par le lemme 1 de 4.1.1. La convergence uniforme
d'une s\'erie de fonctions holomorphes autorisant la d\'erivation terme \`a
terme, on peut majorer, pour $w \in \overset{\circ}{D}(z,\rho)$:
$$
|(l_{q}^{+})'(w)| \leq 
\sum_{r \geq 0} \left|\frac{q^{-r}}{(1 - q^{-r}w)^{2}}\right| \leq
\frac{1}{\delta^{2}} \frac{|q|}{|q| - 1}
$$

Par application de l'in\'egalit\'e des accroissements finis:
$$
|(q - 1) l_{q}^{+}(w) - (q - 1) l_{q}^{+}(z)| \leq
\frac{\rho |q|}{\delta^{2}} \frac{|q - 1|}{|q| - 1}
$$
Ce dernier majorant est un $O(\rho)$ uniforme en $q$. On choisit donc $w$ 
proche de $z$ et de module $\not= 1$, et l'utilisation du m\^eme argument 
qu'\`a la section pr\'ec\'edente permet d'achever la d\'emonstration.
\hfill $\Box$


\subsection{Le comportement ne d\'epend que du premier ordre}

Les r\'esultats de 3.1.2 et de 3.1.5 supposent que la variable est de la forme 
$q^{\alpha} z$, c'est \`a dire varie exactement le long d'une spirale
logarithmique avec $q$. On va ici montrer que cette condition peut \^etre
assouplie, ce qui sera commode (et m\^eme indispensable) pour d\'eformer des
solutions d'\'equations diff\'erentielles (lin\'eaires, \`a coefficients
rationnels et fuchsiennes) arbitrairement prescrites. J'ai d\'ej\`a 
conjectur\'e au d\'ebut de ce chapitre que la variation de $q$ elle-m\^eme 
devait pouvoir \^etre assouplie de fa\c{c}on analogue.\\

Comme dans tout ce chapitre, on n'indique pas toujours explicitement dans les
notations les d\'ependances en $\epsilon$. En particulier, on va consid\'erer
des complexes $c$, $z$, $z_{1}$ et $z_{2}$ qui d\'ependent de $q$, donc de $\epsilon$ (alors que $z_{0}$ et $\alpha$ seront fix\'es).\\

\textbf{\underline{Lemme 1}}\\

\emph{On suppose que $z_{1}$ et $z_{2}$ tendent vers $z_{0} \in \Omega_{0}$.
On suppose de plus que $z_{2}/z_{1} = 1 + o(\epsilon)$. Alors
$$
\underset{\epsilon \to 0^{+}}{\lim} 
\frac{\Theta_{q}^{+}(z_{2})}{\Theta_{q}^{+}(z_{1})} = 1
$$}

\textbf{\underline{Lemme 2}}\\

\emph{
(i) Sous les m\^emes hypoth\`eses,
$$
\underset{\epsilon \to 0^{+}}{\lim} (l_{q}^{+}(z_{2}) - l_{q}^{+}(z_{1})) = 0
$$
(ii) Si l'on suppose seulement que $z_{2}/z_{1} = 1 + O(\epsilon)$,
$$
\underset{\epsilon \to 0^{+}}{\lim} ((q-1)l_{q}^{+}(z_{2}) - (q-1)l_{q}^{+}(z_{1})) = 0
$$
}

Avant de prouver ces deux lemmes, notons tout de suite leurs cons\'equences
qui nous seront utiles pour \'etendre la port\'ee des estimations 
pr\'ec\'edentes:\\

\textbf{\underline{Corollaire 1}}\\

\emph{Soit $z_{0} \in \Omega$. Si, pour $j = 1 \;,\; 2$, 
$z_{j} = z_{0} (1 - 2 \imath \pi \tau_{0} \alpha_{j} \epsilon  + o(\epsilon))$,
alors 
$$
\underset{\epsilon \to 0^{+}}{\lim} 
\frac{\Theta_{q}(z_{2})}{\Theta_{q}(z_{1})} =
(- z_{0})^{\alpha_{2} - \alpha_{1}}
$$}
L'hypoth\`ese dit que $\frac{z_{j}/z_{0} - 1}{q - 1} \to \alpha_{j}$.\\

\newpage

\textbf{\underline{Corollaire 2}}\\

\emph{Soit $z_{0} \in \Omega$. Si $c = 1 - 2 \imath \pi \tau_{0} \gamma \epsilon  + o(\epsilon)$ et si
$z = z_{0} + o(\epsilon)$, alors 
$$
\underset{\epsilon \to 0^{+}}{\lim} e_{q,c}(z) = (-z_{0}) ^{\gamma}
$$}
L'hypoth\`ese dit que $\frac{c - 1}{q - 1} = \gamma$. Ce corollaire justifie
(enfin!) la proposition de 3.1.3.\\

\textbf{\underline{Corollaire 3}}\\

\emph{Soit $z_{0} \in \Omega_{0}$. Si, pour $j = 1 \;,\; 2$, 
$z_{j} = z_{0} (1 - 2 \imath \pi \tau_{0} \alpha_{j} \epsilon  + o(\epsilon))$,
alors
$$\underset{\epsilon \to 0^{+}}{\lim} \Theta_{q}^{+}(z_{2})/\Theta_{q}^{+}(z_{1}) =
(1 - z_{0})^{\alpha_{2} - \alpha_{1}}
$$}

\textbf{\underline{Corollaire 4}}\\

\emph{Soit $z_{0} \in \Omega$. Si $z = z_{0} + o(\epsilon)$, alors 
$$
\underset{\epsilon \to 0^{+}}{\lim} (q-1) l_{q}(z) = \log(-z_{0})
$$}
L'hypoth\`ese dit que $\frac{z/z_{0} - 1}{q - 1} \to 0$.\\

\textbf{\underline{Corollaire 5}}\\

\emph{Soit $z_{0} \in \Omega_{0}$. Si $z = z_{0} + o(\epsilon)$, alors 
$$
\underset{\epsilon \to 0^{+}}{\lim} (q-1) l^{+}_{q}(z) = \log(1-z_{0})
$$}

\textbf{\underline{Preuve du lemme 1}}\\

Soit $u_{q} = z_{1} - z_{2}$: c'est un $o(\epsilon)$.
On part de l'\'egalit\'e:
$$
\frac{\Theta_{q}^{+}(z_{2})}{\Theta_{q}^{+}(z_{1})} =
\prod_{r \geq 0} (1 + \frac{q^{-r}u_{q}}{1 - q^{-r}z_{1}})
$$
Il existe un disque
$\overset{\circ}{D}(z_{0},\rho)$ sur lequel, quel que soit $\lambda \geq 0$,
$|1 - \frac{z}{q_{0}^{\lambda}}| \geq \delta > 0$. Pour $q$ proche de $1$,
$z_{1}$ et $z_{2}$ sont dans le disque et
$$
\left|\log \frac{\Theta_{q}^{+}(z_{2})}{\Theta_{q}^{+}(z_{1})}\right| \leq
\sum_{r \geq 0} \left|\frac{q^{-r} u_{q}}{1 - q^{-r} z_{1}}\right| \leq
\frac{|u_{q}|}{\delta}\frac{|q|}{|q| - 1}
$$
et ce dernier majorant est un $O(\frac{u_{q}}{q - 1})$ uniform\'ement en $q$.\\

Comme $u_{q} = o(\epsilon)$ et $q - 1 \sim - 2 \imath \pi \tau_{0} \epsilon$,
cette estimation est suffisante.
\hfill $\Box$\\

\textbf{\underline{Preuve du lemme 2}}\\

Avec les m\^emes hypoth\`eses que ci-dessus, on part de l'in\'egalit\'e :
$$
|l_{q}^{+}(z_{2}) - l_{q}^{+}(z_{1})| \leq
\sum_{r \geq 0} 
\left|\frac{q^{-r} u_{q}}{(1 - q^{-r} z_{1})(1 - q^{-r} z_{2})}\right| \leq
\frac{|q|}{\delta^{2}}\frac{|q - 1|}{|q| - 1}\frac{|u_{q}|}{|q - 1}
$$

On trouve donc encore une fois que 
$|l_{q}^{+}(z_{2}) - l_{q}^{+}(z_{1})| = O(\frac{u_{q}}{q - 1})$ 
uniform\'ement en $q$. Nous pouvons donc conclure
pour les deux cas cit\'es dans le lemme.
\hfill $\Box$



\section{Comportement \`a valeur initiale fixe}

Selon les conventions d\'etaill\'ees au 3.1,, $q = q_{0}^{\epsilon}$ tend vers 
$1$ le long d'une spirale. Les d\'ependances en $q$ ne seront pas 
n\'ecessairement explicites dans les notations.\\

On consid\`ere l'\'equation aux $q$-diff\'erences \`a coefficients rationnels
fuchsienne en $0$, ``avec conditions initiales'':
$$
\begin{cases} \sigma_{q}X = A X \\ X(0) = X_{0} \end{cases}
$$

L'\'equation est vue ``en famille'', c'est \`a dire que $A$ d\'epend de $q$ (et
donc de $\epsilon$). La condition initiale est, dans cette section,
fixe, c'est \`a dire que $X_{0}$ ne d\'epend pas de $q$. De plus:

\begin{itemize}

\item{$A$ satisfait les hypoth\`eses g\'en\'erales qui font que l'\'equation
est fuchsienne; en particulier, $A_{0} = A(0) \in GL_{n}(\mathbf{C})$. 
Ces hypoth\`eses seront renforc\'ees selon les
besoins au fur et \`a mesure de l'\'etude.
On veut faire confluer ce syst\`eme vers le syst\`eme 
diff\'erentiel \`a coefficients rationnels avec conditions initiales
$$
\begin{cases} \delta \tilde{X} = A \tilde{X} \\ 
\tilde{X}(0) = X_{0} \end{cases}
$$

L'analogie choisie \'etant 
$(\sigma_{q} - Id)/(q-1) \leftrightarrow \delta = z d/dz$, nous supposerons que
la matrice $B = (A - I_{n})/(q-1)$ tend, lorsque $\epsilon \to 0^{+}$,
vers la  matrice rationnelle $\tilde{B} \in M_{n}(\mathbf{k})$. On \'ecrira
\'egalement $A = I_{n} + \eta B \;,\; \eta = q - 1$. La condition de
rationnalit\'e de $\tilde{B}$ peut \^etre affaiblie (dans ce chapitre) en 
$\tilde{B} \in M_{n}(\mathbf{k}_{0})$, mais les hypoth\`eses sur les p\^oles
seraient nettement plus compliqu\'ees \`a formuler. Notons d'ailleurs que nous
\'etudierons au chapitre 4.4 un tel exemple, avec une singularit\'e
essentielle \`a l'infini: nos m\'ethodes s'appliqueront encore, ce qui montre
que les \'equations aux $q$-diff\'erences que nous appelons fuchsiennes ne
sont peut-\^etre pas si fuchsiennes que \c{c}a.}

\item{$X_{0}$ est un vecteur fix\'e non nul de $\mathbf{C}^{n}$ tel que, 
quelque soit $q$, $A_{0}X_{0} = X_{0}$. On suppose de plus qu'il n'y a pas de
r\'esonnance: aucun $q^{m} \;,\; m \in \mathbf{N}^{*}$ n'est valeur propre
de $A_{0}$. Cela sera cons\'equence d'une autre hypoth\`ese sur $\tilde{B}$
\`a partir du 3.2.1.}

\end{itemize}

Les algorithmes du chapitre $1$ attribuent au syst\`eme aux $q$-diff\'erences 
avec condition initiale ci-dessus une unique solution $X$. La th\'eorie
classique des \'equations diff\'erentielles montre l'existence d'une unique 
solution du syst\`eme diff\'erentiel avec conditions initiales. Pour prouver
la convergence de $X$ vers $\tilde{X}$ (qui est en fait une \emph{confluence}),
nous proc\`ederons en deux \'etapes:

\begin{itemize}

\item{Etude locale: il s'agit de reprendre la r\'esolution formelle et la
preuve de convergence ``en famille'', et, en parall\`ele, l'algorithme 
correspondant pour l'\'equation diff\'erentielle. Il faudra donc \'egalement
une hypoth\`ese de non-r\'esonnance pour celle-ci: \`a partir du 3.2.1, nous
supposons donc que
l'\'equation diff\'erentielle $\delta \tilde{X} = \tilde{B} \tilde{X}$ 
est fuchsienne sur $\mathbf{S}$ et non r\'esonnante en $0$, autrement dit 
$Sp(\tilde{B}(0)) \cap \mathbf{N}^{*} = \emptyset$. On verra que cette 
condition entraine que, pour $\epsilon$ assez petit, $A$ est non r\'esonnante 
en $0$, autrement dit: $Sp(A(0)) \cap q^{\mathbf{N}^{*}} = \emptyset$.}

\item{Etude globale: on sort d'un voisinage de $0$ (dans lequel la convergence
a \'et\'e \'etablie) le long de demi-spirales $q_{0}^{\mathbf{R}^{+}}z_{0}$. Ceci
nous ram\`ene \`a un probl\`eme de variable r\'eelle (on pose 
$z = z_{0}q_{0}^{t}$) et \'equivaut \`a l'emploi de la m\'ethode d'Euler pour
approcher les solutions d'une \'equation diff\'erentielle.
}

\end{itemize}

Pour que l'ouvert o\`u cela marche, et donc le r\'esultat, soit non trivial, 
il faudra en outre introduire une hypoth\`ese sur les p\^oles.
Notons $\tilde{z}_{1},...,\tilde{z}_{r}$ les p\^oles de $\tilde{B}$ sur 
$\mathbf{S}$. Ce sont donc en fait des \'el\'ements de $\mathbf{C}^{*}$ 
(l'\'equation diff\'erentielle \'etant fuchsienne). 
L'hypoth\`ese n\'ecessaire est que ``les p\^oles de $A$ (et donc ceux de $B$) 
tendent vers ceux de
$\tilde{B}$'', c'est \`a dire que, quelque soit 
$z \not\in \{\tilde{z}_{1},...,\tilde{z}_{r}\}$, $A$ soit d\'efinie en $z$ pour
$\epsilon$ assez petit. 

Notons alors 
$$
\tilde{\Omega}_{0} =
\mathbf{C} - \bigcup_{1 \leq j \leq r} \tilde{z}_{j} q_{0}^{\mathbf{R}^{+}} =
\mathbf{S} - \bigcup_{1 \leq j \leq r} \tilde{z}_{j} q_{0}^{\overline{\mathbf{R}^{+}}}
$$
C'est donc un ouvert simplement connexe de $\mathbf{S}$. Si $\Omega_{q,0}$ est 
l'ouvert analogue pour les p\^oles de $A$, on voit que $\tilde{\Omega}_{0}$ 
est inclus dans la ``limite'' 
$\underset{\alpha > 0}{\bigcup} \underset{\epsilon \in ]0;\alpha[}{\bigcap}\Omega_{q,0}$. Ceci permettra de donner un sens \`a l'affirmation: ``$B$ converge vers
$\tilde{B}$ sur $\tilde{\Omega}_{0}$'' et \`a des conditions de convergence
uniforme.


\bigskip
\hrule
\bigskip

\unitlength = 1cm

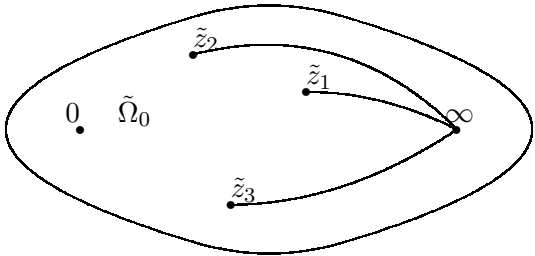
\underline{Figure 2}\\

\bigskip

\begin{picture} (15,3)


\qbezier (4,0)(9,1.5)(4,3)     

\qbezier (2,0)(-3,1.5)(2,3)      

\qbezier (2,3)(3,3.3)(4,3)   

\qbezier (2,0)(3,-0.3)(4,0)  


\put(0.5,1.5){\circle*{0.1}}
\put(0.3,1.6){$0$}

\put(5.5,1.5){\circle*{0.1}}
\put(5.35,1.6){$\infty$}

\put(2,2.5){\circle*{0.1}}
\put(2,2.6){$\tilde{z}_{2}$}

\put(3.5,2){\circle*{0.1}}
\put(3.5,2.1){$\tilde{z}_{1}$}

\put(2.5,0.5){\circle*{0.1}}
\put(2.5,0.6){$\tilde{z}_{3}$}


\qbezier (2,2.5) (4,3) (5.5,1.5)

\qbezier (3.5,2) (4.5,2) (5.5,1.5)

\qbezier (2.5,0.5) (4,0.5) (5.5,1.5)


\put(1,1.6){$\tilde{\Omega}_{0}$}

\end{picture}

\bigskip
\hrule
\bigskip


\subsection{Etude locale}

Outre les hypoth\`eses pr\'ec\'edentes, on suppose que $\tilde{B}$ est 
holomorphe en $0$, ainsi que $B$ pour $\epsilon$ assez petit:
$$
\begin{cases}
\tilde{B} = \tilde{B}_{0} + \tilde{B}_{1}z + ... \\
B = B_{0} + B_{1}z + ...
\end{cases}
$$ 
(les $B_{m}$ d\'ependent donc de $q$),
et que la convergence de $B$ vers $\tilde{B}$ est normale dans un disque
ouvert non vide $\overset{\circ}{D}(0,\rho) \subset \tilde{\Omega}_{0}$: il
existe une s\'erie enti\`ere $b_{0} + b_{1} z + ...$ de rayon de convergence
$\rho > 0$ telle que, pour $\epsilon$ assez petit\footnote{Dans tout ce chapitre, on notera $||-||$ une norme quelconque
sur $\mathbf{C}^{n}$ ainsi que la norme subordonn\'ee sur $M_{n}(\mathbf{C})$,
par exemple 
$||M|| = 
\underset{1 \leq j \leq n}{\max} \underset{1 \leq j \leq n}{\sum} |m_{i,j}|$
}:
$$
\forall m \in \mathbf{N} \;:\; ||B_{m}|| \leq b_{m}
$$

\textbf{\underline{Proposition}}\\

\emph{$X$ converge uniform\'ement vers $\tilde{X}$ sur un disque ouvert non 
vide $\overset{\circ}{D}(0,r) \subset \tilde{\Omega}_{0}$.}\\

Pour le d\'emontrer, nous aurons besoin du\\

\textbf{\underline{Lemme}}\\

\emph{Il existe une constante $C$ telle que\\
(i) $\forall m \geq 1 \;,\; ||(m I_{n} - \tilde{B}(0))^{-1}|| \leq C$\\
(ii) Pour $\epsilon$ suffisamment petit,
$\forall m \geq 1 \;,\; ||(\frac{q^{m}-1}{q-1} I_{n} - B(0))^{-1}|| \leq C$}\\

\textbf{\underline{Preuve du lemme}}\\

Pour $m \geq 1$, $(m I_{n} - \tilde{B}(0))^{-1}$ existe (condition de
non-r\'esonnance); comme sa norme tend vers $0$ pour $m \to \infty$, la 
premi\`ere conclusion est assur\'ee.\\

On a vu au 3.1.1 (corollaire du lemme 2) que
$|\frac{q^{m} - 1} {q - 1}| \geq m \frac{Im(\tau_{0})} {|\tau_{0}|}$. Pour
$\epsilon$ petit, c'est \`a dire pour $q$ proche de $1$, on a donc une 
majoration uniforme de $||(\frac{q^{m}-1}{q-1} I_{n} - B(0))^{-1}||$ valable
pour $m \geq m_{0}$, $m_{0}$ \'etant convenablement choisi (assez grand).
Comme la limite de $\frac{q^{m}-1}{q-1} I_{n} - B(0)$ est 
$m I_{n} - \tilde{B}(0)$ et que
l'on n'a besoin de consid\'erer que le cas o\`u $m < m_{0}$, la seconde
conclusion est \'egalement assur\'ee.
\hfill $\Box$\\

\textbf{\underline{Preuve de la proposition}}\\

La r\'esolution formelle des deux syst\`emes se m\`ene en parall\`ele : leurs
solutions respectives sont $X = X_{0} + zX_{1} + ...$ et
$\tilde{X} = \tilde{X}_{0} + z\tilde{X}_{1} + ...$, o\`u 
$\tilde{X}_{0} = X_{0}$ est une donn\'ee et o\`u les autres termes se calculent
par r\'ecurrence:
$$
\begin{cases}
X_{m} = (q^{m}I_{n} - A_{0})^{-1}(A_{1}X_{m-1} +...+ A_{m}X_{0}) \\
\tilde{X}_{m} = (mI_{n} - \tilde{B}_{0})^{-1}
(\tilde{B}_{1}\tilde{X}_{m-1} +...+ \tilde{B}_{m}\tilde{X}_{0})
\end{cases}
$$
Mais, compte tenu de ce que $A_{0} = I_{n} + \eta B_{0}$ et, pour $m \geq 1$,
$A_{m} = \eta B_{m}$, la premi\`ere relation se r\'e\'ecrit:
$$
X_{m} = \frac{q^{m}-1}{q-1}I_{n} - B_{0})^{-1}(B_{1}X_{m-1} +...+ B_{m}X_{0})
$$

La convergence uniforme sur un disque de la fonction analytique $B$ vers la 
fonction analytique $\tilde{B}$ entraine la convergence terme \`a terme:
$\forall m \in \mathbf{N} \;,\; \tilde{B}_{m} \to B_{m}$. On en d\'eduit par 
r\'ecurrence le fait que, quelque soit $m$, $\tilde{X}_{m} \to X_{m}$. Ceci 
est \'evidemment valable sur $\tilde{\Omega}_{0}$.\\

On applique la m\'ethode des s\'eries majorantes \`a ces deux r\'esolutions
formelles: on introduit 
$\overline{\mathcal{X}}(z) = \overline{x}_{0} + \overline{x}_{1} z + ...$,
o\`u $\overline{x}_{0} = ||X_{0}|| = ||\tilde{X}_{0}||$ et o\`u
$\overline{x}_{m} = C(b_{1} \overline{x}_{m-1} + ... + b_{m} \overline{x}_{0})$
($C$ est la constante fournie par le lemme ci-dessus).\\

La s\'erie enti\`ere $\overline{\mathcal{X}}(z)$ a un rayon de convergence
$\rho' > 0$ (voir la preuve du lemme 1 de 1.2.2). 
Prenant $r = \min(\rho,\rho')$, on conclut que la convergence de $X$ vers
$\tilde{X}$ est domin\'ee sur $\overset{\circ}{D}(0,r)$, plus pr\'ecis\'ement:
$\forall m \in \mathbf{N} \;,\; ||X_{m} z^{m}|| \leq b_{m} r^{m}$, avec
$\underset{m \geq 0}{\sum} b_{m} r^{m} < \infty$. L'invocation du th\'eor\`eme 
de convergence domin\'ee ach\`eve la d\'emonstration.
\hfill $\Box$\\

Cette d\'emonstration contient implicitement la preuve de la non-r\'esonnance 
du syst\`eme au $q$-diff\'erences pour $q$ proche de $1$.


\subsection{Etude globale}

On va maintenant sortir du disque $\overset{\circ}{D}(0,r)$ le long de
segments de spirales de la forme $q_{0}^{[a;b]}$ qui ne rencontrent pas les 
p\^oles $\tilde{z}_{i}$.\\ 

On suppose donc dor\'enavant que, pour tout tel segment de spirale inclus dans
$\tilde{\Omega}_{0}$, $A$ est d\'efini sur tout ce segment pour $q$ assez voisin
de $1$. Cette hypoth\`ese sera subsum\'ee par l'hypoth\`ese C que l'on 
introduira au 3.2.3 (convergence uniforme de $B$ vers $\tilde{B}$ sur tout
tel segment).\\

Fixons un point $z_{1} \in \tilde{\Omega}_{0}$ en lequel on veut d\'emontrer 
que $X \to \tilde{X}$. Pour un certain $\lambda > 0$, 
$z_{0} = q_{0}^{-\lambda}z_{1} \in \overset{\circ}{D}(0,r)$. Si
$z_{1} \in \overset{\circ}{D}(0,r)$, on peut m\^eme prendre $\lambda = 0$, mais
dans ce cas il n'y a rien \`a prouver !\\

Puisque $|z_{0}| < r$, on peut aussi fixer un certain $\alpha > 0$, tel que
$|q_{0}^{\alpha} z_{0}| < r$, de sorte que $q_{0}^{[-\infty;\alpha]}z_{0}$ 
est une partie compacte de $\overset{\circ}{D}(0,r)$.\\

Notons d'autre part $\Sigma$ le lieu polaire de $A$ (et de $B$) et
$\Omega_{q} = \mathbf{C} - q_{0}^{\mathbf{R}} \Sigma$: c'est un ouvert (non connexe
en g\'en\'eral) contenu dans $\Omega_{q,0}$. D'apr\`es les hypoth\`eses sur les
p\^oles, pour $q$ assez proche de $1$, $z_{1} \in \Omega_{q}$ et l'on a m\^eme
en fait l'inclusion 
$q_{0}^{\overline{\mathbf{R}}} z_{1} = q_{0}^{\overline{\mathbf{R}}} z_{0}
\subset \Omega_{q}$. On le supposera d\'esormais.\\

On passe aux variables r\'eelles en posant 
$z = q_{0}^{t} z_{0} \;,\; t \in \mathbf{R}$ et:
$$
\begin{cases}
Y_{\epsilon}(t) = X(q_{0}^{t} z_{0}) \\
\tilde{Y}(t) = \tilde{X}(q_{0}^{t} z_{0}) \\
C_{\epsilon}(t) = \frac{q - 1}{\epsilon} B(q_{0}^{t} z_{0}) \\
\tilde{C}(t) = - 2 \imath \pi \tau_{0} \tilde{B}(q_{0}^{t} z_{0})
\end{cases}
$$

Notre but est donc de
d\'emontrer que $Y_{\epsilon}(\lambda) \to \tilde{Y}(\lambda)$, c'est \`a dire
que $X(z_{1}) \to \tilde{X}(z_{1})$.\\

Les conditions pos\'ees et les conclusions d\'ej\`a obtenues se r\'e\'ecrivent:
$$
\begin{cases}
Y_{\epsilon}(-\infty) = \tilde{Y}(-\infty) = X_{0}  & 
\text{(dont on ne se servira pas)}\\
Y_{\epsilon}(t + \epsilon) = (I_{n} + \epsilon C_{\epsilon}(t)) Y_{\epsilon}(t) &
\text{(l'\'equation fonctionnelle)}\\
\tilde{Y}'(t) = \tilde{C}(t) \tilde{Y}(t) &
\text{(l'\'equation diff\'erentielle)}\\
\underset{\epsilon \to 0^{+}}{\lim} Y_{\epsilon}(t) = \tilde{Y}(t) &
\text{(convergence uniforme sur $[-\infty;\alpha]$, obtenue au 3.2.1)}\\
\underset{\epsilon \to 0^{+}}{\lim} C_{\epsilon}(t) = \tilde{C}(t) &
\text{(hypoth\`ese qui sera renforc\'ee plus loin)}
\end{cases}
$$

On invoquera maintenant la m\'ethode d'Euler pour la r\'esolution approch\'ee
d'\'equations diff\'erentielles lin\'eaires (vectorielles, d'ordre $1$), 
dans sa version multiplicative et \`a pas non constant\footnote{On peut consid\'erer tout ce travail comme une vision g\'eom\'etrique
de la m\'ethode d'Euler: la g\'eom\'etrie \'etant celle des z\'eros, des 
p\^oles et de la ramification des solutions complexes exactes et approch\'ees}.
Faute d'une r\'ef\'erence exactement appropri\'ee \`a nos besoins, nous
invoquerons \cite{Agreg} (probl\`eme d'analyse, III.5 et IV.1).\\

\textbf{\underline{Th\'eor\`eme}}\\

\emph{Soit $\mathcal{A} \in \mathcal{C}(]a;b[,M_{n}(\mathbf{C}))$ et soient
$t_{0} < t_{1}$ dans $]a;b[$. On se donne, pour $p$ entier $\geq 1$, une
subdivision 
$t_{0} = s_{0}^{(p)} < s_{1}^{(p)} < ... < s_{p-1}^{(p)} < s_{p}^{(p)} = t_{1}$
de diam\`etre
$\epsilon_{p} := \underset{0 \leq i \leq p-1}{\max} \Delta s_{i}^{(p)}$
(on note $\Delta s_{i}^{(p)} = s_{i+1}^{(p)} - s_{i}^{(p)}$). On suppose que
$\underset{p \to \infty}{\lim} \epsilon_{p} = 0$. Alors, la r\'esolvante de
l'\'equation diff\'erentielle $\mathcal{X}' = \mathcal{A} \mathcal{X}$ est
donn\'ee par la formule:
$$
U(t_{1},t_{0}) = \underset{p \to \infty}{\lim}
(I_{n} + \mathcal{A}(s_{p-1}^{(p)})\Delta s_{p-1}^{(p)}) ...
(I_{n} + \mathcal{A}(s_{0}^{(p)})\Delta s_{0}^{(p)})
$$}

Appliquons ce th\'eor\`eme pour calculer la r\'esolvante 
$\tilde{U}(\lambda,\mu)$ de
l'\'equation diff\'erentielle $\tilde{Y}' = \tilde{C} \tilde{Y}$: nous aurons
donc $\tilde{Y}(\lambda) = \tilde{U}(\lambda,0) \tilde{Y}(0)$. On subdivise
$[0;\lambda]$ avec le pas $\nu \epsilon$ puis $N$ fois le pas $\epsilon$, o\`u
$N = E(\frac{\lambda}{\epsilon})$ et $\nu = \frac{\lambda}{\epsilon} - N$.
On trouve:
$$
\tilde{U}(\lambda,0) = 
\underset{\epsilon \to 0^{+}}{\lim}
(I_{n} + \epsilon \tilde{C}(\lambda - \epsilon)) ...
(I_{n} + \epsilon \tilde{C}(\lambda - N \epsilon))
(I_{n} + \nu \epsilon \tilde{C}(0))
$$
Mais nous voulons utiliser $C$ et non $\tilde{C}$. Nous remplacerons donc
\cite{Agreg}, III.6 par une version am\'elior\'ee que l'on d\'eduit du\\

\textbf{\underline{Lemme}}\\

\emph{Soient $A_{1},...,A_{m},E_{1},...,E_{m}$ des \'el\'ements de
$M_{n}(\mathbf{C})$ et notons 
$A'_{1} = A_{1} + E_{1}$,...,$A'_{m} = A_{m} + E_{m}$. Alors
$$
||(I_{n} + A_{1})...(I_{n} + A_{m}) - (I_{n} + A'_{1})...(I_{n} + A'_{m})||
\leq (\sum_{i = 1}^{m} ||E_{i}||)
\exp(\sum_{i = 1}^{m} (||A_{i}|| + ||E_{i}||))
$$
}

\textbf{\underline{Preuve}}\\

On pose, pour $t \in [0;1]$, 
$f(t) = (I_{n} + A_{1} + t E_{1})...(I_{n} + A_{m} + t E_{m})$. Alors
$$
f'(t) = \sum_{i = 1}^{m} 
(I_{n} + A_{1} + t E_{1})...E_{i}...(I_{n} + A_{m} + t E_{m})
$$
et l'on tire de l'in\'egalit\'e triangulaire:
$$
||f'(t)|| \leq 
\sum_{i = 1}^{m} ||E_{i}|| \prod_{j \not= i} (1 + ||A_{j}|| + ||E_{j}||)
$$
Cette majoration permet de conclure, puisque $1 + x \leq e^{x}$.
\hfill $\Box$\\

\textbf{\underline{Corollaire}}\\

\emph{Soient $(A_{p,i})$ et $(B_{p,i})$ deux familles de matrices index\'ees
par des couples d'entiers $(p,i)$ tels que $1 \leq i \leq p$. On consid\`ere
les suites de terme g\'en\'eral
$U_{p} = (I_{n} + A_{p,1})...(I_{n} + A_{p,p})$ et
$V_{p} = (I_{n} + B_{p,1})...(I_{n} + B_{p,p})$. On suppose que, lorsque
$p \to \infty$:\\
(i) La suite de terme g\'en\'eral 
$\underset{1 \leq i \leq p}{\sum} ||A_{p,i}||$ est born\'ee.\\
(ii) La suite de terme g\'en\'eral 
$\underset{1 \leq i \leq p}{\sum} ||A_{p,i} - B_{p,i}||$ converge vers $0$.\\
Alors $\underset{p \to \infty}{\lim} ||U_{p} - V_{p}|| = 0$.}\\

\textbf{\underline{Preuve}}\\

On applique le lemme avec $A_{p,i}$ pour $A_{i}$ et $B_{p,i} - A_{p,i}$ pour
$E_{i}$; il vient:
$$
||U_{p} - V_{p}|| \leq 
(\underset{1 \leq i \leq p}{\sum} ||A_{p,i} - B_{p,i}||)
\exp(\underset{1 \leq i \leq p}{\sum} ||A_{p,i}|| +
\underset{1 \leq i \leq p}{\sum} ||A_{p,i} - B_{p,i}||)
$$
La conclusion est imm\'ediate.
\hfill $\Box$\\

Nous utiliserons ce corollaire avec pour matrices $A_{p,i}$ et $B_{p,i}$ les
valeurs respectives de $\tilde{C}$ et $C$ aux points de la subdivision 
introduite plus haut. La condition (i) est v\'erifi\'ee:
$$
\underset{1 \leq i \leq p}{\sum} ||A_{p,i}|| =
\sum_{j = 1}^{N} \epsilon ||\tilde{C}(\lambda - j \epsilon)|| +
\nu \epsilon ||\tilde{C}(0)||
$$
C'est une somme de Riemann, elle converge vers 
$\int_{0}^{\lambda} ||\tilde{C}(t)|| \; dt$ ($\tilde{C}$ est une fonction 
continue); la suite correspondante est donc born\'ee.\\

De m\^eme, pour v\'erifier la condition (ii), on consid\`ere la somme:
$$
\underset{1 \leq i \leq p}{\sum} ||A_{p,i} - B_{p,i}|| =
\sum_{j = 1}^{N} 
\epsilon ||\tilde{C}(\lambda - j \epsilon) - C_{\epsilon}(\lambda - j \epsilon)|| +
\nu \epsilon ||\tilde{C}(0) - C_{\epsilon}(0)||
$$

Majorant uniform\'ement tous les termes, et rempla\c{c}ant
$(N \epsilon + \nu \epsilon)$ par sa valeur $\lambda$, il vient:
$$
\underset{1 \leq i \leq p}{\sum} ||A_{p,i} - B_{p,i}||\leq
\lambda \; \underset{[0;\lambda]}{\sup} ||\tilde{C} - C_{\epsilon}||
$$
Ce dernier majorant aura pour limite $0$, et donc, la condition (ii) sera \`a
coup sur v\'erifi\'ee, si $C_{\epsilon}$ converge uniform\'ement vers
$\tilde{C}$ sur $[0;\lambda]$; on le suppose donc. On voit alors que
$$
\tilde{U}(\lambda,0) = 
\underset{\epsilon \to 0^{+}}{\lim}
(I_{n} + \epsilon C_{\epsilon}(\lambda - \epsilon)) ...
(I_{n} + \epsilon C_{\epsilon}(\lambda - N \epsilon))
(I_{n} + \nu \epsilon C_{\epsilon}(0))
$$

Par ailleurs, it\'erant l'\'equation fonctionnelle satisfaite par 
$Y_{\epsilon}$, on trouve:
$$
Y_{\epsilon}(\lambda) = 
(I_{n} + \epsilon C_{\epsilon}(\lambda - \epsilon)) ...
(I_{n} + \epsilon C_{\epsilon}(\lambda - N \epsilon))
Y_{\epsilon}(\nu \epsilon)
$$
Le produit ci-dessus est, \`a un facteur pr\`es, celui dont la limite est 
$\tilde{U}(\lambda,0)$; le facteur manquant est 
$I_{n} + \nu \epsilon C_{\epsilon}(0)$, qui tend vers $I_{n}$.
La convergence uniforme de $Y_{\epsilon}$ vers $\tilde{Y}$ sur 
$[-\infty;\alpha]$ implique que 
$\underset{\epsilon \to 0^{+}}{\lim}Y_{\epsilon}(\nu \epsilon) = \tilde{Y}(0)$.
Ainsi, la limite de $Y_{\epsilon}(\lambda)$ est 
$\tilde{U}(\lambda,0) \tilde{Y}(0)$, c'est \`a dire $\tilde{Y}(\lambda)$ 
par d\'efinition de la r\'esolvante $\tilde{U}$. Autrement dit :
$\underset{\epsilon \to 0^{+}}{\lim} X(z_{1}) = \tilde{X}(z_{1})$: ceci 
ach\`eve l'\'etude globale.\\

\textbf{\underline{Proposition}}\\

\emph{Les hypoth\`eses : $\tilde{C}$ continue sur $[0;\lambda]$ et:
$C_{\epsilon} \to \tilde{C}$ uniform\'ement sur $[0;\lambda]$ entrainent:
$Y_{\epsilon}(\lambda) \to \tilde{Y}(\lambda)$, c'est \`a dire:
$X(z_{1}) \to \tilde{X}(z_{1})$.}
\hfill $\Box$


\subsection{Synth\`ese des r\'esultats}

On se donne $\tilde{B} \in M_{n}(\mathbf{k})$, matrice d'une \'equation
diff\'erentielle lin\'eaire \emph{fuchsienne}: les p\^oles
$\tilde{z}_{1}$,...,$\tilde{z}_{r}$ de $\tilde{B}$ sur $\mathbf{S}$
sont donc en fait des \'el\'ements de $\mathbf{C}^{*}$.\\

On se donne \'egalement une condition initiale 
$X_{0} \in \mathbf{C}^{n} - \{0\}$ telle que $\tilde{B} X_{0} = 0$.
On suppose que l'\'equation avec condition initiale:
$$
\begin{cases}
\delta \tilde{X} = \tilde{B} \tilde{X}\\
\tilde{X}(0) = X_{0}
\end{cases}
$$
est non-r\'esonnante, c'est \`a dire qu'aucun entier non nul n'est valeur
propre de $\tilde{B}(0)$. Elle admet alors une unique solution formelle,
celle-ci est convergente et admet un unique prolongement  analytique sur
l'ouvert simplement connexe
$$
\tilde{\Omega}_{0} =
\mathbf{C} - \bigcup_{1 \leq j \leq r} \tilde{z}_{j} q_{0}^{\mathbf{R}^{+}} =
\mathbf{S} - \bigcup_{1 \leq j \leq r} \tilde{z}_{j} q_{0}^{\overline{\mathbf{R}^{+}}}
$$

On se donne enfin une matrice $B \in M_{n}(\mathbf{k})$ 
\emph{dont les coefficients
d\'ependent de $q$} et soumise aux conditions:

\begin{itemize}

\item{(A) Les p\^oles de $B$ tendent vers ceux de $\tilde{B}$, au sens suivant:
$\forall z \not\in \{\tilde{z}_{1},...,\tilde{z}_{r}\}$, $B$ est d\'efinie en 
$z$ pour $\epsilon$ assez petit.}

\item{(B) Il y a un disque ferm\'e non trivial
$\overline{D}(0,\rho) \subset \tilde{\Omega}_{0}$ sur lequel $B$ converge
uniform\'ement vers $\tilde{B}$.}

\item{(C) Il y a \'egalement convergence uniforme de $B$ vers $\tilde{B}$ sur
tout segment de spirale $q_{0}^{[a;b]}z_{0} \subset \tilde{\Omega}_{0}$.}

\end{itemize}

Les conditions (B) et (C) sont des renforcements de (A); (C) garantit
l'hypoth\`ese faite au d\'ebut de l'\'etude globale.\\

Si $z_{1}$,...,$z_{s}$ sont les p\^oles de $B$, on peut encore d\'efinir un
ouvert simplement connexe (d\'ependant de $q$)
$$
\Omega_{q,0} =
\mathbf{C} - \bigcup_{1 \leq j \leq s} z_{j} q_{0}^{\mathbf{R}^{+}} =
\mathbf{S} - \bigcup_{1 \leq j \leq s} z_{j} q_{0}^{\overline{\mathbf{R}^{+}}}
$$
D'apr\`es les hypoth\`eses, tout ouvert relativement compact de 
$\tilde{\Omega}_{0}$ est contenu dans $\Omega_{q,0}$ pour $q$ assez proche de 
$1$.\\

\textbf{\underline{Th\'eor\`eme}}\\

\emph{On pose $\eta = q - 1$ et $A = I_{n} + \eta B$. Alors, pour $\epsilon$ 
assez petit, l'\'equation aux $q$-diff\'erences avec condition initiale
$$
\begin{cases}
\sigma_{q}X = A X \\
X(0) = X_{0}
\end{cases}
$$
admet une unique solution formelle $X$, qui est convergente et se prolonge
analytiquement \`a $\Omega_{q,0}$. Lorsque $\epsilon$ tend vers $0$, $X$
tend vers $\tilde{X}$.}\\

\textbf{\underline{Preuve}}\\

La remarque qui pr\'ec\`ede imm\'ediatement l'\'enonc\'e du th\'eor\`eme
justifie le sens de l'expression ``$X$ tend vers $\tilde{X}$''. Les calculs de
l'\'etude locale ont montr\'e que le syst\`eme aux $q$-diff\'erences avec
condition initiale ci-dessus est non-r\'esonnant pour $\epsilon$ assez petit,
d'o\`u les affirmations concernant $X$ (par invocation du chapitre 1).
L'hypoth\`ese (B) est celle qui a servi \`a justifier l'\'etude locale 
elle-m\^eme (pour prouver la convergence domin\'ee). L'hypoth\`ese (C) est
celle qui nous a servi \`a justifier l'\'etude globale au 3.2.2.
\hfill $\Box$\\

Nous nommerons la conjonction des hypoth\`eses (B) et (C) ``convergence
uniforme sur suffisamment de compacts de ${\Omega}_{0}$''. Elle est
notamment impliqu\'ee par la convergence uniforme sur tout compact.


\subsection{Application}

Le th\'eor\`eme sera appliqu\'e de fa\c{c}on quelque peu contourn\'ee \`a une
\'etape du d\'evissage des \'equations $\sigma_{q}X = AX$ et 
$\delta \tilde{X} = \tilde{B} \tilde{X}$ au 3.4: celle qui concerne la partie
$H$ dans l'\'ecriture $X = Q_{0} H L C$ des solutions canoniques obtenue
au 1.2, ou la partie $F$ dans l'\'ecriture $X = F e_{q,A}$ obtenue au 1.3.
On donne ici une br\`eve description de ces applications, qui sera reprise
dans son contexte au 3.4.

\subsubsection{Solutions canoniques premi\`ere mani\`ere}

Dans le d\'evissage premi\`ere mani\`ere, on consid\`ere l'\'equation aux
$q$-diff\'erences \emph{matricielle}:
$$
\begin{cases}
\sigma_{q} Y = K Y K_{0}^{-1} \\
Y(0) = I_{n}
\end{cases}
$$
Ici, $K \in GL_{n}(\mathbf{k})$ et $K_{0} = K(0) \in GL_{n}(\mathbf{C})$.\\

De fa\c{c}on analogue, on a affaire \`a l'\'equation diff\'erentielle
(matricielle aussi):
$$
\begin{cases}
\delta \tilde{Y} = \tilde{J} \tilde{Y} - \tilde{Y} \tilde{J}_{0} \\
\tilde{Y}(0) = I_{n}
\end{cases}
$$
Ici, $\tilde{J} \in M_{n}(\mathbf{k})$ et 
$\tilde{J}_{0} = \tilde{J}(0) \in M_{n}(\mathbf{C})$.\\

On suppose cette derni\`ere \'equation non r\'esonnante, au sens o\`u deux
valeurs propres de $\tilde{J}_{0}$ ne peuvent diff\'erer d'un entier non nul.
Cela revient \`a dire qu'aucune valeur propre de l'endomorphisme 
$M \mapsto \tilde{J}(0) \tilde{M} - \tilde{M} \tilde{J}_{0}$ de 
$M_{n}(\mathbf{C})$ n'est un entier non nul (cf 1.3.1.2). On suppose aussi que
$K = I_{n} + \eta J$ o\`u $J \to \tilde{J}$ avec ``convergence uniforme
sur suffisamment de compacts'' de $\tilde{\Omega}_{0}$,
ce dernier d\'efini par les p\^oles de $\tilde{J}$; et que les p\^oles de
$K$ tendent vers ceux de $\tilde{J}$.\\

Ces \'equations matricielles peuvent \^etre vues comme des \'equations
vectorielles dans $\mathbf{C}^{n^{2}}$ auxquelles on applique ce qui 
pr\'ec\`ede: on instancie $X \mapsto AX$ par $Y \mapsto K Y K_{0}^{-1}$ et
$\tilde{X} \mapsto \tilde{B} \tilde{X}$ par 
$\tilde{Y} = \tilde{J} \tilde{Y} - \tilde{Y} \tilde{J}_{0}$. Alors $BX$ est
instanci\'e par
$$
\frac{K Y K_{0}^{-1} - Y}{q - 1} = (J Y - Y J_{0}) K_{0}^{-1}
$$
On a en effet \'ecrit $K = I_{n} + \eta J$ et $K_{0} = I_{n} + \eta J_{0}$.
Comme $J \to \tilde{J}$ et $K_{0} \to I_{n}$ (uniform\'ement), on peut prouver
que les hypoth\`eses du th\'eor\`eme sont satisfaites et l'on conclut
que $Y \to \tilde{Y}$ sur $\tilde{\Omega}_{0}$.

\subsubsection{Solutions canoniques deuxi\`eme mani\`ere}

Ici, l'\'etape correspondante du d\'evissage consiste \`a se ramener \`a une
\'equation \`a coefficients constants par une transformation de jauge 
holomorphe et tangente \`a l'identit\'e.\\

Pour l'\'equation $\sigma_{q}X = AX$, la transformation de jauge est 
sp\'ecifi\'ee par les conditions:
$$
\begin{cases}
\sigma_{q}F = A F A(0)^{-1} \\ 
F(0) = I_{n}
\end{cases}
$$

Pour l'\'equation
$\delta \tilde{X} = \tilde{B} \tilde{X}$, la transformation de jauge est 
sp\'ecifi\'ee par les conditions:
$$
\begin{cases}
\delta \tilde{F} = \tilde{B} \tilde{F} - \tilde{F} \tilde{B}(0) \\
\tilde{F}(0) = I_{n}
\end{cases}
$$

Ces \'equations \emph{matricielles} peuvent \^etre consid\'er\'ees comme des
\'equations \emph{vectorielles} dans $M_{n}(\mathbf{C})$ et on peut leur
appliquer les r\'esultats pr\'ec\'edents: les hypoth\`eses sont manifestement 
v\'erifi\'ees et l'on conclut que $F \to \tilde{F}$ sur $\tilde{\Omega}_{0}$.



\section{Convergence de la partie log-car}

Nous traiterons successivement (et quasi-ind\'ependamment) les deux formes
canoniques obtenues au chapitre 1. La premi\`ere ($X = Q_{0} H L C$) 
n\'ecessitera une renormalisation, tandis que la seconde ($X = F e_{q,A(0)}$) 
non: c'est l\`a un de ses nombreux charmes.


\subsection{Renormalisation et convergence des parties log-car $L C$}

La solution canonique obtenue en 1.2 a la forme $X = Q_{0} H N$, o\`u:

\begin{itemize}

\item{$K_{0} = Q_{0}^{-1} A_{0} Q_{0}$ est une matrice de Jordan, form\'ee de
blocs $c \xi_{1,m}$.}

\item{$N = L C$, partie log-car, est form\'ee de blocs $e_{q,c} L_{m}$ et 
v\'erifie
l'\'equation fonctionnelle matricielle $\sigma_{q} N = K_{0} N$.}

\item{$H$ est l'unique solution du syst\`eme aux $q$-diff\'erences matriciel 
avec conditions initiales:
$$
\begin{cases}
\sigma_{q} H = K H K_{0}^{-1} \\
H(0) = I_{n}
\end{cases}
$$
dans lequel on a pos\'e $K = Q_{0}^{-1} A Q_{0}$.}

\end{itemize}

\subsubsection{Renormalisation des solutions canoniques}

La th\'eorie des fonctions de 3.1 a mis en lumi\`ere le fait que c'est
$(q-1) l_{q}$, et non $l_{q}$ elle-m\^eme, qui se comporte agr\'eablement 
lorsque $q \to 1$. Ainsi nos matrices $L_{m}$ ne convergeront pas. \\

Au lieu des blocs $e_{q,c} L_{m}$, on pr\'ef\`erera donc faire apparaitre des
blocs $e_{q,c} L'_{c,m}$ d\'efinis comme suit:
$$
L'_{c,m}=
\begin{pmatrix}
l_{q}^{(0)} & (\eta/c)^{1} l_{q}^{(1)} & ... & (\eta/c)^{m-1} l_{q}^{(m-1)}\\
0 & l_{q}^{(0)} & ... & (\eta/c)^{m-2} l_{q}^{(m-2)}\\
... & ... & ... & ...\\
0 & 0 & ... & l_{q}^{(0)}\\
\end{pmatrix}
= (\xi'_{1,\eta/c,m})^{l_{q}} 
$$
$$
\xi'_{\lambda,\mu,m} =
\begin{pmatrix}
\lambda & \mu & 0 & ... & 0 & 0\\
0 & \lambda & \mu & ... & 0 & 0\\
... & ... & ... & ... & .. & ...\\
... & ... & ... & ... & .. & ...\\
0 & 0 & 0 & ... & \lambda & \mu\\
0 & 0 & 0 & ... & 0 & \lambda\\
\end{pmatrix}
= \mu \xi_{\lambda/\mu,m} 
$$
Comme d'habitude, on a pos\'e $\eta = q-1$. Ces blocs satisfont l'\'equation
fonctionnelle:
$$
\sigma_{q} (e_{q,c} L'_{c,m}) = \xi'_{c,\eta,m} (e_{q,c} L'_{c,m})
$$

On passe des anciens blocs aux nouveaux par une conjugaison. Notant
$E_{m}(x) = \Diag(1,x,...,x^{m-1})$, on a les relations:
$$
\xi'_{\lambda,\mu,m} = \
(E_{m}(\mu/\lambda))^{-1} (\lambda \xi_{1,m}) E_{m}(\mu/\lambda)
$$
$$
L'_{c,m} = (E_{m}(\eta/c))^{-1} (L_{m}) E_{m}(\eta/c)
$$

La r\'eduite de Jordan de $A_{0}$ sous forme de diagonale de blocs 
renormalis\'es $\xi'$ sera 
alors $K'_{0} = (Q'_{0})^{-1} A_{0} Q'_{0}$, o\`u $Q'_{0} = Q_{0} E_{0}$ et
o\`u $K'_{0} = E_{0}^{-1} K_{0} E_{0}$; la ``matrice de renormalisation''
$E_{0}$ est form\'ee des blocs $E_{m}(\eta/c)$ correspondant aux exposants
$c$. La partie log-car de la solution renormalis\'ee est
$N' = E_{0}^{-1} N E_{0}$, form\'ee des blocs $e_{q,c}$ et v\'erifiant
l'\'equation fonctionnelle $\sigma_{q} N' = K'_{0} N'$.\\

La solution renormalis\'ee est alors $X' = Q'_{0} H' N'$, o\`u $H'$ est
sp\'ecifi\'e par le syst\`eme avec condition initiale:
$$
\begin{cases}
\sigma_{q} H' = K' H' (K'_{0})^{-1} \\
H'(0) = I_{n}
\end{cases}
$$
dans lequel on a pos\'e $K' = (Q'_{0})^{-1} A Q'_{0}$.\\

On voit facilement qu'en fait, $H' = E_{0}^{-1} H E_{0}$, puisque cette
derni\`ere matrice satisfait aux conditions qui sp\'ecifient $H'$ de mani\`ere
univoque. La solution
canonique renormalis\'ee est donc simplement $X' = X E_{0}$.\\

Autrement dit, en termes de la construction vectorielle de 1.2, la solution 
renormalis\'ee s'obtient en multipliant chaque
vecteur d'exposant $c$ et de ``poids logarithmique'' $k$ (c'est \`a dire, dont
l'expression fait intervenir $l_{q}^{(0)}$,...,$l_{q}^{(k)}$) par le facteur
$(\eta/c)^{k}$. Les logarithmes ``sous-dominants'' sont donc vou\'es \`a
disparaitre \`a la limite.

\subsubsection{Convergence de la partie log-car renormalis\'ee}

Supposons maintenant que la matrice $A$ de l'\'equation aux $q$-diff\'erences
varie ``avec structure de Jordan en $0$ continue'', c'est \`a dire: les blocs 
de la r\'eduction de $A(0)$ sont de taille constante, les exposants de la forme
$c_{i} = 1 + \eta \gamma_{i}$, o\`u $\gamma_{i} \to \tilde{\gamma}_{i}$. Ainsi, on peut
\'ecrire $K = I_{n} + \eta J$, avec $J \to \tilde{J}$.\\

Dans ces conditions, on peut d\'eduire des r\'esultats de 3.1.3 et 3.1.4 que:
$$
\begin{cases}
e_{c_{i}}(z) \to (-z)^{\tilde{\gamma}_{i}} \\
(\eta/c_{i})^{k} l_{q}^{(k)}(z) \to \frac{(\log(-z))^{k}}{k!}
\end{cases}
$$
Il s'agit, bien sur, des d\'eterminations principales sur $\Omega$ d\'ej\`a
mentionn\'ees. La partie log-car renormalis\'ee tend donc vers la matrice
$\tilde{N}$ form\'ee des blocs:
$$
(-z)^{\tilde{\gamma}_{i}}
\begin{pmatrix}
1 & \log(-z) & ... & \frac{(log(-z))^{m_{i} - 1}}{(m_{i} - 1)!} \\
0 & 1 & ... & \frac{(log(-z))^{m_{i} - 2}}{(m_{i} - 2)!} \\
... & ... & ... & ... \\
0 & 0 & ... & 1
\end{pmatrix}
= (-z)^{\tilde{\gamma}_{i}} e^{\xi_{0,m_{i}} \log(-z)} =
= (-z)^{\xi_{\tilde{\gamma}_{i},m_{i}}}
$$
qui n'est autre que la solution canonique de l'\'equation diff\'erentielle
$\delta \tilde{L} = \tilde{J} \tilde{L}$ produite par la m\'ethode de Frobenius
(\emph{voir} \cite{Ince}, \cite{Wasow}); \`a ceci pr\`es que nous utilisons 
$(-z)^{\mu}$ et $\log(-z)$ au lieu de $z^{\mu}$ et $\log(z)$.

            
\subsection{Parties log-car de la forme $e_{q,A}$}

Dans le cas non-r\'esonnant, on a introduit au 1.3 une autre version des
solutions fondamentales canoniques, via l'\'equivalence avec une \'equation
\`a coefficients constants: $A$ est \'equivalente \`a $A(0)$ par une 
transformation de jauge $F$ m\'eromorphe sur $\mathbf{C}$, holomorphe et
m\^eme tangente \`a l'identit\'e en $0$. Il y a une solution canonique
sp\'eciale $e_{q,A(0)}$ \`a l'\'equation \`a coefficients constants, et donc
une solution canonique $F e_{q,A(0)}$ \`a l'\'equation de d\'epart. La partie
log-car de celle-ci est donc $e_{A(0)}$, et on va ici \'etudier son 
comportement lorsque $q \to 1$.\\ 

Pour simplifier les notations, on se placera directement dans le cas d'une
famille (param\'etr\'ee par $q$) d'\'equations \`a coefficients constants,
confluant en une \'equation diff\'erentielle \`a coefficients constants.
$A$ est donc une matrice de $GL_{n}(\mathbf{C})$ mais d\'ependant de $q$.

\subsubsection{Remarque pr\'eliminaire}

Consid\'erons les matrices suivantes:
$$ 
\tilde{B} =
\begin{pmatrix}
0 & 0\\
0 & 0\\
\end{pmatrix}
\qquad
B_{\epsilon} = \epsilon
\begin{pmatrix}
\cos^{2} \frac{1}{\epsilon} & 
- \sin \frac{1}{\epsilon} \cos \frac{1}{\epsilon}\\
- \sin \frac{1}{\epsilon} \cos \frac{1}{\epsilon} & 
\sin^{2} \frac{1}{\epsilon}\\
\end{pmatrix}
$$
On voit donc que $\underset{\epsilon \to 0^{+}}{\lim} B_{\epsilon} = \tilde{B}$. Mais
$B_{\epsilon}/\epsilon$ est la matrice de la projection orthogonale sur le
vecteur de coordonn\'ees $(\cos 1/\epsilon, -\sin 1/\epsilon)$. Toute
trigonalisation de $B_{\epsilon}$ sera donc obtenue dans une base tournant 
infiniment vite lorsque $\epsilon \to 0^{+}$, et l'on ne peut d\'eployer
\emph{aucune} trigonalisation de $\tilde{B}$ en une trigonalisation de 
$B_{\epsilon}$. Ceci justifie l'introduction d'une hypoth\`ese particuli\`ere
dans la d\'efinition qui suit.\\

\textbf{\underline{D\'efinition}}\\

\emph{On dira que $B$ tend vers $\tilde{B}$ avec sa triangularisation s'il 
existe une triangularisation $S^{-1}BS = T$ (d\'ependant donc de $q$) et une
triangularisation $\tilde{S}^{-1}\tilde{B}\tilde{S} = \tilde{T}$ telles que
$S \to \tilde{S}$ et donc $T \to \tilde{T}$.}\\

C'est par exemple le cas si $\tilde{B}$ est diagonalisable.

\subsubsection{Convergence de $e_{q,A}$}

On suppose ici que, pour chaque valeur de $q$, $B \in M_{n}(\mathbf{C})$. On
note $A = I_{n} + (q-1)B$ et l'on se propose de pr\'eciser le comportement de
$e_{q,A}$ lorsque $\epsilon \to 0^{+}$ et que $B \to \tilde{B}$ avec sa
triangularisation. On reprend les notations pr\'ec\'edentes.\\

Soit $\overline{D} + \overline{N}$ la d\'ecomposition de Dunford additive de 
$B$. La d\'ecomposition de Dunford multiplicative $DU$ de $A$ s'obtient en
prenant $D = I_{n} + (q-1) \overline{D}$ et 
$U = I_{n} + (q-1) D^{-1} \overline{N}$. Alors la matrice diagonale
$\Delta = S D S^{-1}$ et la matrice diagonale 
$\overline{\Delta} = S \overline{D} S^{-1}$ (qui est la partie diagonale de 
$T$) sont li\'ees par la relation: $\Delta = I_{n} + (q-1)\overline{\Delta}$.\\

Soient $\tilde{D} + \tilde{N}$ la d\'ecomposition de Dunford additive de
$\tilde{B}$ et $\tilde{\Delta} = \tilde{S} \tilde{D} \tilde{S}^{-1}$ la partie 
diagonale de $\tilde{T}$. Donc $(\Delta - I_{n})/(q-1)$ tend vers 
$\tilde{\Delta}$. Il est
alors imm\'ediat (d'apr\`es les r\'esultats de 3.1.3) que $e_{q,\Delta}$ tend
vers $(-z)^{\tilde{\Delta}}$ et donc $e_{q,D}$ tend vers $(-z)^{\tilde{D}}$.\\

Comme $D^{-1} \overline{N}$ tend vers $\tilde{N}$, il d\'ecoule pareillement
des r\'esultats de 3.1.4 que la limite de
$U^{l_{q}} = (I_{n} + (q-1)D^{-1} \overline{N})^{l_{q}}$ est
$(-z)^{\tilde{N}}$: cela revient en effet \`a dire que 
$(q-1)^{k}l_{q}^{(k)} \to (\log(-z))^{k}/k!$.\\

On voit donc que $e_{q,A}$ tend vers $(-z)^{\tilde{B}}$.



\section{D\'evissages en parall\`ele}

On se donne maintenant

\begin{itemize}

\item{Une \'equation diff\'erentielle lin\'eaire \`a coefficients rationnels 
fuchsienne non r\'esonnante: $\delta \tilde{X} = \tilde{B} \tilde{X}$. On nomme
$\tilde{z}_{1},...,\tilde{z}_{r}$ les p\^oles de $\tilde{B}$: ce sont des
\'el\'ements de $\mathbf{C}^{*}$. On note
$$
\tilde{J}_{0} = (\tilde{Q}_{0})^{-1} \tilde{B}_{0} \tilde{Q}_{0} =
\begin{pmatrix}
\xi_{\tilde{\gamma}_{1},m_{1}} & 0 & ... & 0 \\
0 & \xi_{\tilde{\gamma}_{2},m_{2}} & ... & 0 \\
... & ... & ... & ... \\
0 & 0 & ... & \xi_{\tilde{\gamma}_{r},m_{r}}
\end{pmatrix}
$$
une r\'eduction de Jordan de $\tilde{B}_{0} = \tilde{B}(0)$.}

\item{Une \'equation aux $q$-diff\'erences lin\'eaire \`a coefficients
rationnels fuchsienne en $0$: $\sigma_{q} X = A X$, o\`u 
$A \in GL_{n}(\mathbf{k})$ d\'epend de $q$. On suppose que l'on peut \'ecrire
la matrice $A$ sous la forme $I_{n} + \eta B \;\;(\eta = q - 1)$, avec 
$B \underset{\epsilon \to 0^{+}}{\to} \tilde{B}$.}

\end{itemize}

On supposera que la convergence $B \to \tilde{B}$ est astreinte \`a satisfaire
les conditions suivantes:

\begin{itemize}

\item{(A) Les p\^oles de $A$ (et donc de $B$) tendent vers ceux de 
$\tilde{B}$.}

\item{(B) $B \to \tilde{B}$ et il y a ``convergence uniforme sur suffisamment
de compacts'' de $\tilde{\Omega}_{0}$.}

\item{(C) La variation de $B_{0}$ est \`a structure de Jordan constante: 
$\tilde{Q}_{0}$
se prolonge en $Q_{0}$, fonction continue de $q$ dont elle est la limite
lorsque $\epsilon \to 0^{+}$, de sorte que l'on ait une r\'eduction de Jordan
en famille:
$$
J_{0} = (Q_{0})^{-1} B_{0} Q_{0} =
\begin{pmatrix}
\xi_{\gamma_{1},m_{1}} & 0 & ... & 0 \\
0 & \xi_{\gamma_{2},m_{2}} & ... & 0 \\
... & ... & ... & ... \\
0 & 0 & ... & \xi_{\gamma_{r},m_{r}}
\end{pmatrix}
$$
et $\gamma_{i} \to \tilde{\gamma}_{i} \;\;(i = 1,...,r)$.}

\end{itemize}

Toutes ces conditions sont, par exemple, v\'erifi\'ees si l'on prend
$B = \tilde{B}$!\\

Dans ces conditions, la variation de $A_{0}$ est aussi \`a structure de Jordan
constante:
$$
K_{0} = (Q_{0})^{-1} A_{0} Q_{0} =
\begin{pmatrix}
\xi'_{c_{1},\eta,m_{1}} & 0 & ... & 0 \\
0 & \xi'_{c_{2},\eta,m_{2}} & ... & 0 \\
... & ... & ... & ... \\
0 & 0 & ... & \xi'_{c_{r},\eta,m_{r}}
\end{pmatrix}
$$
avec $c_{i} = 1 + \eta \gamma_{i}  \;\;(i = 1,...,r)$.\\

La r\'esolution canonique de l'\'equation diff\'erentielle par la m\'ethode de 
Frobenius comporte les \'etapes suivantes 
(\emph{voir} \cite{Ince},\cite{Wasow}):

\begin{enumerate}

\item{On choisit une matrice ``log-exposants'' $\tilde{N}$ telle que
$\delta \tilde{N} = \tilde{J}_{0} \tilde{N}$. Le choix canonique est celui
d'une matrice diagonale par blocs, dont les blocs sont d\'etermin\'es par
la structure de Jordan et ont la forme:
$$
(-z)^{\tilde{\gamma}_{i}}
\begin{pmatrix}
1 & \log(-z) & ... & \frac{(log(-z))^{m_{i} - 1}}{(m_{i} - 1)!} \\
0 & 1 & ... & \frac{(log(-z))^{m_{i} - 2}}{(m_{i} - 2)!} \\
... & ... & ... & ... \\
0 & 0 & ... & 1
\end{pmatrix}
= (-z)^{\tilde{\gamma}_{i}} e^{\xi_{0,m_{i}} \log(-z)} =
= (-z)^{\xi_{\tilde{\gamma}_{i},m_{i}}}
$$
A ceci pr\`es, bien sur, que le choix classique comporte plut\^ot $\log z$ et
$z^{\tilde{\gamma}}$, que nous avons remplac\'es par les d\'eterminations
principales de $\log(-z)$ et de $(-z)^{\tilde{\gamma}}$ sur $\Omega$.}

\item{On cherche alors la solution $\tilde{X}$ sous la forme
$\tilde{Q}_{0} \tilde{Y} \tilde{N}$, $\tilde{Y}$ \'etant sp\'ecifi\'e par
l'\'equation diff\'erentielle matricielle avec condition initiale:
$$
\begin{cases}
\delta \tilde{Y} = \tilde{J} \tilde{Y} - \tilde{Y} \tilde{J}_{0} \\
\tilde{Y}(0) = I_{n}
\end{cases}
$$
On a pos\'e $\tilde{J} = (\tilde{Q}_{0})^{-1} \tilde{B} \tilde{Q}_{0}$.
Ceci d\'etermine alors un unique $\tilde{Y}$ sur $\tilde{\Omega}_{0}$, grace 
\`a l'hypoth\`ese de non-r\'esonnance.}

\end{enumerate}


\subsection{Convergence des solutions du type $Q_{0} H L C$}

La r\'esolution canonique de l'\'equation aux $q$-diff\'erences obtenue en 1.2,
telle qu'elle a \'et\'e renormalis\'ee en 3.3 comporte les \'etapes suivantes:

\begin{enumerate}

\item{Une matrice log-car $N = L C$ v\'erifie 
$\sigma_{q} N = K_{0} N$. Le choix
canonique (une fois choisie une r\'eduction de Jordan) est celui de la matrice
diagonale par blocs dont les blocs sont d\'etermin\'es par la structure de
Jordan et ont la forme:
$$
e_{q,c_{i}} 
\begin{pmatrix}
l_{q}^{(0)} & (\eta/c)^{1} l_{q}^{(1)} & ... & (\eta/c)^{m-1} l_{q}^{(m-1)}\\
0 & l_{q}^{(0)} & ... & (\eta/c)^{m-2} l_{q}^{(m-2)}\\
... & ... & ... & ...\\
0 & 0 & ... & l_{q}^{(0)}\\
\end{pmatrix}
$$
On a vu que, dans ces conditions, $N \to \tilde{N}$ sur $\Omega$.}

\item{On cherche alors la solution $X$ sous la forme $Q_{0} Y N$, $Y$ \'etant
sp\'ecifi\'e par le syst\`eme aux $q$-diff\'erences matriciel avec condition
initiale:
$$
\begin{cases}
\sigma_{q} Y = K Y K_{0}^{-1} \\
Y(0) = I_{n}
\end{cases}
$$
On a pos\'e $K = Q_{0}^{-1} A Q_{0}$. Ceci d\'etermine alors, pour $\epsilon$
assez petit, un unique $Y$ sur $\Omega_{q,0}$, grace \`a la non-r\'esonnance. 
On a vu au 3.2 que $Y \to \tilde{Y}$ sur $\tilde{\Omega}_{0}$.}

\end{enumerate}

On r\'eunit ces conclusions en un\\

\textbf{\underline{Th\'eor\`eme}}\\

\emph{La solution canonique renormalis\'ee $X$ tend vers la solution canonique
$\tilde{X}$ sur $\tilde{\Omega}'_{0} = \Omega \cap \tilde{\Omega}_{0}$.}
\hfill $\Box$


\bigskip
\hrule
\bigskip

\unitlength = 1cm

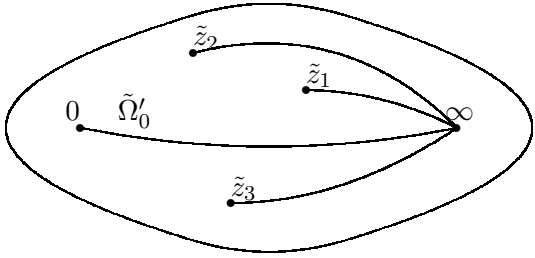
\underline{Figure 3}\\

\bigskip

\begin{picture} (15,3)


\qbezier (4,0)(9,1.5)(4,3)     

\qbezier (2,0)(-3,1.5)(2,3)      

\qbezier (2,3)(3,3.3)(4,3)   

\qbezier (2,0)(3,-0.3)(4,0)  


\put(0.5,1.5){\circle*{0.1}}
\put(0.3,1.6){$0$}

\put(5.5,1.5){\circle*{0.1}}
\put(5.35,1.6){$\infty$}

\put(2,2.5){\circle*{0.1}}
\put(2,2.6){$\tilde{z}_{2}$}

\put(3.5,2){\circle*{0.1}}
\put(3.5,2.1){$\tilde{z}_{1}$}

\put(2.5,0.5){\circle*{0.1}}
\put(2.5,0.6){$\tilde{z}_{3}$}


\qbezier (2,2.5) (4,3) (5.5,1.5)

\qbezier (3.5,2) (4.5,2) (5.5,1.5)

\qbezier (2.5,0.5) (4,0.5) (5.5,1.5)

\qbezier (0.5,1.5) (3,1) (5.5,1.5)


\put(1,1.6){$\tilde{\Omega}'_{0}$}

\end{picture}

\bigskip
\bigskip
\bigskip
\hrule
\bigskip


\subsection{Convergence des solutions du type $F e_{q,A(0)}$}

Le raisonnement est en tout point analogue \`a celui qui pr\'ec\`ede. 
Pour l'\'equation $\sigma_{q}X = AX$, la transformation de jauge $F$ est 
sp\'ecifi\'ee par les conditions:
$$
\begin{cases}
\sigma_{q}F = A F A(0)^{-1} \\
F(0) = I_{n}
\end{cases}
$$

Pour l'\'equation $\delta \tilde{X} = \tilde{B} \tilde{X}$, la transformation 
de jauge correspondante est sp\'ecifi\'ee par les conditions:
$$
\begin{cases}
\delta \tilde{F} = \tilde{B} \tilde{F} - \tilde{F} \tilde{B}(0) \\
\tilde{F}(0) = I_{n}
\end{cases}
$$

Ces \'equations \emph{matricielles} peuvent \^etre consid\'er\'ees comme des
\'equations \emph{vectorielles} dans $M_{n}(\mathbf{C})$ et on peut leur
appliquer les r\'esultats de 3.2. Les hypoth\`eses sont manifestement 
v\'erifi\'ees et l'on conclut que $F \to \tilde{F}$ sur $\tilde{\Omega}_{0}$.
On a vu au 3.3 que (en r\'eadaptant les notations) $e_{q,A(0)}$ tend vers 
$(-z)^{\tilde{B(0)}}$, de sorte que:\\

\textbf{\underline{Th\'eor\`eme}}\\

\emph{La solution canonique $F e_{q,A(0)}$ tend vers la solution canonique
$X = \tilde{F} (-z)^{\tilde{B(0)}}$ sur 
$\tilde{\Omega}'_{0} = \Omega \cap \tilde{\Omega}_{0}$.}
\hfill $\Box$\\

En effet, ce d\'evissage l\'eg\`erement diff\'erent produit bien la solution
canonique de Frobenius du syst\`eme diff\'erentiel.\\

\textbf{\underline{Corollaire}}\\

\emph{Ces r\'esultats s'appliquent en particulier \`a la d\'eformation
lin\'eaire $A = I_{n} + (q-1) \tilde{B}$.}



\chapter{Matrice de connexion et monodromie}



On part maintenant d'une \'equation diff\'erentielle fuchsienne
$$
\delta \tilde{X} = \tilde{B} \tilde{X}
$$
o\`u $\tilde{B} \in M_{n}(\mathbf{k})$. Les p\^oles de $\tilde{B}$ sont not\'es
$\tilde{z}_{1},...,\tilde{z}_{r}$: ce sont des \'el\'ements de
$\mathbf{C}^{*}$ puisque l'\'equation est fuchsienne en $0$ et en $\infty$.
On suppose en outre l'\'equation non-r\'esonnante en $0$ et en $\infty$, c'est
\`a dire que deux valeurs propres distinctes de $\tilde{B}(0)$ (resp. de
$\tilde{B}(\infty)$) ne peuvent diff\'erer d'un entier.\\

On se donne alors une d\'eformation $B$ de $\tilde{B}$, autrement dit une 
matrice $B \in M_{n}(\mathbf{k})$ dont les coefficients d\'ependent de $q$
et qui tend vers $\tilde{B}$ lorsque $q \to 1$.
La variation de $q$ ob\'eit aux m\^emes r\`egles que pr\'ec\'edemment: elle
a lieu le long d'une spirale logarithmique $q_{0}^{\mathbf{R_{+}^{*}}}$ et on 
la param\`etre par un $\epsilon > 0$ tel que $q = q_{0}^{\epsilon}$. On
d\'enotera de telles familles par l'indice $\epsilon$, par exemple 
$B_{\epsilon}$ pour celle que l'on vient d'introduire.\\

On suppose de plus de la famille des $B_{\epsilon}$ qu'elle satisfait aux 
contraintes suivantes:

\begin{itemize}

\item{(A) Les p\^oles de $B_{\epsilon}$ tendent vers ceux de $\tilde{B}$; 
ceci signifie
qu'en tout point de $\mathbf{S}$ qui n'est pas un p\^ole de $\tilde{B}$, la
matrice $B_{\epsilon}$ est d\'efinie (n'a pas de p\^ole) d\`es que $\epsilon$ 
est assez
petit.}

\item{(B) 
$\underset{\epsilon \to 0^{+}}{\lim} B_{\epsilon} = \tilde{B}$, la convergence
\'etant de plus uniforme sur: (i) un disque ferm\'e non trivial de centre $0$
ne contenant aucun des $\tilde{z}_{i}$; (ii) un disque ferm\'e non trivial de 
centre $\infty$ ne contenant aucun des $\tilde{z}_{i}$; (iii) tout segment de
spirale logarithmique inclus dans 
$\mathbf{C}^{*} - \{ \tilde{z}_{1},...,\tilde{z}_{r} \}$. Par exemple, la 
convergence uniforme sur tout compact de 
$\mathbf{C}^{*} - \{ \tilde{z}_{1},...,\tilde{z}_{r} \}$ entraine cette
hypoth\`ese.}

\item{(C) 
$B_{\epsilon}(0)$ converge vers $\tilde{B}(0)$ et $B_{\epsilon}(\infty)$ 
converge vers
$\tilde{B}(\infty)$ \`a structure de Jordan constante, avec matrice de passage
continue (voir 3.4).}

\end{itemize}

L'une ou l'autre condition de convergence de (B) a bien un sens, car la 
condition (A) entraine que, sur tout compact de
$\mathbf{C}^{*} - \{ \tilde{z}_{1},...,\tilde{z}_{r} \}$,
$B_{\epsilon}$ n'a pas de p\^ole pour $\epsilon$ assez petit.\\

\emph{Toutes ces hypoth\`eses sont en particulier satisfaites par la
d\'eformation ``lin\'eaire'', d\'efinie par $B_{\epsilon} = \tilde{B}$
quelque soit $\epsilon > 0$.}


\section{Les hypoth\`eses sont sym\'etriques en $0$, $\infty$}

On va d'abord \'etudier l'effet du changement de variable $w = 1/z$ sur 
l'\'equation diff\'erentielle:
$$
\delta \tilde{X} = \tilde{B} \tilde{X}
$$

Posant $\tilde{Y}(w) = \tilde{X}(w) = \tilde{X}(z)$, on sait que:
$$
\delta_{w} \tilde{Y} = \tilde{D} \tilde{Y} \qquad 
\delta_{w} := w \frac{d}{dw}
$$
o\`u l'on note $\tilde{D}(w) = - \tilde{B}(1/w)$. Les p\^oles de $\tilde{D}$
sont donc les $\tilde{w}_{i} = 1/\tilde{z}_{i}$ et
$$
\begin{cases}
\tilde{D}(0) = \tilde{B}(\infty) \\
\tilde{D}(\infty) = \tilde{B}(0)
\end{cases}
$$
Il y a bien sym\'etrie des hypoth\`eses.\\

On consid\`ere ensuite l'\'equation aux $q$-diff\'erences qui doit d\'eformer
notre \'equation diff\'erentielle:
$$
\sigma_{q} X = A_{\epsilon} X
$$
o\`u $A_{\epsilon} = I_{n} + \eta B_{\epsilon} \; (\eta = q-1)$. 
Pour $q$ voisin de $1$, elle est
fuchsienne puisque $\det A_{\epsilon}(0)$ et $\det A_{\epsilon}(\infty)$ 
ont pour limite $1$. On a
vu au chapitre 2 l'effet du changement de variable $w = 1/z$: Posant
$Y(w) = X(1/w) = X(z)$, on a:
$$
\sigma_{q} Y = C_{\epsilon} Y
$$
o\`u l'on note $C_{\epsilon}(w) = (A_{\epsilon}(1/qw))^{-1}$. 
Cette \'equation est \'evidemment 
encore fuchsienne, puisque 
$$
\begin{cases}
C_{\epsilon}(0) = A_{\epsilon}(\infty) \\
A_{\epsilon}(0) = C_{\epsilon}(\infty)
\end{cases}
$$
De plus, $C_{\epsilon}(w) = I_{n} + \eta D_{\epsilon}(w)$, 
o\`u $D_{\epsilon}(w) = B_{\epsilon}(1/qw) (A_{\epsilon}(1/qw))^{-1}$. On
v\'erifie alors ais\'ement les hypoth\`eses sym\'etriques de celles que nous
avons faites:

\begin{itemize}

\item{(A') Les p\^oles de $D_{\epsilon}$ sont les $w$ tels que $1/qw$ est 
p\^ole de $B_{\epsilon}$
ou bien z\'ero de $\det A_{\epsilon}$. Mais il est clair que, si
$z \in \mathbf{C}^{*} - \{ \tilde{z}_{1},...,\tilde{z}_{r} \}$, 
$\det A_{\epsilon}(z)$
tend vers $1$. Donc les p\^oles de $D_{\epsilon}$ tendent vers les 
$\tilde{w}_{i}$, c'est 
\`a dire vers ceux de $\tilde{D}$.}

\item{(B') La m\^eme formule montre que $D_{\epsilon}$ 
converge uniform\'ement vers
$\tilde{D}$ sur un disque centr\'e en $0$, sur un disque centr\'e en $\infty$
(\'evitant tous deux les $\tilde{w}_{i}$) et les segments de spirales inclus
dans $\mathbf{C}^{*} - \{ \tilde{w}_{1},...,\tilde{w}_{r} \}$.}

\item{(C') Enfin, des hypoth\`eses sur $B_{\epsilon}(0)$ et 
$B_{\epsilon}(\infty)$ et des formules
$$
\begin{cases}
D_{\epsilon}(0) = \frac{C_{\epsilon}(0) - I_{n}}{q - 1} = 
\frac{(A_{\epsilon}(\infty))^{-1} - I_{n}}{q - 1} \\
D_{\epsilon}(\infty) = \frac{C_{\epsilon}(\infty) - I_{n}}{q - 1} = 
\frac{(A_{\epsilon}(0))^{-1} - I_{n}}{q - 1}
\end{cases}
$$
on d\'eduit que $D_{\epsilon}(0)$ converge vers $\tilde{D}(0)$ et 
$D_{\epsilon}(\infty)$ converge vers
$\tilde{D}(\infty)$ \`a structure de Jordan constante, avec matrice de passage
continue.}

\end{itemize}

Au voisinage de $\infty$, c'est donc 
$C_{\epsilon} = I_{n} + \eta D_{\epsilon}$ que l'on 
consid\`erera comme une d\'eformation de $\tilde{D}$.\\


\section{La matrice de connexion}

Notons $\tilde{X}^{(0)}$ et $\tilde{X}^{(\infty)}$ les solutions canoniques
en $z$ et en $w$ de l'\'equation diff\'erentielle obtenue par la m\'ethode
de Frobenius. Sur un ouvert de $\mathbf{S}$ o\`u toutes deux sont des 
syst\`emes fondamentaux (uniformes, de d\'eterminant non nul), il y a une 
unique matrice $\tilde{P}(z)$ telle que 
$$
\tilde{X}^{(0)}(z) = \tilde{X}^{(\infty)}(z) \tilde{P}(z)
$$

En d\'erivant des deux c\^ot\'es et en utilisant l'\'equation diff\'erentielle,
on voit que $\tilde{P}' = 0$: $\tilde{P}$ est localement constante sur cet
ouvert. Par exemple, tout ouvert simplement connexe ne contenant aucun des
$\tilde{z}_{i}$ conviendrait et $\tilde{P}$ y serait constante.\\

Notons de m\^eme $X_{\epsilon}^{(0)}$ et $X_{\epsilon}^{(\infty)}$ 
les solutions canoniques de
l'\'equation aux $q$-diff\'erences de matrice $A_{\epsilon}$. 
Ce peuvent \^etre celles obtenues par la
m\'ethode du 1.2, renormalis\'ees comme au 3.3; ou bien celles obtenues par
la m\'ethode du 1.3, puisque l'on sait que, pour $q$ voisin de $1$,
l'\'equation est non-r\'esonnante (en $0$ aussi bien qu'en $\infty$ d'apr\`es
la sym\'etrie des hypoth\`eses). Sur un ouvert $\sigma_{q}$-invariant
sur lequel $X_{\epsilon}^{(0)}$ et 
$X_{\epsilon}^{(\infty)}$ sont des solutions fondamentales, 
il y a une unique matrice
$P_{\epsilon}(z)$ telle que 
$$
X_{\epsilon}^{(0)}(z) = X_{\epsilon}^{(\infty)}(z) P_{\epsilon}(z)
$$

C'est, bien sur, un choix de la matrice de connexion.\\
 
En appliquant $\sigma_{q}$ des deux
c\^ot\'es et en utilisant l'\'equation aux $q$-diff\'erences, on voit que
$\sigma_{q} P_{\epsilon} = P_{\epsilon}$: $P_{\epsilon}$ est elliptique, 
``constante'' pour notre th\'eorie.\\


\bigskip \hrule \bigskip 

\unitlength = 1cm

\underline{Figure 4 a}\\

\bigskip

\begin{picture} (15,3)


\qbezier (4,0)(9,1.5)(4,3)     

\qbezier (2,0)(-3,1.5)(2,3)      

\qbezier (2,3)(3,3.3)(4,3)   

\qbezier (2,0)(3,-0.3)(4,0)  


\put(0.5,1.5){\circle*{0.1}}
\put(0.3,1.6){$0$}

\put(5.5,1.5){\circle*{0.1}}
\put(5.35,1.6){$\infty$}

\put(2,2.5){\circle*{0.1}}
\put(2,2.6){$z_{2}$}

\put(3.5,2){\circle*{0.1}}
\put(3.5,2.1){$z_{1}$}

\put(2.5,0.5){\circle*{0.1}}
\put(2.5,0.6){$z_{3}$}


\qbezier [20] (2,2.5) (4,3) (5.5,1.5)

\qbezier [20] (3.5,2) (4.5,2) (5.5,1.5)

\qbezier [20] (2.5,0.5) (4,0.5) (5.5,1.5)

\qbezier [50] (0.5,1.5) (3,1) (5.5,1.5)

\qbezier [50] (0.5,1.5) (3,0.5) (5.5,1.5)

\qbezier [50] (0.5,1.5) (3,1.5) (5.5,1.5)


\put(1.5,-0.8) {Singularit\'es de $X_{\epsilon}^{(0)}$}

\end{picture}

\bigskip \bigskip \bigskip \hrule \bigskip

Pour que $X_{\epsilon}^{(0)}$ y soit d\'efini, un tel ouvert ne doit
rencontrer aucune 
demi-spirale discr\`ete partant des p\^oles de $A_{\epsilon}$; 
pour que $X_{\epsilon}^{(0)}$ y soit
solution fondamentale, il faut en outre qu'il ne contienne aucune 
demi-spirale discr\`ete partant des z\'eros de $\det A_{\epsilon}$.
Il s'agit, bien sur, des demi-spirales engendr\'ees par $q$. 
Enfin, il faut exclure les spirales de p\^oles et de z\'eros des caract\`eres,
et, \'eventuellement, la spirale des p\^oles de $l_{q}$, selon la structure
de Jordan de $A_{\epsilon}(0)$. On a des \'enonc\'es analogues en $\infty$.\\


\bigskip \hrule \bigskip 

\unitlength = 1cm

\underline{Figure 4 b}\\

\bigskip

\begin{picture} (15,3)


\qbezier (4,0)(9,1.5)(4,3)     

\qbezier (2,0)(-3,1.5)(2,3)      

\qbezier (2,3)(3,3.3)(4,3)   

\qbezier (2,0)(3,-0.3)(4,0)  


\put(0.5,1.5){\circle*{0.1}}
\put(0.3,1.6){$0$}

\put(5.5,1.5){\circle*{0.1}}
\put(5.35,1.6){$\infty$}

\put(2,2.5){\circle*{0.1}}
\put(2,2.6){$z_{2}$}

\put(3.5,2){\circle*{0.1}}
\put(3.5,2.1){$z_{1}$}

\put(2.5,0.5){\circle*{0.1}}
\put(2.5,0.6){$z_{3}$}


\qbezier [20] (0.5,1.5) (0.75,2) (02,2.5)

\qbezier [20] (0.5,1.5) (1.5,2) (03.5,2)

\qbezier [20] (0.5,1.5) (1.5,0.5) (02.5,0.5)

\qbezier [50] (0.5,1.5) (03,1) (05.5,1.5)

\qbezier [50] (0.5,1.5) (03,0.5) (05.5,1.5)

\qbezier [50] (0.5,1.5) (03,1.5) (05.5,1.5)


\put(1.5,-0.8) {Singularit\'es de $X_{\epsilon}^{(\infty)}$}

\end{picture}

\bigskip \bigskip \bigskip \hrule \bigskip

On introduit maintenant les ouverts simplement connexes:
$$
\begin{cases}
\tilde{\Omega}_{0} = 
\mathbf{S} - 
\underset{1 \leq i \leq r}{\bigcup} \tilde{z}_{i} q_{0}^{\overline{\mathbf{R}}_{+}} \\
\tilde{\Omega}_{\infty} = 
\mathbf{S} - 
\underset{1 \leq i \leq r}{\bigcup} \tilde{z}_{i} q_{0}^{\overline{\mathbf{R}}_{-}}
\end{cases}
\qquad
\begin{cases}
\tilde{\Omega}'_{0} = \tilde{\Omega}_{0} - q_{0}^{\overline{\mathbf{R}}} \\
\tilde{\Omega}'_{\infty} = \tilde{\Omega}_{\infty} - q_{0}^{\overline{\mathbf{R}}}
\end{cases}
\qquad
\tilde{\Omega}' = \tilde{\Omega}'_{0} \cap \tilde{\Omega}'_{\infty}
$$
On sait alors que $X_{\epsilon}^{(0)}$ tend vers $\tilde{X}^{(0)}$ sur 
$\tilde{\Omega}'_{0}$
et que $X_{\epsilon}^{(\infty)}$ tend vers $\tilde{X}^{(\infty)}$ sur 
$\tilde{\Omega}'_{\infty}$:
cela a \'et\'e prouv\'e au chapitre 3. De plus, quelque soit 
$z \in \tilde{\Omega}'$, 
$P_{\epsilon}$ y est d\'efinie d\`es que $q$ est suffisamment proche de $1$. 


\bigskip \hrule \bigskip

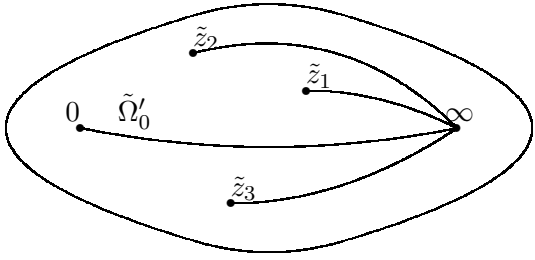
\unitlength = 1cm

\underline{Figure 5 a}\\

\bigskip

\begin{picture} (15,3)


\qbezier (4,0)(9,1.5)(4,3)     

\qbezier (2,0)(-3,1.5)(2,3)      

\qbezier (2,3)(3,3.3)(4,3)   

\qbezier (2,0)(3,-0.3)(4,0)  


\put(0.5,1.5){\circle*{0.1}}
\put(0.3,1.6){$0$}

\put(5.5,1.5){\circle*{0.1}}
\put(5.35,1.6){$\infty$}

\put(2,2.5){\circle*{0.1}}
\put(2,2.6){$\tilde{z}_{2}$}

\put(3.5,2){\circle*{0.1}}
\put(3.5,2.1){$\tilde{z}_{1}$}

\put(2.5,0.5){\circle*{0.1}}
\put(2.5,0.6){$\tilde{z}_{3}$}


\qbezier (2,2.5) (4,3) (5.5,1.5)

\qbezier (3.5,2) (4.5,2) (5.5,1.5)

\qbezier (2.5,0.5) (4,0.5) (5.5,1.5)

\qbezier (0.5,1.5) (3,1) (5.5,1.5)


\put(1,1.6){$\tilde{\Omega}'_{0}$}


\put(1,-0.8) {Domaine de d\'efinition de $\tilde{X}^{(0)}$}

\end{picture}

\bigskip \bigskip \bigskip \hrule \newpage


\bigskip \hrule \bigskip

\unitlength = 1cm

\underline{Figure 5 b}\\

\bigskip

\begin{picture} (15,3)


\qbezier (4,0)(9,1.5)(4,3)     

\qbezier (2,0)(-3,1.5)(2,3)      

\qbezier (2,3)(3,3.3)(4,3)   

\qbezier (2,0)(3,-0.3)(4,0)  


\put(0.5,1.5){\circle*{0.1}}
\put(0.3,1.6){$0$}

\put(5.5,1.5){\circle*{0.1}}
\put(5.35,1.6){$\infty$}

\put(2,2.5){\circle*{0.1}}
\put(2,2.6){$\tilde{z}_{2}$}

\put(3.5,2){\circle*{0.1}}
\put(3.5,2.1){$\tilde{z}_{1}$}

\put(2.5,0.5){\circle*{0.1}}
\put(2.5,0.6){$\tilde{z}_{3}$}


\qbezier (0.5,1.5) (0.75,2) (02,2.5)

\qbezier (0.5,1.5) (1.5,2) (03.5,2)

\qbezier (0.5,1.5) (1.5,0.5) (02.5,0.5)

\qbezier (0.5,1.5) (03,1) (05.5,1.5)


\put(04.5,1.6){$\tilde{\Omega}'_{\infty}$}


\put(1,-0.8) {Domaine de d\'efinition de $\tilde{X}^{(\infty)}$}

\end{picture}

\bigskip \bigskip \bigskip \bigskip \hrule \bigskip

On obtient donc le\\

\textbf{\underline{Th\'eor\`eme}}\\

\emph{$P_{\epsilon}$ converge vers $\tilde{P}$ sur $\tilde{\Omega}'$.}\\

Posons, pour simplifier les notations, $\tilde{z}_{0} = 1$. L'ouvert 
$\tilde{\Omega}' = \tilde{\Omega}'_{0} \cap \tilde{\Omega}'_{\infty} = 
\mathbf{S} - \bigcup_{0 \leq i \leq r} \tilde{z}_{i} q_{0}^{\overline{R}}$
est la sph\`ere de Riemann priv\'ee d'un nombre fini de coupures spirales
compactes reliant $0$ \`a $\infty$. Il a autant de composantes connexes 
$\tilde{U}_{i}$
qu'il y a de spirales distinctes et chacune est hom\'eomorphe \`a un disque
(et m\^eme diff\'eomorphe). La matrice $\tilde{P}$ est constante sur chacune 
de ces ``tranches de melon''.\\


\bigskip \hrule \bigskip

\unitlength = 1cm

\underline{Figure 6}\\

\bigskip

\begin{picture} (15,3)


\qbezier (4,0)(9,1.5)(4,3)     

\qbezier (2,0)(-3,1.5)(2,3)      

\qbezier (2,3)(3,3.3)(4,3)   

\qbezier (2,0)(3,-0.3)(4,0)  


\put(0.5,1.5){\circle*{0.1}}
\put(0.3,1.6){$0$}

\put(5.5,1.5){\circle*{0.1}}
\put(5.35,1.6){$\infty$}

\put(2,2.5){\circle*{0.1}}
\put(2,2.6){$\tilde{z}_{2}$}

\put(3.5,2){\circle*{0.1}}
\put(3.5,2.1){$\tilde{z}_{1}$}

\put(2.5,0.5){\circle*{0.1}}
\put(2.5,0.6){$\tilde{z}_{3}$}


\qbezier (2,2.5) (4,3) (5.5,1.5)

\qbezier (3.5,2) (4.5,2) (5.5,1.5)

\qbezier (2.5,0.5) (4,0.5) (5.5,1.5)

\qbezier (0.5,1.5) (3,1) (5.5,1.5)

\qbezier (0.5,1.5) (0.75,2) (2,2.5)

\qbezier (0.5,1.5) (1.5,2) (3.5,2)

\qbezier (0.5,1.5) (1.5,0.5) (2.5,0.5)


\put(3,1.6){$\tilde{U}_{0}$}

\put(3,2.2){$\tilde{U}_{1}$}

\put(3,2.7){$\tilde{U}_{2}$}

\put(3,0.8){$\tilde{U}_{3}$}


\put(4.5,0.5){$\tilde{\Omega}$}


\put(1,-0.8) {Domaine de d\'efinition de $\tilde{P}$}

\end{picture}

\bigskip \bigskip \bigskip \bigskip \hrule \bigskip



\section{Calcul de la monodromie}

On va montrer comment le th\'eor\`eme pr\'ec\'edent permet de calculer la 
monodromie de l'\'equation diff\'erentielle. On d\'eterminera celle-ci sous
la forme de l'effet sur $\tilde{X}^{(0)}$ du prolongement analytique le long 
de petit lacets positifs autour de $\tilde{z}_{1},...,\tilde{z}_{r}$. On ne
s'occupera pas de la monodromie en $0$ et en $\infty$, qui se lit directement 
sur les structures de Jordan de $\tilde{B}(0)$ et de $\tilde{B}(\infty)$
(exposants venant des valeurs propres, blocs logarithmiques venant des blocs
unipotents). Ceci est conforme \`a la tradition (depuis Riemann !) qui 
consid\`ere la monodromie de $z^{\gamma}$ et celle de $\log z$ comme des
donn\'ees immanentes et ne cherche la combinatoire g\'eom\'etrique qu'aux
autres singularit\'es.\\

Pour la premi\`ere fois, nous devrons choisir $q_{0}$: ceci, dans le but
d'emp\^echer la superposition des coupures. La condition que nous imposons
est que les $(r+1)$ spirales 
$\tilde{z}_{i} q_{0}^{\overline{R}}$ ($0 \leq i \leq r$) soient deux \`a
deux distinctes, donc ne se rencontrent qu'en $0$ et en $\infty$; de 
fa\c{c}on \'equivalente: aucun $\frac{\tilde{z}_{i}}{\tilde{z}_{j}}$
($0 \leq i,j \leq r$) n'appartient \`a $q_{0}^{\mathbf{R}}$. Autrement dit,
si $q_{0} = e^{- 2 \imath \pi \tau_{0}}$, et si l'on note $\log z$ 
l'ensemble des logarithmes de $z$ dans $\mathbf{C}$,
$\tau_{0} \not \in \underset{0 \leq i,j \leq r}{\bigcup}
\mathbf{R} \frac{\log \tilde{z}_{i} - \log \tilde{z}_{j}}{2 \imath \pi}$.
Cette condition est \'evidemment v\'erifi\'ee pour $q$ g\'en\'erique.\\

Nous ordonnons les $\tilde{z}_{i}$ en les r\'eindexant de la fa\c{c}on
suivante: lorsque l'on parcourt dans le sens positif un petit cercle 
autour de $0$, on rencontre les spirales logarithmiques 
$\tilde{z}_{0} q_{0}^{\overline{\mathbf{R}}} = q_{0}^{\overline{\mathbf{R}}}$,
$\tilde{z}_{1} q_{0}^{\overline{\mathbf{R}}}$, ...,
$\tilde{z}_{r} q_{0}^{\overline{\mathbf{R}}}$ dans cet ordre.\\

Le m\^eme parcours fait rencontrer successivement les composantes connexes
de $\tilde{\Omega}'$: on nomme $\tilde{U}_{i}$ celle dont la fronti\`ere est
form\'ee de
$\tilde{z}_{i} q_{0}^{\overline{\mathbf{R}}}$ et de
$\tilde{z}_{i+1} q_{0}^{\overline{\mathbf{R}}}$ (en convenant que
$\tilde{z}_{r+1} = \tilde{z}_{0}$).



\bigskip \hrule \bigskip

\unitlength = 1cm

\underline{Figure 7 a}\\

\bigskip

\begin{picture} (15,3)


\qbezier (4,0)(9,1.5)(4,3)     

\qbezier (2,0)(-3,1.5)(2,3)      

\qbezier (2,3)(3,3.3)(4,3)   

\qbezier (2,0)(3,-0.3)(4,0)  


\put(0.5,1.5){\circle*{0.1}}
\put(0.3,1.6){$0$}

\put(5.5,1.5){\circle*{0.1}}
\put(5.35,1.6){$\infty$}

\put(2,2.5){\circle*{0.1}}
\put(2,2.6){$\tilde{z}_{2}$}

\put(3.5,2){\circle*{0.1}}
\put(3.5,2.1){$\tilde{z}_{1}$}

\put(2.5,0.5){\circle*{0.1}}
\put(2.5,0.6){$\tilde{z}_{3}$}


\qbezier (2,2.5) (4,3) (5.5,1.5)

\qbezier (3.5,2) (4.5,2) (5.5,1.5)

\qbezier (2.5,0.5) (4,0.5) (5.5,1.5)

\qbezier (0.5,1.5) (3,1) (5.5,1.5)


\put (2.5,0.1){\circle*{0.1}}
\put (2.6,0.1){$a$}
\put (2.5,0.9){\circle*{0.1}}
\put (2.6,0.9){$b$}
\qbezier (2.5,0.1) (1.5,0.5) (2.5,0.9)
\put (2,0.5){\vector(0,-1){0.1}}
\put (1.5,0.5){$\gamma_{2}$}


\put(1,1.6){$\tilde{\Omega}'_{0}$}


\put(1,-0.8) {Promenade de $\tilde{X}^{(0)}$}

\end{picture}

\bigskip \bigskip \bigskip \hrule \bigskip \newpage


\bigskip \hrule \bigskip

\unitlength = 1cm

\underline{Figure 7 b}\\

\bigskip

\begin{picture} (15,3)


\qbezier (4,0)(9,1.5)(4,3)     

\qbezier (2,0)(-3,1.5)(2,3)      

\qbezier (2,3)(3,3.3)(4,3)   

\qbezier (2,0)(3,-0.3)(4,0)  


\put(0.5,1.5){\circle*{0.1}}
\put(0.3,1.6){$0$}

\put(5.5,1.5){\circle*{0.1}}
\put(5.35,1.6){$\infty$}

\put(2,2.5){\circle*{0.1}}
\put(2,2.6){$\tilde{z}_{2}$}

\put(3.5,2){\circle*{0.1}}
\put(3.5,2.1){$\tilde{z}_{1}$}

\put(2.5,0.5){\circle*{0.1}}
\put(2.5,0.6){$\tilde{z}_{3}$}


\qbezier (0.5,1.5) (0.75,2) (02,2.5)

\qbezier (0.5,1.5) (1.5,2) (03.5,2)

\qbezier (0.5,1.5) (1.5,0.5) (02.5,0.5)

\qbezier (0.5,1.5) (03,1) (05.5,1.5)


\put (02.5,0.1){\circle*{0.1}}
\put (02.3,0.1){$a$}
\put (02.5,0.9){\circle*{0.1}}
\put (02.3,0.9){$b$}
\qbezier (02.5,0.1) (03.5,0.5) (02.5,0.9)
\put (03,0.5){\vector(0,1){0.1}}
\put (03.2,0.5){$\gamma_{1}$}


\put(04.5,1.6){$\tilde{\Omega}'_{\infty}$}


\put(1,-0.8) {Promenade de $\tilde{X}^{(\infty)}$}

\end{picture}

\bigskip \bigskip \bigskip \hrule \bigskip \bigskip

Consid\'erons un petit lacet positif $\mathbf{\gamma}$ autour de 
$\tilde{z}_{j}$ ($1 \leq j \leq r$) de base $a \in \tilde{U}_{j-1}$; on peut le
supposer compos\'e d'un chemin $\mathbf{\gamma_{1}}$ dans 
$\tilde{\Omega}_{\infty}$ partant de $a$ et arrivant en $b \in \tilde{U}_{j}$,
suivi d'un chemin $\mathbf{\gamma_{2}}$ dans $\tilde{\Omega}_{0}$
partant de $b$ et arrivant en $a$. Le prolongement analytique le long de 
$\mathbf{\gamma_{1}}$ transforme $\tilde{X}^{(0)}$ en
$\tilde{X}^{(\infty)}\tilde{P}_{j-1}$; celui le long de $\mathbf{\gamma_{2}}$ 
transforme $\tilde{X}^{(\infty)}$ en $\tilde{X}^{(0)}\tilde{P}_{j}^{-1}$ et 
donc $\tilde{X}^{(\infty)}\tilde{P}_{j-1}$ en 
$\tilde{X}^{(0)}\tilde{P}_{j}^{-1}\tilde{P}_{j-1}$.\\

\textbf{\underline{Th\'eor\`eme}}\\ 

\emph{La matrice de monodromie de (2) 
autour de $\tilde{z}_{j}$ dans la base $\tilde{X}^{(0)}$ est
$\tilde{P}_{j}^{-1}\tilde{P}_{j-1}$.}\hfill $\Box$\\



\bigskip \hrule \bigskip

\unitlength = 1cm

\underline{Figure 8}\\

\bigskip

\begin{picture} (15,3)


\qbezier (4,0)(9,1.5)(4,3)     

\qbezier (2,0)(-3,1.5)(2,3)      

\qbezier (2,3)(3,3.3)(4,3)   

\qbezier (2,0)(3,-0.3)(4,0)  


\put(0.5,1.5){\circle*{0.1}}
\put(0.3,1.6){$0$}

\put(5.5,1.5){\circle*{0.1}}
\put(5.35,1.6){$\infty$}

\put(2,2.5){\circle*{0.1}}
\put(2,2.6){$\tilde{z}_{2}$}

\put(3.5,2){\circle*{0.1}}
\put(3.5,2.1){$\tilde{z}_{1}$}

\put(2.5,0.5){\circle*{0.1}}
\put(2.5,0.6){$\tilde{z}_{3}$}


\qbezier (2,2.5) (4,3) (5.5,1.5)

\qbezier (3.5,2) (4.5,2) (5.5,1.5)

\qbezier (2.5,0.5) (4,0.5) (5.5,1.5)

\qbezier (0.5,1.5) (3,1) (5.5,1.5)

\qbezier (0.5,1.5) (0.75,2) (2,2.5)

\qbezier (0.5,1.5) (1.5,2) (3.5,2)

\qbezier (0.5,1.5) (1.5,0.5) (2.5,0.5)


\put(3,1.6){$\tilde{U}_{0}$}

\put(3,2.2){$\tilde{U}_{1}$}

\put(3,2.7){$\tilde{U}_{2}$}

\put(3,0.8){$\tilde{U}_{3}$}


\put (2.5,0.5){\circle{0.8}}
\put (2.5,0.9){\vector(-1,0){0.1}}
\put (3,0.2){$\gamma$}
\put (2.5,0.1){\circle*{0.1}}
\put (2.8,0){$a$}


\put(4.5,0.5){$\tilde{\Omega}$}


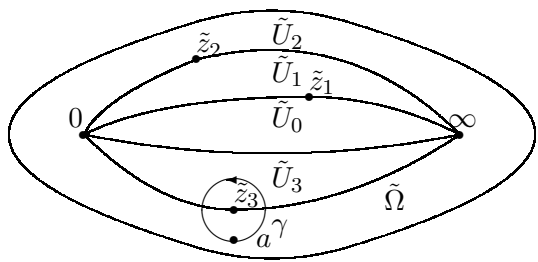
\put(1,-0.8) {Monodromie de $\tilde{X}^{(0)}$}

\end{picture}

\bigskip \bigskip \bigskip \hrule \bigskip

\textbf{\underline{Remarque}}\\

Dans \cite{Fruchard}, il est prouv\'e dans le cadre de l'analyse non-standard 
que, si $\epsilon$ est infiniment petit $>0$, et si $f$ est une fonction
analytique complexe $\epsilon$-p\'eriodique limit\'ee au voisinage d'un point 
$x_{0}$, il existe une bande horizontale bord\'ee par des droites
d'ordonn\'ees $\alpha$ et $\beta$ 
telles que $\alpha << Im(x_{0}) << \beta$ sur laquelle $f$ est constante \`a 
$\pounds e^{-C/\epsilon}$ pr\`es, o\`u $\pounds$ est limit\'e et o\`u
$C >> 0$.\\

Cela semble pouvoir \^etre appliqu\'e \`a $P_{\epsilon} - \tilde{P}$ pour
conclure que la convergence est exponentiellement rapide: ceci rejoint
l'estimation d'erreur de 3.1.2.



\section{Trois exemples complets de d\'eformations}

On tire ici des informations g\'eom\'etriques des comportements asymptotiques
obtenus au 3. Pr\'ecis\'ement, partant d'une \'equation diff\'erentielle
fuchsienne, ou d'une de ses solutions particuli\`erement int\'eressante,
on montre comment les formules de connexion d'une famille d'\'equations aux
$q$-diff\'erences qui en sont une d\'eformation permettent d'en calculer la
monodromie.\\

Ces exemples sont plus une illustration qu'une application de la th\'eorie.
Ils ont en effet \'et\'e \'etudi\'es les premiers (\`a l'exception du dernier,
``non fuchsien'') et ce sont eux qui ont r\'ev\'el\'e le ph\'enom\`ene que
nous avons d\'ecrit plus haut; l'\'enonc\'e plus g\'en\'eral est venu
ensuite. Aussi, sont ils \'etudi\'es ``\`a la main''. Tout d'abord, on utilise
plut\^ot des \'equations scalaires d'ordre $n$ que des \'equations 
vectorielles de dimension $n$
et d'ordre $1$: dans la pratique les calculs sont plus faciles ainsi. On
choisit les solutions canoniques et donc la matrice de connexion selon le
contexte (par exemple historique dans le cas hyperg\'eom\'etrique basique !).
Enfin, on passe librement de la forme scalaire \`a la forme matricielle (par 
exemple des formules de connexion \`a la matrice de connexion) et
r\'eciproquement.\\

Nous avons d\'ej\`a remarqu\'e \`a plusieurs reprises que les propri\'et\'es
d\'emontr\'ees sont ind\'ependantes du choix d'une famille particuli\`ere de
caract\`eres soumis \`a certaines conditions g\'en\'erales. Nous jouerons
sur cette libert\'e de choix chaque fois que cela sera commode pour obtenir
un comportement g\'eom\'etrique plus frappant.


\subsection{D\'eformation de $(1-z/z_{0})^{\alpha}$}

On consid\`ere ici l'\'equation diff\'erentielle fuchsienne:
$$
\delta \tilde{f} = \tilde{b} \tilde{f} \qquad
\tilde{b}(z) = \frac{- \alpha z/z_{0}}{1 - z/z_{0}} = 
\frac{\alpha}{1 - w/w_{0}} \qquad
(w = 1/z , w_{0} = 1/z_{0})
$$

La m\'ethode de Frobenius (o\`u l'on remplace $z^{\gamma}$ par $(-z)^{\gamma}$
et $\log(z)$ par $\log(-z)$) donne les solutions canoniques en $0$ et en
$\infty$:
$$
\begin{cases}
\tilde{f}^{(0)}(z) = (1 - z/z_{0})^{\alpha} \\
\tilde{f}^{(\infty)} = (-w)^{- \alpha} (1 - w/w_{0})^{\alpha} =
(-z)^{\alpha} (1 - z_{0}/z)^{\alpha}
\end{cases}
$$

Pour les rendre uniformes, on pratique des coupures d\'etermin\'ees par les 
exposants et par les p\^oles de $\tilde{b}$ sur $\mathbf{C}^{*}$: on
consid\`ere $\tilde{f}^{(0)}$ comme d\'efinie sur 
$U_{0} = \mathbf{S} - z_{0} q_{0}^{\overline{\mathbf{R}}^{+}}$ et 
$\tilde{f}^{(\infty)}$ comme d\'efinie sur 
$U_{\infty} = 
\mathbf{S} - z_{0}q_{0}^{\overline{\mathbf{R}}^{-}}- q_{0}^{\overline{\mathbf{R}}}$: ce sont bien deux ouverts simplement connexes \'evitant respectivement
les singularit\'es de $\tilde{f}^{(0)}$ et de $\tilde{f}^{(\infty)}$.
On supposera $z_{0} \not \in q_{0}^{\overline{\mathbf{R}}}$ pour que les
coupures ne se superposent pas (c'est la raison pour laquelle on n'a pas
pris tout simplement $z_{0} = 1$). Cela nous permettra d'utiliser notre
famille habituelle de caract\`eres.\\

\newpage


\hrule \bigskip

\unitlength = 1cm

\underline{Figure 9 a}\\

\bigskip

\begin{picture} (15,3)


\qbezier (4,0)(9,1.5)(4,3)     

\qbezier (2,0)(-3,1.5)(2,3)      

\qbezier (2,3)(3,3.3)(4,3)   

\qbezier (2,0)(3,-0.3)(4,0)  


\put(0.5,1.5){\circle*{0.1}}
\put(0.3,1.6){$0$}

\put(5.5,1.5){\circle*{0.1}}
\put(5.35,1.6){$\infty$}

\put(3.5,2){\circle*{0.1}}
\put(3.5,2.1){$z_{0}$}


\qbezier (3.5,2) (4.5,2) (5.5,1.5)


\put(1,1.6){$U_{0}$}


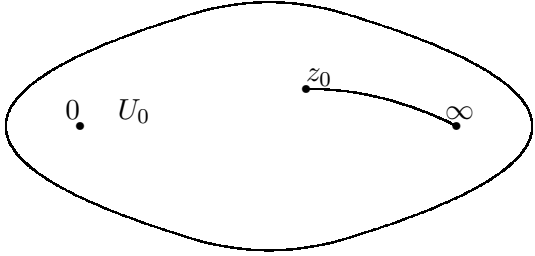
\put(1,-0.8) {Domaine de $\tilde{f}^{(0)}$}

\end{picture}

\bigskip \bigskip \bigskip \hrule \bigskip


\unitlength = 1cm

\underline{Figure 9 b}\\

\bigskip

\begin{picture} (15,3)


\qbezier (4,0)(9,1.5)(4,3)     

\qbezier (2,0)(-3,1.5)(2,3)      

\qbezier (2,3)(3,3.3)(4,3)   

\qbezier (2,0)(3,-0.3)(4,0)  


\put(0.5,1.5){\circle*{0.1}}
\put(0.3,1.6){$0$}

\put(5.5,1.5){\circle*{0.1}}
\put(5.35,1.6){$\infty$}

\put(3.5,2){\circle*{0.1}}
\put(3.5,2.1){$z_{0}$}


\qbezier (0.5,1.5) (1.5,2) (03.5,2)

\qbezier (0.5,1.5) (03,1) (05.5,1.5)


\put(04.5,1.6){$U_{\infty}$}

\put (2.4,0.8) {$q_{0}^{\overline{\mathbf{R}}}$}


\put(1,-0.8) {Domaine de $\tilde{f}^{(\infty)}$}

\end{picture}

\bigskip \bigskip \bigskip \hrule \bigskip

On d\'eformera l'\'equation diff\'erentielle ci-dessus en la famille 
d'\'equations aux $q$-diff\'erences fuchsiennes:
$$
\sigma_{q} f = a_{\epsilon} f
$$
avec $a_{\epsilon} \in \mathbf{k}^{*}$ et m\^eme
$a_{\epsilon}(0) , a_{\epsilon}(\infty) \in \mathbf{C}^{*}$. Les conditions
naturelles \`a imposer pour mieux sp\'ecifier cette famille sont les
suivantes:

\begin{enumerate}

\item{$\frac{a_{\epsilon} - 1}{q - 1} \underset{q \to 1}{\to} \tilde{b}$: ceci
parce que l'on a en t\^ete la formule heuristique
$\delta_{q} := \frac{\sigma_{q} - 1}{q - 1} \underset{q \to 1}{\to} \delta$.
Cela permet d'esp\'erer qu'une famille de solutions $f_{\epsilon}$ bien
choisies v\'erifieront
$\frac{f_{\epsilon} - 1}{q - 1} \underset{q \to 1}{\to} \tilde{f}$, solution
de l'\'equation diff\'erentielle. Il y a en particulier des contraintes au 
niveau des exposants:
$\frac{a_{\epsilon}(0) - 1}{q - 1} \underset{q \to 1}{\to} \tilde{b}(0)$
et similairement en $\infty$.}

\item{Les p\^oles de $a_{\epsilon}$ doivent en outre tendre vers ceux de
$\tilde{b}$ sur $\mathbf{C}^{*}$ (ici: vers $z_{0}$); ceci, afin que les
demi-spirales logarithmiques discr\`etes de singularit\'es des solutions 
aillent se condenser sur les coupures, comme on l'a vu se produire pour les
caract\`eres au 3.1.3.}

\end{enumerate}

Nous d\'ecidons de prendre 
$a_{\epsilon}(z) = \frac{1 - q^{\alpha} z/z_{0}}{1 - z/z_{0}}$. Cela revient 
\`a poser 
$b_{\epsilon}(z) = \frac{q^{\alpha} - 1}{q - 1} 
\frac{- z/z_{0}}{1 - z/z_{0}}$ et \`a \'etudier la famille d'\'equations
$\delta_{q} f = b_{\epsilon} f$. Les conditions pr\'ec\'edentes sont
manifestement satisfaites. De plus, les calculs seront simples
grace aux r\'esultats de 3.1.5.

\subsubsection{Solution canonique en $0$}

Il n'y a pas d'exposant (autrement dit, il vaut $1$). On trouve donc
imm\'ediatement:
$$
f_{\epsilon}^{(0)}(z) = 
\frac{\Theta_{q}^{+}(q^{\alpha}z/z_{0})}{\Theta_{q}^{+}(z/z_{0})} \div
\frac{1 - q^{\alpha} z/z_{0}}{1 - z/z_{0}}
$$
Cette solution est d\'efinie sur $\mathbf{C} - q_{0}^{\mathbf{N}^{*}}$, et
en particulier sur $U_{0}$. De plus, ses z\'eros disparaissent de tout compact
de $U_{0}$ lorsque $q \to 1$.

\subsubsection{Solution canonique en $\infty$}

La fonction inconnue $g(w) = f(1/w)$ doit y v\'erifier:
$$
\sigma_{q}g(w) = (a_{\epsilon}(1/qw))^{-1} g(w) =
q^{-\alpha} \frac{1 - q w/w_{0}}{1 - q w/q^{\alpha} w_{0}} g(w)
$$
Il y a l'exposant $q^{-\alpha}$ et l'on trouve la solution
$$
f_{\epsilon}^{(\infty)}(z) = 
e_{q,q^{-\alpha}}(w) 
\frac{\Theta_{q}^{+}(w/w_{0})}{\Theta_{q}^{+}(w/q^{\alpha}w_{0})}
$$
qui est d\'efinie sur $U_{\infty}$. Lorsque $q \to 1$,
ses z\'eros disparaissent de tout compact de $U_{\infty}$.

\subsubsection{Matrice de connexion}

On devrait ici dire nombre de connexion ! C'est
$$
p_{\epsilon}(z) = \frac{f_{\epsilon}^{(0)}(z)}{f_{\epsilon}^{(\infty)}(z)} =
\frac{1}{e_{q^{-\alpha}}(1/z)} 
\frac{\Theta_{q}(q^{\alpha} z/z_{0})}{\Theta_{q}(z/z_{0})}
$$
Elle est d\'efinie sur tout 
$U = U_{0} \cap U_{\infty} = 
\mathbf{S} - z_{0}q_{0}^{\overline{\mathbf{R}}}- q_{0}^{\overline{\mathbf{R}}}$,
qui a deux composantes connexes: $U'$, \`a droite de la spirale 
$z_{0}q_{0}^{\overline{\mathbf{R}}}$ (lorsqu'on la parcourt de $0$ \`a 
$\infty$), et $U''$ \`a sa gauche.

\subsubsection{Confluence}

Lorsque $q \to 1$, $f_{\epsilon}^{(0)} \to \tilde{f}^{(0)}$ d'apr\`es 3.1.5,
$f_{\epsilon}^{(\infty)} \to \tilde{f}^{(\infty)}$ d'apr\`es 3.1.3 et 3.1.5,
et donc 
$p_{\epsilon} \to \tilde{p} = \frac{\tilde{f}^{(0)}}{\tilde{f}^{(\infty)}}$;
par application de 3.1.2 et 3.1.3, on voit que
$\tilde{p}(z) = (-z)^{-\alpha} (-z/z_{0})^{\alpha}$. Cette fonction est
localement constante sur $U$ et prend donc une valeur $\tilde{p}'$ sur
$U'$ et une valeur $\tilde{p}''$ sur $U''$. Ces valeurs sont li\'ees par
la relation $\tilde{p}'' = \tilde{p}' e^{- 2 \imath \pi \alpha}$: en effet,
seul le facteur $(-z/z_{0})^{\alpha}$ est affect\'e lors du franchissement de
la coupure $z_{0}q_{0}^{\overline{\mathbf{R}}}$ de droite \`a gauche.

\newpage


\hrule \bigskip

\unitlength = 1cm

\underline{Figure 9 c}\\

\bigskip

\begin{picture} (15,3)


\qbezier (4,0)(9,1.5)(4,3)     

\qbezier (2,0)(-3,1.5)(2,3)      

\qbezier (2,3)(3,3.3)(4,3)   

\qbezier (2,0)(3,-0.3)(4,0)  


\put(0.5,1.5){\circle*{0.1}}
\put(0.3,1.6){$0$}

\put(5.5,1.5){\circle*{0.1}}
\put(5.35,1.6){$\infty$}

\put(3.5,2){\circle*{0.1}}
\put(3.5,2.1){$z_{0}$}


\qbezier (0.5,1.5) (1.5,2) (03.5,2)

\qbezier (0.5,1.5) (03,1) (05.5,1.5)


\put (03.5,1.6){\circle*{0.1}}
\put (03.3,1.6){$a$}
\put (03.5,2.4){\circle*{0.1}}
\put (03.3,2.4){$b$}
\qbezier (03.5,1.6) (04.5,2.0) (03.5,2.4)
\put (04,2.0){\vector(0,1){0.1}}
\put (04.2,2.0){$\gamma_{1}$}

\put (2.4,0.8) {$q_{0}^{\overline{\mathbf{R}}}$}
\put (2,2.5) {$U''$}
\put (2,1.5) {$U'$}


\put(04.5,1.6){$U_{\infty}$}


\put(1,-0.8) {Promenade de $\tilde{f}^{(0)}$}

\end{picture}

\bigskip \bigskip \bigskip \hrule \bigskip


\unitlength = 1cm

\underline{Figure 9 d}\\

\bigskip

\begin{picture} (15,3)


\qbezier (4,0)(9,1.5)(4,3)     

\qbezier (2,0)(-3,1.5)(2,3)      

\qbezier (2,3)(3,3.3)(4,3)   

\qbezier (2,0)(3,-0.3)(4,0)  


\put(0.5,1.5){\circle*{0.1}}
\put(0.3,1.6){$0$}

\put(5.5,1.5){\circle*{0.1}}
\put(5.35,1.6){$\infty$}

\put(3.5,2){\circle*{0.1}}
\put(3.5,2.1){$z_{0}$}


\qbezier (3.5,2) (4.5,2) (5.5,1.5)


\put (3.5,1.6){\circle*{0.1}}
\put (3.6,1.6){$a$}
\put (3.5,2.4){\circle*{0.1}}
\put (3.6,2.4){$b$}
\qbezier (3.5,1.6) (2.5,2.0) (3.5,2.4)
\put (3,2.0){\vector(0,-1){0.1}}
\put (2.5,2.0){$\gamma_{2}$}

\put (2,2.5) {$U''$}
\put (2,1.5) {$U'$}


\put(1,1.6){$U_{0}$}


\put(1,-0.8) {Retour de la promenade}

\end{picture}

\bigskip \bigskip \bigskip \hrule \bigskip

\subsubsection{Monodromie}

On calculera la monodromie autour de $z_{0}$: celles en $0$ et en $\infty$
se lisent directement sur les exposants. On regarde l'effet sur 
$\tilde{f}^{(0)}$ d'un petit lacet positif autour de $z_{0}$, partant d'un 
point $a$ de $U'$.\\

La premi\`ere partie $\gamma_{1}$ du lacet est dans $U_{\infty}$. Comme 
$\tilde{f}^{(\infty)}$ y est uniforme et que, au d\'epart,
$\tilde{f}^{(0)} = \tilde{f}^{(\infty)} \tilde{p}'$, le r\'esultat du
prolongement analytique de $\tilde{f}^{(0)}$ le long de cette premi\`ere 
partie du lacet est $\tilde{f}^{(\infty)} \tilde{p}'$.\\

La seconde partie $\gamma_{2}$
du lacet est dans $U_{0}$, o\`u $\tilde{f}^{(0)}$ est
uniforme. Comme $\tilde{f}^{(\infty)} = \tilde{f}^{(0)}(\tilde{p}'')^{-1}$
au d\'epart $b$ de ce second chemin, le r\'esultat du prolongement analytique 
de $\tilde{f}^{(0)}$ le long du lacet complet est
$\tilde{f}^{(0)}(\tilde{p}'')^{-1}\tilde{p}'$.\\

L'effet sur la base canonique $\tilde{f}^{(0)}$ de la monodromie autour de 
$z_{0}$ est donc la multiplication \`a droite par 
$(\tilde{p}'')^{-1}\tilde{p}'$: bien qu'en dimension $1$ tout commute, 
l'ordre correct des multiplications a \'et\'e syst\'ematiquement respect\'e
dans les calculs. On trouve donc que c'est la multiplication par
$e^{2 \imath \pi \alpha}$, comme il \'etait souhaitable que ce le f\^ut !


\subsection{D\'eformation de la fonction hyperg\'eom\'etrique}

\emph{Les formules de connexion pour les 
fonctions hyperg\'eom\'etriques basiques donnent, \`a la limite, celles des 
fonctions hyperg\'eom\'etriques classiques.}\\

La s\'erie hyperg\'eom\'etrique (\emph{voir} \cite{WW},\emph{voir} \cite{GP})
$$
F(\alpha,\beta,\gamma;z) = \sum_{n \geq 0}
\frac{(\alpha)_{n} (\beta)_{n}}{(\gamma)_{n} (1)_{n}}z^{n}
\qquad (x)_{n} = \prod_{0 \leq k \leq n-1}(x+k) = 
\frac{\Gamma(x+n)}{\Gamma(x)}
$$
est solution de l'\'equation diff\'erentielle \`a coefficients rationnels: 
$$
\delta^{2}F - \tilde{\lambda} \delta F + \tilde{\mu} F = 0
\quad \text{avec: }
\begin{cases}
\tilde{\lambda} = \frac{((\alpha + \beta)z + (1 - \gamma))}{(1-z)} \\
\tilde{\mu} = \frac{\alpha \beta z}{(1-z)}
\end{cases}
$$

Les singularit\'es de cette \'equation 
sont $0$,$1$ et $\infty$ et elles sont fuchsiennes. La m\'ethode de Frobenius 
en fournit des syst\`eme fondamentaux de solutions en $0$ et en $\infty$:
$$ 
\begin{cases}
F(\alpha,\beta,\gamma;z) \\
z^{1-\gamma}F(\alpha - \gamma + 1,\beta - \gamma +1,2 - \gamma;z)
\end{cases}
\quad
\begin{cases}
(1/z)^{\alpha}F(\alpha,\alpha - \gamma +1,\alpha - \beta + 1;1/z) \\
(1/z)^{\beta}F(\beta,\beta - \gamma +1,\beta - \alpha + 1;1/z)
\end{cases}
$$
Ces fonctions a priori multivalu\'ees sont respectivement rendues uniformes par
les coupures
$[1,\infty]$ et $[-\infty,0]$ de la sph\`ere de Riemann.\\

Le probl\`eme de d\'eformer cette \'equation diff\'erentielle en une famille
d'\'equations aux $q$-diff\'erences et, r\'eciproquement, de d\'eduire sa
monodromie du comportement asymptotique des matrices de connexion est 
doublement simplifi\'e parce que

\begin{itemize}

\item{Une famille de fonctions d\'eformant $F$ est connue: les fonctions
hyperg\'eom\'etriques basiques.}

\item{La monodromie de $F$ est classiquement calcul\'ee \`a l'aide d'une
formule de connexion reliant les syst\`emes fondamentaux de solutions en
$0$ et en $\infty$; de plus, la formule de Mellin-Barnes, qui permet ce
calcul, a \'et\'e \'etendue au cas hyperg\'eom\'etrique basique par Watson.}

\end{itemize}

La s\'erie hyperg\'eom\'etrique \emph{classique} a \'et\'e g\'en\'eralis\'ee 
par Heine, puis Ramanujan (\emph{voir} \cite{GR}, \cite{Ramanujan}) en la 
s\'erie hyperg\'eom\'etrique \emph{basique}
$$
\Phi(a,b,c;q,z) = \sum_{n \geq 0}
\frac{(a;p)_{n} (b;p)_{n}}{(c;p)_{n} (p;p)_{n}}z^{n}
\quad |q| > 1 \;,\; p = q^{-1} \;,\; 
(x;p)_{n} = \prod_{0 \leq k \leq n-1}(1 - x p^{k})
$$

Lorsque $q \to 1$, la s\'erie hyperg\'eom\'etrique basique de param\`etres
$a = q^{\alpha} , b = q^{\beta} , c = q^{\gamma}$
tend vers la s\'erie hyperg\'eom\'etrique classique de param\`etres
$\alpha , \beta , \gamma$ coefficient par coefficient. Il s'agit de pr\'eciser
un peu la nature de cette ``confluence''.\\

Soient $a,b,c \in \mathbf{C}^{*}$. La s\'erie hyperg\'eom\'etrique basique de 
param\`etres $a , b , c$ est solution de l'\'equation aux 
$q$-diff\'erences \`a coefficients rationnels:
$$
\sigma_{q}^{2}\Phi - \lambda' \sigma_{q}\Phi + \mu' \Phi = 0
\quad \text{avec: }
\begin{cases}
\lambda' = \frac{((a + b)z - (1+c/q))}{(abz - c/q)} \\
\mu' = \frac{(z-1)}{(abz-c/q)}
\end{cases}
$$
Nous pr\'ef\'erons r\'e\'ecrire cette \'equation \`a l'aide de l'op\'erateur
$\delta_{q} = \frac{\sigma_{q} - 1}{q - 1}$ afin de rendre plus visibles les
conditions formelles de la confluence; il suffit de faire 
$\sigma_{q} = 1 + \eta \delta_{q} \; (\eta = q - 1)$.
$$
\delta_{q}^{2}\Phi - \lambda \delta_{q}\Phi + \mu \Phi = 0
\quad \text{avec: }
\begin{cases}
\lambda = \frac{\lambda' - 2}{\eta} \\
\mu = \frac{\mu' - \lambda' + 1}{\eta^{2}}
\end{cases}
$$

Notre m\'ethode attribue \`a son \'equation les syst\`emes fondamentaux
de solutions canoniques en $0$ et en $\infty$: 
$$
\begin{cases}
\Phi(a,b,c;q,z) \\
e_{q/c}(z)\Phi(aq/c,bq/c,q^{2}/c;q,z)
\end{cases}
\quad
\begin{cases}
e_{a}(1/z)\Phi(a,aq/c,aq/b;q,cq/abz) \\
e_{b}(1/z)\Phi(b,bq/c,bq/a;q,cq/abz)
\end{cases}
$$
On raisonne ici directement sur l'\'equation du second ordre, mais des
changements de variable simples montrent que les solutions ci-dessus sont
de la forme (fonction holomorphe valant $1$ en $0$) $\times$ (caract\`ere): 
ce sont bien des solutions canoniques telles que nous les avons d\'efinies 
dans le cadre matriciel. Pour cette m\^eme raison, les formules de connexion 
que nous allons tirer de \cite{GR} sont bien celles que donneraient les 
coefficients de la
matrice de connexion.\\

Le syst\`eme diff\'erentiel correspondant \`a l'\'equation 
hyperg\'eom\'etrique classique de param\`etres $\alpha$, $\beta$ et $\gamma$
a pour matrice
$$
\tilde{B} = 
\begin{pmatrix}
0 & 1\\
-\tilde{\mu} & \tilde{\lambda}\\
\end{pmatrix}
$$

Nous le d\'eformerons donc par le syst\`eme aux $q$-diff\'erences
correspondant \`a l'\'equation
hyperg\'eom\'etrique basique de param\`etres $a = q^{\alpha}$, $b = q^{\beta}$
et $c = q^{\gamma}$,de  matrice 
$$
A_{\epsilon} = 
\begin{pmatrix}
1 & q-1\\
-\frac{(\mu' - \lambda' + 1)}{\eta} & \lambda' - 1
\end{pmatrix}
$$
On a alors $A_{\epsilon} = I_{2} + \eta B_{\epsilon}$, o\`u
$$
B_{\epsilon} = 
\begin{pmatrix}
0 & 1 \\
- \mu & \lambda
\end{pmatrix}
$$
On observe que, lorsque $\epsilon \to 0^{+}$,
$(A_{\epsilon} - I_{2})/(q-1) \to \tilde{B}$. Toutes nos autres conditions sont
d'ailleurs satisfaites.\\

Les formules de connexion de Barnes-Mellin-Watson (\emph{voir} \cite{GR} 
p. 106) relient les deux syst\`emes fondamentaux de solutions:
\begin{eqnarray*}
\Phi(a,b;c;q,z) & = &
\Phi(c/b,a;c;q,q)\frac{\Theta_{q}(bz/cq)}{\Theta_{q}(abz/cq)}
\Phi(a,aq/c,aq/b;q,cq/abz)\\
                & + &
\Phi(c/a,b;c;q,q)\frac{\Theta_{q}(az/cq)}{\Theta_{q}(abz/cq)}
\Phi(b,bq/c,bq/a;q,cq/abz)
\end{eqnarray*}
Sans autre calcul, nous en tirons directement les coefficients de la premi\`ere
colonne de la matrice de connexion:
$$
P_{1,1}(z) = 
\Phi(c/b,a;c;q,q)\frac{\Theta_{q}(bz/cq)}{\Theta_{q}(abz/cq)}
\frac{1}{e_{q,a}(1/z)}
$$
$$
P_{2,1}(z) = 
\Phi(c/a,b;c;q,q)\frac{\Theta_{q}(az/cq)}{\Theta_{q}(abz/cq)}
\frac{1}{e_{q,b}(1/z)}
$$
De m\^eme, en ce qui concerne la seconde solution fondamentale:
\begin{eqnarray*}
\Phi(aq/c,bq/c;q^{2}/c;q,z) & = &
\Phi(q/b,aq/c;q^{2}/c;q,q)\frac{\Theta_{q}(bz/q^{2})}{\Theta_{q}(abz/cq)}
\Phi(a,aq/c,aq/b;q,cq/abz)\\
                & + &
\Phi(q/a,bq/c;q^{2}/c;q,q)\frac{\Theta_{q}(az/q^{2})}{\Theta_{q}(abz/cq)}
\Phi(b,bq/c,bq/a;q,cq/abz)
\end{eqnarray*}
d'o\`u l'on tire les coefficients de la deuxi\`eme colonne de la matrice de 
connexion:
$$
P_{1,2}(z) = 
\Phi(q/b,aq/c;q^{2}/c;q,q)\frac{\Theta_{q}(bz/q^{2})}{\Theta_{q}(abz/cq)}
\frac{e_{q,q/c}(z)}{e_{q,a}(1/z)}
$$
$$
P_{2,2}(z) = 
\Phi(q/a,bq/c;q^{2}/c;q,q)\frac{\Theta_{q}(az/q^{2})}{\Theta_{q}(abz/cq)}
\frac{e_{q,q/c}(z)}{e_{q,b}(1/z)}
$$
On peut alors en d\'eduire la matrice de connexion: 
$P_{\epsilon} = (D_{\epsilon}^{(\infty)})^{-1}M_{\epsilon}D_{\epsilon}^{(0)}$,
o\`u:
$$
D_{\epsilon}^{(0)} =
\begin{pmatrix}
1 & 0\\
0 & e_{q,q/c}(z)\\
\end{pmatrix}
\hspace{3cm}
D_{\epsilon}^{(\infty)} =
\begin{pmatrix}
e_{q,a}(1/z) & 0\\
0 & e_{q,b}(1/z)\\
\end{pmatrix}
$$
et o\`u la matrice $M_{\epsilon}$ est donn\'ee par
$$
\begin{pmatrix}
\Phi(c/b,a;c;q,q)\Theta_{q}(bz/cq)/\Theta_{q}(abz/cq) & 
\Phi(q/b,aq/c;q^{2}/c;q,q)\Theta_{q}(bz/q^{2})/\Theta_{q}(abz/cq) \\
\Phi(c/a,b;c;q,q)\Theta_{q}(az/cq)/\Theta_{q}(abz/cq) & 
\Phi(q/a,bq/c;q^{2}/c;q,q)\Theta_{q}(az/q^{2})/\Theta_{q}(abz/cq) \\
\end{pmatrix}
$$

Conform\'ement \`a la th\'eorie g\'en\'erale du chapitre 2:

\begin{itemize}

\item{La matrice de connexion $P_{\epsilon}$ est elliptique.}

\item{Les singularit\'es de la partie ``interne'' $M_{\epsilon}$ sont sur les
spirales logarithmiques discr\`etes d\'eduites des singularit\'es de la
matrice $A_{\epsilon}$ de l'\'equation. Ici, l'unique p\^ole est $c/qab$ 
(et on voit bien la spirale qu'il engendre) et l'unique z\'ero du 
d\'eterminant est $1$ (mais l\`a, cela ne se voit pas facilement!).}

\item{Les caract\`eres sont \`a l'ext\'erieur.}

\end{itemize}

Des r\'esultats du chapitre 3, on d\'eduit la valeur de la matrice
$\tilde{M} = \underset{\epsilon \to 0^{+}}{\lim}M_{\epsilon}$:
$$
\begin{pmatrix}
(-z)^{-\alpha}
\Gamma(\gamma)\Gamma(\beta - \alpha)/\Gamma(\beta)\Gamma(\gamma - \alpha) &
(-z)^{\gamma - \alpha + 1}
\Gamma(2 - \gamma)\Gamma(\beta - \alpha)/
\Gamma(1 - \alpha)\Gamma(1 - \gamma + \beta)\\
(-z)^{-\beta}
\Gamma(\gamma)\Gamma(\alpha - \beta)/\Gamma(\alpha)\Gamma(\gamma - \beta) &
(-z)^{\gamma - \beta + 1}
\Gamma(2 - \gamma)\Gamma(\alpha - \beta)/
\Gamma(1 - \beta)\Gamma(1 - \gamma + \alpha)\\
\end{pmatrix}
$$

Ces formules sont \emph{ind\'ependantes du choix d'une famille particuli\`ere 
de caract\`eres}.
Pour \'eviter la superposition des coupures, on remplacera les caract\`eres
utilis\'es jusque l\`a par $\Theta_{q}(-z)/\Theta_{q}(-c^{-1}z)$, dont il est
ais\'e de voir que leurs propri\'et\'es sont les m\^emes. En cons\'equence,
pour une fois, nous obtiendrons \`a la limite les vraies solutions canoniques
de Frobenius !\\

On trouve alors
$\tilde{D}^{(0)} = \underset{\epsilon \to 0^{+}}{\lim}D_{\epsilon}^{(0)}$ et
$\tilde{D}^{(\infty)} = \underset{\epsilon \to 0^{+}}{\lim}D_{\epsilon}^{(\infty)}$ sur
$\mathbf{S} - [-\infty,0]$ :
$$
\tilde{D}^{(0)} = 
\begin{pmatrix}
1 & 0\\
0 & z^{1 - \gamma}\\
\end{pmatrix}
\hspace{3cm}
\tilde{D}^{(\infty)} = 
\begin{pmatrix}
z^{-\alpha} & 0\\
0 & z^{-\beta}\\
\end{pmatrix}
$$
Ainsi 
$\tilde{P} = \underset{\epsilon \to 0^{+}}{\lim}P_{\epsilon} =
(\tilde{D}^{(\infty)})^{-1}\tilde{M}\tilde{D}^{(0)}$ prend les valeurs
$\tilde{P}^{+}$ et $\tilde{P}^{-}$
sur le demi-plan de Poincar\'e et sur son oppos\'e, ces valeurs \'etant
respectivement donn\'ees par les matrices suivantes: 
$$
\begin{pmatrix}
e^{i\pi\alpha}
\Gamma(\gamma)\Gamma(\beta - \alpha)/\Gamma(\beta)\Gamma(\gamma - \alpha) &
e^{i\pi(\alpha - \gamma + 1)}
\Gamma(2 - \gamma)\Gamma(\beta - \alpha)/
\Gamma(1 - \alpha)\Gamma(1 - \gamma + \beta)\\
e^{i\pi\beta}
\Gamma(\gamma)\Gamma(\alpha - \beta)/\Gamma(\alpha)\Gamma(\gamma - \beta) &
e^{i\pi(\beta - \gamma + 1)}
\Gamma(2 - \gamma)\Gamma(\alpha - \beta)/
\Gamma(1 - \beta)\Gamma(1 - \gamma + \alpha)\\
\end{pmatrix}
$$
$$
\begin{pmatrix}
e^{-i\pi\alpha}
\Gamma(\gamma)\Gamma(\beta - \alpha)/\Gamma(\beta)\Gamma(\gamma - \alpha) &
e^{-i\pi(\alpha - \gamma + 1)}
\Gamma(2 - \gamma)\Gamma(\beta - \alpha)/
\Gamma(1 - \alpha)\Gamma(1 - \gamma + \beta)\\
e^{-i\pi\beta}
\Gamma(\gamma)\Gamma(\alpha - \beta)/\Gamma(\alpha)\Gamma(\gamma - \beta) &
e^{-i\pi(\beta - \gamma + 1)}
\Gamma(2 - \gamma)\Gamma(\alpha - \beta)/
\Gamma(1 - \beta)\Gamma(1 - \gamma + \alpha)\\
\end{pmatrix}
$$

Notons $\tilde{X}^{(0)}$ et $\tilde{X^{(\infty)}}$ les vecteurs solutions
fondamentales de l'\'equation diff\'erentielle au voisinage de $0$ et
$\infty$ respectivement,d\'efinis sur
$\tilde{\Omega}'_{0} = \mathbf{S} - [1;\infty] - [-\infty;0]$ et sur 
$\tilde{\Omega}'_{\infty} = \mathbf{S} - [0;1] - [-\infty;0]$ 
(la premi\`ere coupure est due aux p\^oles
de l'\'equation, la seconde aux caract\`eres). Sur l'intersection
$\mathbf{S} - \overline{\mathbf{R}}$ de leurs domaines,
$\tilde{X}^{(0)} = \tilde{X^{(\infty)}} \tilde{P}$. Pour calculer la
monodromie autour de $1$, on prom\`enera $\tilde{X}^{(0)}$ le long d'un petit
lacet positif autour de $1$ et partant d'un point $a$ du demi-plan de 
Poincar\'e $\mathcal{H}$.\\

La premi\`ere partie $\gamma_{1}$ du lacet est dans le domaine de 
$\tilde{X}^{(0)}$, qui se
prolonge donc en $\tilde{X}^{(0)}$. On arrive alors en $b \in - \mathcal{H}$, 
o\`u $\tilde{X}^{(0)} = \tilde{X}^{(\infty)} \tilde{P}^{-}$.\\

La deuxi\`eme partie $\gamma_{2}$ du lacet est dans le domaine de 
$\tilde{X}^{(\infty)}$;
aussi, $\tilde{X}^{(\infty)} \tilde{P}^{-}$ se prolonge en
$\tilde{X}^{(\infty)} \tilde{P}^{-}$. On arrive dans $\mathcal{H}$, o\`u
$\tilde{X}^{(\infty)} = \tilde{X}^{(0)} (\tilde{P}^{+})^{-1}$. On obtient
ainsi, comme prolongement de $\tilde{X}^{(0)}$, 
$\tilde{X}^{(0)} (\tilde{P}^{+})^{-1}\tilde{P}^{-}$.\\

\newpage


\hrule \bigskip

\unitlength = 1cm

\underline{Figure 10 a}\\

\bigskip

\begin{picture} (15,3)


\qbezier (4,0)(9,1.5)(4,3)     

\qbezier (2,0)(-3,1.5)(2,3)      

\qbezier (2,3)(3,3.3)(4,3)   

\qbezier (2,0)(3,-0.3)(4,0)  


\put(0.5,1.5){\circle*{0.1}}
\put(0.3,1.6){$0$}

\put(5.5,1.5){\circle*{0.1}}
\put(5.35,1.6){$\infty$}

\put(3.5,2){\circle*{0.1}}
\put(3.5,2.1){$1$}


\qbezier (3.5,2) (4.5,2) (5.5,1.5)

\qbezier (0.5,1.5) (3,1) (5.5,1.5)


\put (3.5,1.6){\circle*{0.1}}
\put (3.6,1.6){$b$}
\put (3.5,2.4){\circle*{0.1}}
\put (3.6,2.4){$a$}
\qbezier (3.5,1.6) (2.5,2.0) (3.5,2.4)
\put (3,2.0){\vector(0,-1){0.1}}
\put (2.5,2.0){$\gamma_{1}$}

\put (2.4,0.8) {$[-\infty,0]$}
\put (2,2.5) {$\mathcal{H}$}
\put (2,1.5) {-$\mathcal{H}$}


\put(1,1.6){$\tilde{\Omega}'_{0}$}


\put(1,-0.8) {Promenade de $\tilde{X}^{(0)}$}

\end{picture}

\bigskip \bigskip \bigskip \hrule \bigskip 


\unitlength = 1cm

\underline{Figure 10 b}\\

\bigskip

\begin{picture} (15,3)


\qbezier (4,0)(9,1.5)(4,3)     

\qbezier (2,0)(-3,1.5)(2,3)      

\qbezier (2,3)(3,3.3)(4,3)   

\qbezier (2,0)(3,-0.3)(4,0)  


\put(0.5,1.5){\circle*{0.1}}
\put(0.3,1.6){$0$}

\put(5.5,1.5){\circle*{0.1}}
\put(5.35,1.6){$\infty$}

\put(3.5,2){\circle*{0.1}}
\put(3.5,2.1){$1$}


\qbezier (0.5,1.5) (1.5,2) (03.5,2)

\qbezier (0.5,1.5) (03,1) (05.5,1.5)


\put (03.5,1.6){\circle*{0.1}}
\put (03.3,1.6){$b$}
\put (03.5,2.4){\circle*{0.1}}
\put (03.3,2.4){$a$}
\qbezier (03.5,1.6) (04.5,2.0) (03.5,2.4)
\put (04,2.0){\vector(0,1){0.1}}
\put (04.2,2.0){$\gamma_{2}$}

\put (2.4,0.8) {$[-\infty,0]$}
\put (2,2.5) {$\mathcal{H}$}
\put (2,1.5) {-$\mathcal{H}$}


\put(04.5,1.6){$\tilde{\Omega}'_{\infty}$}


\put(1,-0.8) {Retour de la promenade}

\end{picture}

\bigskip \bigskip \bigskip \hrule \bigskip

La matrice de monodromie en le p\^ole $1$ est donc
$(\tilde{P}^{+})^{-1}\tilde{P}^{-}$, ce qui est la valeur
classique: voir, par exemple, \cite{GP} pp 111-114\footnote{
Il faut cependant prendre garde \`a une erreur de signe dans cette 
r\'ef\'erence: avec la 
base canonique utilis\'ee, d\'efinie p. 38 de loc.cit., on
a bien la formule de connexion 4.6.3, p. 111; mais, avec le choix des chemins
de la figure 2.5, p. 78, l'expression de $A_{\infty}$ dans la
formule 4.6.12, p. 114 devrait faire intervenir $\epsilon(\alpha)$ et
$\epsilon(\beta)$ au lieu de $\epsilon(-\alpha)$ et $\epsilon(-\beta)$.}.


\subsection{D\'eformation de $\log(1 - z/z_{0})$}

On va \`a pr\'esent \'etudier le cas d'un syst\`eme dont la monodromie est 
unipotente, 
mais qui rentre tout de m\^eme dans le cas non-r\'esonnant. L'\'equation
diff\'erentielle de d\'epart est
$$ 
\frac{d^{2}f}{dz^{2}} = \frac{1}{z_{0} - z} \frac{df}{dz}
\qquad \Longleftrightarrow \qquad
\delta^{2}f = \frac{1}{1 - \frac{z}{z_{0}}} \delta f
$$

Les singularit\'es de cette \'equation 
sont $z_{0}$ et $\infty$ et elles sont fuchsiennes. La m\'ethode de Frobenius 
en fournit des syst\`emes fondamentaux de solutions en $0$ et en $\infty$: 
$$
\begin{cases}
1 \\
z_{0}\log(1-z/z_{0})
\end{cases}
\quad
\begin{cases}
1 \\
-\log(1 - w/w_{0}) + \log(-w)
\end{cases}
$$
On a pos\'e $w = 1/z$ et $w_{0} = 1/ z_{0}$.\\

On les d\'eformera \`a l'aide de la fonction 
$l_{q}^{+}(z) = z (\Theta_{q}^{+})'(z)/\Theta_{q}^{+}(z)$, o\`u l'on a pos\'e
$\Theta_{q}^{+}(z) = (z;p)_{\infty} = 
\underset{0 \leq k}{\prod}(1 - z q^{-k})$. La fonction $l_{q}^{+}(z)$
v\'erifie
l'\'equation aux $q$-diff\'erences fuchsienne:
$$ 
\delta_{q}^{2}f = \frac{1}{1 - q^{2}\frac{z}{z_{0}}} \delta_{q} f
\qquad \delta_{q} = \frac{\sigma_{q} - Id}{q - 1}
$$

Les syst\`emes fondamentaux canoniques \emph{renormalis\'es} de solutions 
en $0$ et en $\infty$ sont
$$
\begin{cases}
1 \\
(q-1)z_{0} l_{q}^{+}(z/z_{0})
\end{cases}
\quad
\begin{cases}
1 \\
(q-1)(l_{q}(w) - l_{q}^{+}(w/qw_{0}))
\end{cases}
$$
La matrice de connexion est donc
$$
P(z) =
\begin{pmatrix}
1 & (q-1)z_{0}(1 - l_{q}(z) + l_{q}(z/z_{0}))\\
0 & -z_{0}
\end{pmatrix}
$$

Pour calculer sa limite lorsque $\epsilon \to 0^{+}$, on invoque les 
estimations suivantes, ais\'ement obtenues de fa\c{c}on analogue \`a celles de 
3.1:
$$
\lim_{\epsilon \to 0^{+}} \Theta_{q}^{+}(z_{2})/\Theta_{q}^{+}(z_{1}) =
(1 - z)^{\alpha_{2} - \alpha_{1}}
$$
$$
\lim_{\epsilon \to 0^{+}} (q-1) l^{+}_{q}(z) = \log(1-z_{0})
$$
On trouve alors
$$
\tilde{P}(z) =
\begin{pmatrix}
1 & z_{0}(\log(-z/z_{0}) - \log(-z))\\
0 & -z_{0}
\end{pmatrix}
$$
Celle-ci est localement constante sur la sph\`ere de Riemann priv\'ee de la
coupure triangulaire
$[0;\infty] \cup [0;z_{0}] \cup [z_{0};\infty]$. Ses deux 
d\'eterminations
$\tilde{P}_{0}$ et $\tilde{P}_{1}$ satisfont
$$
(\tilde{P}_{1})^{-1} \tilde{P}_{0} =
\begin{pmatrix}
1 & z_{0}(a_{0} - a_{1})\\
0 & 1
\end{pmatrix}
$$
o\`u $a_{1}$ et $a_{0}$ sont les valeurs de la fonction 
$(\log(-z/z_{0}) - \log(-z))$, qui est localement constante
$(\log(-z/z_{0}) - \log(-z))$ sur les ouverts correspondants: donc
$a_{0} - a_{1} = 2 \imath \pi$ et l'on obtient la matrice de la monodromie
autour de $z_{0}$:
$$
(\tilde{P}_{1})^{-1} \tilde{P}_{0} =
\begin{pmatrix}
1 & 2 \imath \pi z_{0}\\
0 & 1
\end{pmatrix}
$$


\subsection{Un exemple incomplet non fuchsien}

Voici maintenant un exemple \emph{incomplet} de d\'eformation d'une \'equation
diff\'erentielle \emph{non fuchsienne}, ceci, \`a titre exp\'erimental: on
constate qu'elle peut \^etre d\'eform\'ee en une famille d'\'equations aux
$q$-diff\'erences \emph{fuchsiennes}, et m\^eme que la th\'eorie qui 
pr\'ec\`ede fonctionne encore dans ce cas. Cet exemple a \'et\'e d\'ecouvert
en collaboration avec Zhang Changgui.\\

On part de l'\'equation scalaire:

$$\delta \tilde{f} = \tilde{b} \tilde{f}$$

avec $\tilde{b}(z) = e^{\frac{1}{z+1}}$: la ``matrice'' de cette \'equation a 
donc une singularit\'e essentielle en $-1$, c'est donc une \'equation
\emph{tr\`es, tr\`es} irr\'eguli\`ere\footnote{
Il e\^ut suffi d'un p\^ole double pour la rendre irr\'eguli\`ere, alors une
singularit\'e essentielle, pensez donc!}.\\

On la d\'eforme en prenant 
$b_{\epsilon} (z) = 
\underset{0 \leq k \leq N}{\sum} \frac{1}{k!} (\frac{1}{z + 1})^{k}$, avec
$\epsilon = \frac{1}{N}$: ce n'est donc pas une famille continue, mais il
est clair que ce n'est pas le point essentiel. En tout cas, on a bien la
convergence de $b_{\epsilon}$ vers $\tilde{b}$ soumise \`a toutes les 
contraintes du chapitre 4.\\

Pour voir comment la th\'eorie fonctionne ici, on devra proc\'eder de fa\c{c}on
indirecte: on ne connait pas les z\'eros de $a_{\epsilon}$, on ne saura donc
pas \'ecrire les solutions en termes de fonctions $\Theta$, ni non plus la
matrice de connexion. Mais on a tout de m\^eme affaire \`a des fonctions 
m\'eromorphes, auxquelles nos calculs g\'en\'eraux s'appliquent tr\`es bien !
Pour se simplifier la vie, on utilisera des caract\`eres $e_{q,c}$ qui
confluent en $z^{\gamma}$.\\

Les exposants sont:

\begin{itemize}

\item{En $0$, $c_{\epsilon} = a_{\epsilon}(0) = 1 + (q-1) b_{\epsilon}(0)$, 
donc de la forme $1 + (q-1)e + o(q-1)$. Le caract\`ere correspondant va 
tendre vers $z^{e}$: cela correspond bien \`a l'exposant de l'\'equation
diff\'erentielle en $0$.}

\item{En $\infty$, $a_{\epsilon}(\infty) = 1 + (q-1) b_{\epsilon}(\infty) = q$
et la solution correspondante sera de la forme $z \times$ (fonction
valant $1$ en $\infty$). Le caract\`ere correspondant tendra vers $z$ : pas
de monodromie de ce c\^ot\'e. Cela correspond bien \`a l'exposant de
l'\'equation diff\'erentielle en $\infty$.}

\end{itemize}

\subsubsection{Solutions de l'\'equation aux $q$-diff\'erences}

\begin{enumerate}

\item{Solution locale en $0$:\\

On \'ecrit 
$$
f_{\epsilon}^{(0)} = e_{q,c_{\epsilon}} g_{\epsilon}^{(0)}
$$
avec 
$\sigma_{q} g_{\epsilon}^{(0)} = 
\frac{a_{\epsilon}}{c_{\epsilon}} g_{\epsilon}^{(0)}$.
Donc 
$$
g_{\epsilon}^{(0)}(z) = 
\underset{r \geq 1}{\prod} \frac{a_{\epsilon}(q^{-r}z)}{c_{\epsilon}}
$$
}

\item{Solution locale en $\infty$:\\

Les calculs \'etant relativement simples, on les fait en $z$. On \'ecrit
$$
f_{\epsilon}^{(\infty)} = z g_{\epsilon}^{(\infty)}
$$
avec 
$\sigma_{q} g_{\epsilon}^{(\infty)} = 
\frac{a_{\epsilon}}{q} g_{\epsilon}^{(\infty)}$.
Donc 
$$
g_{\epsilon}^{(\infty)}(z) = 
\underset{r \geq 0}{\prod} \frac{q}{a_{\epsilon}(q^{r}z)}
$$
}

\item{Connexion des solutions locales:\\

On obtient, par division:
$$
p_{\epsilon}(z) = 
\frac{e_{q,c_{\epsilon}}}{z}
\underset{r \geq 1}{\prod} \frac{a_{\epsilon}(q^{-r}z)}{c_{\epsilon}}
\underset{r \geq 0}{\prod} \frac{a_{\epsilon}(q^{r}z)}{q}
$$
C'est la partie form\'ee des deux produits que l'on pourrait exprimer \`a
l'aide des fonctions $\Theta$ si l'on connaissait les z\'eros de 
$a_{\epsilon}$. Mais le nombre de facteurs $\Theta$ serait bien sur variable!
}

\end{enumerate}

\subsubsection{Solutions de l'\'equation diff\'erentielle}

\begin{enumerate}

\item{Solution locale en $0$:\\

On \'ecrit 
$$
\tilde{f}^{(0)} = z^{e} \tilde{g}^{(0)}
$$ 
avec
$\delta \tilde{g}^{(0)} = (\tilde{b} - e) \tilde{g}^{(0)}$. Donc
$$
\tilde{g}^{(0)}(z) = e^{\int_{0}^{z} \frac{\tilde{b}(t) - e}{t} \; dt}
$$
}

\item{Solution locale en $\infty$:\\

On \'ecrit 
$$
\tilde{f}^{(\infty)} = z \tilde{g}^{(\infty)}
$$ 
avec
$\delta \tilde{g}^{(\infty)} = (\tilde{b} - 1) \tilde{g}^{(\infty)}$. Donc
$$
\tilde{g}^{(\infty)}(z) = 
e^{\int_{0}^{1/z} \frac{1 - \tilde{b}(1/t)}{t} \; dt}
$$
}

\item{Connexion des solutions locales:\\

On obtient, par division:
$$
\tilde{p}(z) = 
\frac{z^{e}}{z}
e^{\int_{0}^{z} \frac{\tilde{b}(t) - e}{t} \; dt - \int_{0}^{1/z} \frac{1 - \tilde{b}(1/t)}{t} \; dt}
$$
}

\end{enumerate}

\subsubsection{Monodromie de l'\'equation diff\'erentielle}

Les seules singularit\'es sont $0$, $-1$ et $\infty$. La monodromie en $0$
est (la multiplication par) $e^{2 \imath \pi e}$ et il n'y en a pas en 
$\infty$ (singularit\'es r\'eguli\`eres d'exposants $e$ et $1$). 
Donc la monodromie en $-1$ est $e^{ - 2 \imath \pi e}$: on connait ici le
groupe fondamental !

\subsubsection{Confluence}

Le facteur ``log-car'' $\frac{e_{q,c_{\epsilon}}}{z}$ de la matrice de
connexion $p_{\epsilon}$ tend vers le facteur $\frac{z^{e}}{z}$ de $\tilde{p}$.
Le facteur interne 
$\underset{r \geq 1}{\prod} \frac{a_{\epsilon}(q^{-r}z)}{c_{\epsilon}}
\underset{r \geq 0}{\prod} \frac{a_{\epsilon}(q^{r}z)}{q}$
tend vers
$e^{\int_{0}^{z} \frac{\tilde{b}(t) - e}{t} \; dt - \int_{0}^{1/z} \frac{1 - \tilde{b}(1/t)}{t} \; dt}$
par la th\'eorie sur les \'equations avec conditions initiales fixes (3.2)
qui s'applique sans probl\`eme ici. Il y a peut-\^etre une preuve directe
via la th\'eorie des int\'egrales multiplicatives.\\

La limite de la matrice de connexion est donc la fonction localement constante
$\tilde{p}$ sur $\mathbf{S} - \overline{R}$, avec une valeur sur
$\mathcal{H}$ et une valeur sur $- \mathcal{H}$, reli\'ees par le prix de
franchissement de coupure: $e^{ - 2 \imath \pi e}$.


\tableofcontents

\end{document}